\documentclass[a4paper,12pt]{elsarticle}

\usepackage{latexsym}

\usepackage{setspace}
\usepackage{amsmath}
\usepackage{amssymb}
\usepackage{subfigure}
\usepackage{latexsym}
\usepackage{amsmath}
\usepackage{amssymb}
\usepackage{bigints}
\usepackage{natbib}
\usepackage{wrapfig,lipsum,booktabs}
\usepackage{caption}
\usepackage{graphicx}
\usepackage{subcaption}
\usepackage{placeins}
\usepackage{pgfplots} 
\usepackage{mwe} 
\usepackage{epstopdf}
\usepackage{color}
\setcitestyle{authoryear,open={(},close={)}}
\usepackage{hyperref}
\hypersetup{colorlinks,citecolor=blue}
\usepackage{graphics}
\usepackage{multirow}
\usepackage{pdflscape}
\setcitestyle{numbers}
\setcitestyle{square}
\makeatletter
\def\ps@pprintTitle{%
  \let\@oddhead\@empty
  \let\@evenhead\@empty
  \let\@oddfoot\@empty
  \let\@evenfoot\@oddfoot
}
\makeatother
\makeatletter \oddsidemargin  -.5cm \evensidemargin -.5cm \textwidth
17.5cm \topmargin -.65in \textheight 24cm

\newtheorem{Lemma}{Lemma}

\newtheorem{theorem}{Theorem}[section]

\newtheorem{lemma}[theorem]{Lemma}

\newtheorem{corollary}[theorem]{Corollary}
\newtheorem{remark}{Remark}[section]
\newtheorem{example}{Example}[section]

\numberwithin{equation}{section}
\numberwithin{equation}{section}

\title{On the improved estimation of ordered scale parameters based on doubly type-II censored sample}
\author{Shrajal Bajpai and Lakshmi Kanta Patra\footnote
	{\baselineskip=10pt
		~lkpatra@iitbhilai.ac.in; patralakshmi@gmail.com}    \\
\it  Department of Mathematics\\
 \it Indian Institute of Technology Bhilai, India, 491001}
\spacing{1.1}
\begin{document}
	\begin{frontmatter}
\date{}
\begin{abstract}
	A doubly type-II censored scheme is an important sampling scheme in the life testing experiment and reliability engineering. In the present commutation, we have considered estimating ordered scale parameters of two exponential distributions based on doubly type-II censored samples with respect to a general scale invariant loss function. We have obtained several estimators that improve upon the BAEE. We also propose a class of improved estimators. It is shown that the boundary estimator of this class is generalized Bayes. As an application, we have derived improved estimators with respect to three special loss functions, namely quadratic loss, entropy loss, and symmetric loss function. We have applied these results to special life-testing sampling schemes. Finally, we conducted a simulation study to compare the performance of the improved estimators. A real-life data analysis has been considered for implementation purposes.  
  \\
\noindent\textbf{Keywords}:  Doubly type-II censored sample, Best affine equivariant estimator, Scale invariant loss function, 
Generalized Bayes estimator, Risk function, Improved estimator. 
\end{abstract}
\end{frontmatter}
\section{Introduction}
In practical life testing scenarios, observing the lifetimes of all tested items is often difficult due to constraints such as time, cost, and resource limitations during data collection. Censoring plays a vital role in such situations. A censored sample contains observations, but only some pertinent details are available. The underlying concept of censoring involves misreporting data either below or above a certain threshold. 
Several censoring schemes exist in the literature, one of the most popular censoring methodologies is a doubly type-II censoring scheme. Suppose in a life testing experiment,  $X_1<X_2<\dots<X_n$ are the observed $n$ ordered lifetimes of  experimental units. If the first $r$ lifetimes and the last $s$ lifetimes are censored then we  get a doubly type-II censored sample as   $X_{r} \le X_{r+1} \le \dots \le X_{n-s}$. For a detailed discussion on censoring, we refer to \cite{meeker2014statistical}, \cite{balakrishnan2014art}. 
A two-parameter exponential distribution is used to analyse lifetime data in life testing and reliability theory, for example \cite{lawless1977prediction} applied the two-parameter exponential distribution to analyse lifetime data. Numerous authors have explored the problem of estimating unknown parameters of various statistical models based on censored samples in the past. \cite{elfessi1997estimation} used a doubly censored sample for estimating the location parameter, rate parameter, and a linear function of location and scale parameter of a two-parameter exponential distribution. The author used \cite{stein1964} techniques to derive an estimator better than the best affine equivariant estimator (BAEE).   \cite{madi2002invariant} studied invariant estimation of the scale parameter of a shifted exponential distribution with respect to a scale-invariant loss function. He obtained a dominating estimator using the techniques of \cite{brewster1974improving}. Estimating scale parameter of a shifted exponential distribution based on record statistics and progressive type censored sample has been studied by \cite{madi2006decision} and \cite{madi2010note} respectively. In these works, he used the sequential technique of \cite{brewster1974improving} to propose a smooth, improved estimator. \cite{tripathi2017estimating} considered the estimation of linear parametric function of a two-parameter exponential distribution based on doubly type II censored sample with respect to arbitrary strictly convex loss function. They have proved the inadmissibility of BAEE by proposing improved estimators. The authors used the IERD approach of \cite{kubokawa1994unified} to obtain a class of improved estimators. Some other works based on censored sample are  \cite{madi2008improved}, \cite{fernandez2000bayesian}, \cite{zhu2015exact}, \cite{patra2020improved} and reference therein.
	
In the literature, we have seen that the estimation of parameters based on a censored sample has studied only one population. In this work, we have studied the estimation of ordered scale parameters of two exponential populations based on doubly type-II censored sample. We have used the techniques of \cite{brewster1974improving}, \cite{strawderman1974minimax} \cite{kubokawa1994unified}, \cite{maruyama1998minimax} to derive improved estimators. Several authors have used these techniques to derive improved estimators for the unknown parameters of various probability models. Some recent work on these direction are \cite{bobotas2011improved}, \cite{bobotas2017estimation}, \cite{bobotas2019improved}, \cite{patra2021componentwise}, \cite{patra2021minimax}, \cite{garg2023componentwise}.

Let $(X_{i1},\ldots, X_{in_{i}})$ be a random sample drawn from the
$i$-th population $\Pi_{i},$ $i=1,2$. The probability density function of the $i$-th population $\Pi_{i}$ is given by
\begin{eqnarray}\label{model}
	f_{i}(x;\mu_{i},\sigma_{i})=\left\{\begin{array}{ll}
		\frac{1}{\sigma_{i}}~\exp{\Big(-\frac{x-\mu_{i}}{\sigma_{i}}\Big)},
		& \textrm{if $x>\mu_{i},$}\\
		0,& \textrm{otherwise,}
	\end{array} \right.
\end{eqnarray}
where $-\infty<\mu_{i}<\infty$ and $\sigma_{i}>0.$ We assume that $\sigma_{i}$'s and $\mu_{i}$'s, $i=1,2$ are unknown.  It is known in advance that $\sigma_1\le \sigma_2$.  A doubly type-II censored sample from $i$-th sample is 
\begin{eqnarray*}
	X_{a_{i}} \le X_{a_{i}+1} \le \dots \le X_{b_{i}},~~1 \le a_{i} \le b_{i} \le n_{i}.
\end{eqnarray*}
From \cite{elfessi1997estimation} we have a complete sufficient statistics based on the doubly type-II censored sample is  $(X_{a_{1}}, X_{a_{2}}, V_1,V_2)$ are independent with densities
\begin{eqnarray*}
	h_{X_{a_i}}(x)=A_i\left[1-e^{-\frac{x-\mu_i}{\sigma_i}}\right]^{a_i-1}\frac{1}{\sigma_i}
	\left[e^{-\frac{x-\mu_i}{\sigma}}\right]^{n_i-a_i+1}~x>\mu_i,~\sigma_i>0
\end{eqnarray*}
with $A_i=\frac{n_{i}!}{(a_{i}-1)!(n_{i}-a_{i})!}$ and
\begin{eqnarray*}
	g_{V_i}(v)=\frac{1}{\sigma_i^{b_i-a_i}\Gamma(b_i-a_i)}e^{-v/\sigma_i}v^{b_i-a_i-1},~v>0,
	\sigma_i>0,
\end{eqnarray*}
$i=1,2$. From now onwards  we denote $\underline{X_{a}}=(X_{a_1},X_{a_2})$, $\underline{V}=(V_1,V_2)$,  $b=b_2+b_1,a=a_2+a_1$ and $\theta=(\mu_1,\mu_2,\sigma_1,\sigma_2)$. In this manuscript, we study component wise estimation of ordered scale parameters of two exponential distributions based on doubly type-II censored sample from decision theoretic point of view.  For this estimation problem we consider a class of bowl-shaped loss function $L(\delta/\theta)$, where $\delta$ is an estimator of $\theta$. We assume that the loss function $L(t)$ satisfies the following conditions:  (i) $L(t)$ is a bowl-shaped function with the minimum value 0 attained at $t = 1.$ So $L(t)$ is decreasing when $t \le 1$ and increasing for $t \ge1.$  (ii) Integrals involving $L(t)$ are finite and differentiation under the integral sign is permissible. 

We invoke the  principle of invariance to find the BAEE.   For this we  consider the affine group of transformations $G_{p_i,q_i}=\{g_{p_i,q_i}(x)=p_ix_{ij}+q_{i}, j=1,\dots,n_i,  i=1,2\}$. Under this group, we have 
$$(X_{a_i},V_i) \longrightarrow (p_iX_{a_i}+q_i,p_iV_i)~~~\mbox{ and } (\mu_{i},\sigma_i) \longrightarrow (p_i\mu_i+q_i,p_i\sigma_i)$$
After some simplification the form of an affine equivariant estimators of $\sigma_i$ can be obtained as 
\[\delta_{ic}(\underline{X}_{a},\underline{V})=cV_{i},\] where $c>0$ is constant. The following theorem provides the BAEE of $\sigma_i$. The proof is simple and hence omitted for the shake of brevity.
\begin{lemma}
	Under the scale invariant loss function given as L(t) the (unrestricted) best affine equivariant estimator for $\sigma_i$ is $\delta_{0i}(\underline{X}_{a},\underline{V})=c_{0i}V_i$, where $c_{0i}$ is the unique solution of the equation 
	\begin{equation} \label{baee}
	E_{\sigma_i=1}(L^{'}(c_{0i}V_{i})V_{i})=0, ~i=1,2. 
\end{equation}
	\end{lemma}
	\begin{example}\rm
\begin{enumerate}
 \item[(i)]Consider the quadratic loss $L_1(t)=(t-1)^2$. The BAEE of $\sigma_i$ is $\delta_{0i} ^{1}=\frac{V_1}{b_{i}-a_{i}+1}$. 
		\item[(ii)] Consider the entropy loss function $L_2(t)=t-\ln t-1$. The BAEE  of $\sigma_i$ is 
		$\delta_{0i}^2=\frac{V_1}{b_{i}-a_{i}}$. 
		\item[(iii)]Consider the symmetric loss function $L_3(t)=t+\frac{1}{t}-2$. The BAEE for $\sigma_i$ is $\delta_{01}^{3}=\frac{V_i}{\sqrt{(b_i-a_i)(b_i-a_i-1)}}$.
\end{enumerate}
\end{example}
We aim to propose various estimators that improve upon the BAEE. The following lemma will be helpful in proving improvement results. 
\begin{lemma}[\cite{tripathi2017estimating}]\label{tripathi}
	Suppose that $Y$ is random variable with some distribution $g(y)$ and support on the  distribution is on $[0,\infty)$. Let $\psi(y)$ be a function that changes its sign once from negative to positive on the real line such that $\psi(y) < 0$ for $y < y_0$ and $\psi(y)>0$ for $ y > y_0$. If $h(y)$ be an increasing function with non negative value such that  $\int_{0}^{\infty} \psi(y) g(y) h(y)dy=0$ then $\int_{0}^{\infty} \psi(y)g(y) dy \le 0 .$
\end{lemma} 
The rest of the paper is organized as follows. In Section \ref{sect2}, we have proposed several estimators that improve upon the BAEE with respect to a general scale invariant loss function. As an application, we have derived the improved estimators for three special loss functions in Subsection \ref{subsect2.1}. A class of improved estimators has been obtained using the IERD approach of \cite{kubokawa1994unified} in Subsection \ref{subsec2.3}. In addition, we have shown that the boundary estimator of this class is a generalized Bayes estimator. Also, we have obtained the Maruyama type and Strawderman type improved estimators for quadratic loss and entropy loss function. Similarly, in section \ref{sect3}, we have discussed the improved estimation of $\sigma_2$. Section \ref{sect4} discusses the same estimation problem for special life testing sampling schemes. Finally, in Section \ref{sec5}, we have conducted a simulation study and data analysis. 
\section{Improved estimation of $\sigma_1$} \label{sect2}
\noindent In this section, we propose estimators that dominate the BAEE of $\sigma_1$. To derive a dominating estimator, we consider a class of estimators of the form
\begin{equation}\label{imp1}
\mathcal{C}_1=\left\{\delta_{\phi_1}=V_1\phi_{1}(Z_1): Z_1=\frac{V_2}{V_1}, \mbox{ and } \phi_{1} \mbox{ is a positive measurable function }\right\}
\end{equation}
For finding  improved estimator, in this section we will use the techniques of \cite{brewster1974improving}. The risk function of the estimator $\delta_{\phi_{1}}(\underline{X_a},\underline{V})$ is $$R\left(\theta,\delta_{\phi_{1}}\right)=E\left[E\left(L\left(\frac{V_1\phi_{1}(Z_1)}{\sigma_1}\big| Z_1=z_1\right)\right)\right].$$ 
The inner conditional risk is  $R_{1}(\theta,c)=E\left(L\left(\frac{V_1c}{\sigma_1}\right)\big| Z_1=z_1\right)=E_{\eta}(L(T_1c)|Z_1=z_1)$ with $T_1|Z_1=z_1 \sim Gamma(b-a,(1+\eta z_1)^{-1})$, where $\eta =\frac{\sigma_1}{\sigma_2}\le1$. Suppose $\zeta_\eta(t_1|z_1)$ is the conditional density of $T_1$ given $Z_1=z_1$. We can easily get $\frac{\zeta_{\eta}(t_1|z_1)}{\zeta_{1}(t_1|z_1)}$ is non decreasing in $t_1$. Now, using Lemma 3.4.2 in \cite{lehmann2006testing}, we get $E_{\eta}(L'(T_1c(1,z_1))T_1\big|Z_1=z_1)\geq E_{1}(L'(T_1c(1,z_1))T_1\big|Z_1=z_1)=0$. Consequently we obtain $c(\eta,z_1)\le c(1,z_1)$ with  $c(1,z_1)$ is unique solution of  the equation $E_{1}(L'(T_1c(1,z_1))T_1)=0$. Using the transformation  $z=t_1(1+z_1)$ we get
$EL'\left(\frac{Z c(1,z_1)}{1+z_1}\right)=0$ and comparing with (\ref{s1st1}) we have $c(1,z_1)= \beta(1+z_1)$. Consider a function $\phi_1^*(z_1)=\min\{\phi_1(z_1),c(1,z_1)\}$.  So we have $c(\eta,z_1) \le c(1,z_1) = \phi_1^*(z_1) <\phi_1(z_1)$ with positive probability provided $P_{\eta}(\phi_1(Z_1)>c(1,Z_1)) \ne 0$. It follows that $R_1(\theta,\delta_{\phi^*_{1}})< R_1( \theta,\delta_{\phi_1})$.  This shows that $R(\theta,\delta_{\phi^*_{1}})< R( \theta,\delta_{\phi_1})$. Hence we state the following result. 
\begin{theorem} \label{Th1sigma1}
Let $\beta$ be the unique solution of 
\begin{equation}\label{s1st1}
	E(L'(\beta 
	Z))=0
\end{equation}
where Z follows $Gamma(b-a+1,1)$ distribution. Define a function $\phi_1^*(z_1)=\min\{\phi_1(z_1),c(1,z_1)\}$ with $c(1,z_1)= \beta(1+z_1)$. Then the estimator $\delta_{\phi_1^*}$ has uniformly smaller risk then $\delta_{\phi_1}$ under a scale  invariant loss function $L(t)$, provided $P(\phi_1(Z_1)>c(1,Z_1)) \ne 0$ for some $\eta>0$. 
\end{theorem}
\begin{corollary}
	The risk of the estimator  $\delta_{1S1}=\min\{c_{01}, \beta(1+z_1)\}V_1,$
 is nowhere larger than $\delta_{01}$ under a bowl shaped loss function $L(t)$ provided $\beta\le c_{01}.$
\end{corollary}
\noindent\textbf{Proof:}  Taking $\phi_1(z_1)=c_{01}$ we have 
 $c(\eta, z_1) < \phi_{01}(z_1) < c_{01} $ on a set of positive probability as we have $\beta\le c_{01}.$ Now $R_{1}(\theta,c)$ is strictly bowl shaped. So $R_1(\theta,c)$ is strictly increasing in $\left(c(\eta,z_1),c\right)$. Hence $R_1(\theta,\delta_{\phi_{1S1}})\le R_1( \theta,\delta_{01})$ and strict inequality for some $z_1>0$ . Hence we get the result. \hfill $\blacksquare$\\

Consider the another class of estimators
 \begin{equation}
 	\mathcal{C}_2=\left\{\delta_{\phi_2}=\phi_2(Z_1,Z_2)V_1: Z_1=\frac{V_2}{V_1},\mbox{ and } Z_2=\frac{X_{a_1}}{V_1}\right\}
 \end{equation}
\begin{theorem}\label{th2sigma1}
Suppose $Y\sim Gamma(b-a+2,1)$ distribution and let $\beta_1$ is unique solution of 
	\begin{equation}\label{s2t21}
			EL'( Y \beta_1)=0 .
	\end{equation}
Define a function $\phi_2^{*}=\min\{\phi_{2}(z_1,z_2),\beta_1(1+ z_1+(n_1-a_1+1)z_2)\}$ for $z_1>0$ and $z_2>0$. Then the estimator 
\begin{eqnarray*}
	\delta_{\phi^*_{2}}=\left\{
	\begin{array}{ll}
		\phi_2^{*}(Z_1,Z_2)V_1,~Z_1>0,~Z_2>0\\\\
		\phi_{2}(Z_1,Z_2)V_1,~~~~~~~~~~~~~~~~~~~\mbox{ otherwise. }
	\end{array}
	\right.
\end{eqnarray*}
dominates the estimator $\delta_{\phi_2}$ with respect to a scale invariant loss function $L(t)$ provided $P(\phi_{2}(Z_1,Z_2)>\beta_1(1+ Z_1+(n_1-a_1+1)Z_2)) \ne 0$ for $Z_2>0$. 
\end{theorem}
\textbf{Proof:}
Risk function for this estimator is 
\begin{align*}
	R(\theta,\delta_{\phi_2})=E\left[E \left(L\left(\frac{V_1\phi_{1}(Z_1,Z_2)}{\sigma_1}\right)\big| Z_1=z_1,Z_2=z_2 \right)\right]
\end{align*}
Consider $z_2>0$, the inner conditional risk is obtained as  $R_2(\theta, c)=E\left(L(T_1c)\big| Z_1=z_1,Z_2=z_2\right)$, the expectation has been taken with respect to  $T_1|Z_1=z_1,Z_2=z_2 \sim f_{\rho,\eta}(t_1)$, where 
\begin{equation*}
	f_{\rho,\eta}(t_1)\propto e^{-t_1(1+\eta z_1+(n_1-a_1+1) z_2)}t_{1}^{b-a}\eta^{b_2-a_2}[1-e^{-t_{1}z_2}e^\rho]^{a_1-1}e^{\rho(n_1-a_1+1)},   z_1 >0, t> \max\left\{0,\frac{\rho}{z_2}\right\} ,
\end{equation*}
where $\rho=\frac{\mu_1}{\sigma_1}$. Let 
 $c_{\rho,\eta}(z_1,z_2)$ be the unique minimize of $R_2(\theta, c)$, so $c_{\rho,\eta}(z_1,z_2)$ is unique solution of 
\begin{equation}\label{st2th2}
	E_{\rho,\eta}(L'(T_1c_{\rho,\eta}(z_1,z_2))T_1\big|Z_1=z_1,Z_2=z_2)=0.
\end{equation}
Now, $\rho>0$, if $c^{-1}_{\rho,\eta}(z_1,z_2)<\frac{\rho}{z_2}$ implies $t_1c_{\rho,\eta}(z_1,z_2)>1$. Consequently we have
\begin{equation}
		E_{\rho,\eta}(L'(T_1c_{\rho,\eta}(z_1,z_2))T_1\big|Z_1=z_1,Z_2=z_2)>0.
\end{equation}
which is a contradiction to (\ref{st2th2}). Therefore this confirms $c^{-1}_{\rho,\eta}(z_1,z_2)>\frac{\rho}{z_2}$. Again, $c_{0,\eta}(z_1,z_2)$ uniquely satisfies
\begin{equation}\label{st2th3}
E_{0,\eta}(L'(T_1c_{0,\eta}(z_1,z_2))T_1\big|Z_1=z_1,Z_2=z_2)=0.
\end{equation}
Now we have,
\begin{align}\nonumber\label{st2th4}
	E_{0,\eta}(L'(T_1c_{\rho,\eta}(z_1,z_2))T_1\big|Z_1=z_1,Z_2=z_2) 
	= &\int_{0}^{\infty}L'(t_1 c_{\rho,\eta}(z_1,z_2))t_1 f_{}(t_1|z_1,z_2)dt_1 \\
	=&\int_{0}^{\frac{\rho}{z_2}}L'(t_1 c_{\rho,\eta}(z_1,z_2))t_1 f_{}(t_1|z_1,z_2)dt_1<0.
\end{align}
Therefore, from (\ref{st2th3}) and (\ref{st2th4}), we have $c_{\rho,\eta}(z_1,z_2)<c_{0,\eta}(z_1,z_2)$
Again, It can be easily verified that $\frac{f_{0,\eta}(t_1|z_1,z_2)}{f_{0,1}(t_1|z_1,z_2)}$ is non-decreasing in $t_1$. This shows that $c_{\rho,\eta}(z_1,z_2)<c_{0,\eta}(z_1,z_2)<c_{0,1}(z_1,z_2)$.
For $\rho\le 0$ we have $c_{\rho,\eta}(z_1,z_1)=c_{0,\eta}(z_1,z_2)<c_{0,1}(z_1,z_2)$ and  $c_{0,1}(z_1,z_2)$ is uniquely satisfies $E_{0,1}(L'(T_1c_{0,1}(z_1,z_2))T_1\big|Z_1=z_1,Z_2=z_2)=0$
that is, 
\begin{equation*}
	\int_{0}^{\infty}L'(t_1c_{0,1}(z_1,z_2))t_1^{b-a+1}e^{-t_1(1+\eta z_1+(n_1-a_1+1)z_2}[1-e^{-t_1z_2}]^{a_1-1}dt_1=0.
\end{equation*}
Substituting $y=t_1(1+ z_1+(n_1-a_1+1)z_2)$ in the above integral we obtain 
\begin{equation*}
	\int_{0}^{\infty}L'\left( \frac{yc_{0,1}(z_1,z_2)}{1+z_1+(n_1-a_1+1)z_2}\right) y^{b-a+1}e^{-y}[1-e^{\frac{yz_2}{1+z_1+(n_1-a_1+1)z_2}}]^{a_1-1} =0 .
\end{equation*}
Now using Lemma \ref{tripathi} we get
\begin{equation}
	\int_{0}^{\infty}L'\left( \frac{yc_{0,1}(z_1,z_2)}{1+z_1+(n_1-a_1+1)z_2}\right) y^{b-a+1}e^{-y}\le 0 .
\end{equation}
Comparing with (\ref{s2t21}) we have 
\begin{equation*}
	c_{0,1}(z_1,z_2) \le \beta_1(1+z_1+(n_1-a_1+1)z_2)
\end{equation*}
Now consider a function $\phi^*_{2}=\min\left\{\phi_2(z_1,z_2), \beta_1(1+z_1+(n_1-a_1+1)z_2)\right\}$ then $c_{\rho,\eta}(z_1, z_2) <  \beta_1(1+z_1+(n_1-a_1+1)z_2)=\phi^*_{2}<c_{01}$. Since $R_2(\theta,c)$ is strictly increasing in $\left (c_{\rho,\eta}(z_1,z_2), \infty\right)$. Hence we have $R(\theta,\delta_{\phi^*_{2}}) < R(\theta,\delta_{01})$. This proves the result.  \hfill $\blacksquare$
\begin{corollary}
The estimator 
\begin{eqnarray*}
	\delta_{1S2}=\left\{
	\begin{array}{ll}
		\min\{c_{01},\beta_1(1+ z_1+(n_1-a_1+1)z_2)\}V_1,~z_1>0,~z_2>0\\\\
		c_{01}V_1,~~~~~~~~~~~~~~~~~~~\mbox{ otherwise. }
	\end{array}
	\right.
\end{eqnarray*}
dominates $\delta_{01}$ under a scale invariant loss function $L(t)$ provided provided $\beta_1 <c_{01}.$
\end{corollary}
Consider the another class of estimators
\begin{equation*}
	\mathcal{C}_3=\left\{\delta_{\phi_3}=\phi_{3}(Z_1,Z_2,Z_3)V_1: ~Z_1=\frac{V_2}{V_1},~ Z_2=\frac{X_{a_1}}{V_1}, ~ Z_3=\frac{X_{a_2}}{V_1}\right\}
\end{equation*}
\begin{theorem}\label{th3sigma1}
	Suppose $U\sim Gamma(b-a+3,1)$ distribution and let $\beta_2$ is unique solution of 
	\begin{equation}\label{s2t31}
		EL'( U\beta_2)=0 .
	\end{equation}
	Define a function $\phi_3^{*}=\min\{\phi_{3}(z_1,z_2,z_3),\beta_2(1+ z_1+(n_1-a_1+1)z_2+(n_2-a_2+1)z_3)\}$ for $z_2>0$ and $z_3>0$. Then the estimator 
	\begin{eqnarray*}
		\delta_{\phi^*_{3}}=\left\{
		\begin{array}{ll}
			\phi_3^{*}(Z_1,Z_2,Z_3)V_1,~Z_i>0,~i=1,2,3\\\\
			\phi_{3}(Z_1,Z_2,Z_3)V_1,~~~~~~~~~\mbox{ otherwise. }
		\end{array}
		\right.
	\end{eqnarray*}
	dominates the estimator $\delta_{\phi_3}$ with respect to a scale invariant loss $L(t)$ function provided $P(\phi_{3}(Z_1,Z_2,Z_3)>\beta_2(1+ Z_1+(n_1-a_1+1)Z_2+(n_2-a_2+1)Z_3)) \ne 0$ for $Z_i>0,i=1,2,3$. 
\end{theorem}
\textbf{Proof:} Proof of this theorem is similar to the proof of Theorem \ref{th2sigma1}. 
\begin{corollary}
	The estimator 
	\begin{eqnarray*}\label{exp}
	\delta_{1S3}(z_1,z_2,z_3)=\left\{
	\begin{array}{ll}
		\min\{c_{01},\beta_2(1+ Z_1+(n_1-a_1+1)Z_2+(n_2-a_2+1)Z_3)\}V_1,~Z_i>0,~i=1,2,3\\\\
		c_{01}V_1,~~~~~~~~~~~~~~~~~~~\mbox{ otherwise. }
	\end{array}
		\right.
	\end{eqnarray*}
	dominates $\delta_{01}$ under a scale invariant loss function $L(t)$ provided provided $\beta_2 <c_{01}.$
\end{corollary}
\subsection{Improved estimators of $\sigma_1$ for special loss functions} \label{subsect2.1}
In this subsection, we will derive improved estimators for various special loss functions as an application Theorem \ref{Th1sigma1}, \ref{th2sigma1} and \ref{th3sigma1}. Improved estimators for quadratic loss function $L_1(t)$ are
\begin{eqnarray*}
	\delta_{1S1}=\left\{
	\begin{array}{ll}
		\frac{V_1}{b-a+1}(1+ Z_1),~ \text{if}\ \ 0< Z_1 < \frac{b-a}{b_1-a_1+1}-1\\
		\frac{1}{b_1-a_1+1}V_1,~~~~~~~~~~~~~~~~~~~\mbox{ otherwise. }
	\end{array}
	\right.
\end{eqnarray*}

\begin{eqnarray*}
\delta_{1S2}=\left\{
\begin{array}{ll}
	\min\left\{\frac{1}{b_1-a_1+1},\frac{(1+ Z_1+(n_1-a_1+1)Z_2)}{b-a+2}\right\}V_1,~Z_i>0,~i=1,2\\\\
	\frac{1}{b_1-a_1+1}V_1,~~~~~~~~~~~~~~~~~~~\mbox{ otherwise. }
\end{array}
\right.
\end{eqnarray*}

\begin{eqnarray*}
	\delta_{1S3}=\left\{
	\begin{array}{ll}
		\min\{c_{01},\beta_2(1+ Z_1+(n_1-a_1+1)Z_2+(n_2-a_2+1)Z_3)\}V_1,~Z_i>0,i=1,2,3\\\\
		c_{01}V_1,~~~~~~~~~~~~~~~~~~~\mbox{ otherwise. }
	\end{array}
	\right.
\end{eqnarray*}
where, $\beta_2=\frac{1}{b-a+3}$ respectively.

Using the Theorem \ref{Th1sigma1}, \ref{th2sigma1} and \ref{th3sigma1} we obtained the improved estimators for entropy loss function $L_2(t)$ are
\begin{eqnarray*}
	\delta_{1S1}=\left\{
	\begin{array}{ll}
		\frac{1}{b-a}(1+ Z_1)V_1, \text{if}\ \  0< Z_1<\frac{b-a}{b_1-a_1}-1\\
		\frac{1}{b_1-a_1}V_1,~~~~~~~~~~~~~~~~~~~\mbox{ otherwise. }
	\end{array}
	\right.
\end{eqnarray*} \\
 \begin{eqnarray*}
	\delta_{1S2}=\left\{
	\begin{array}{ll}
		\min\left\{\frac{(1+ Z_1+(n_1-a_1+1)Z_2)}{b-a+1},	\frac{1}{b_1-a_1}\right\}V_1,~Z_i>0,~i=1,2\\\\
		\frac{V_1}{b_1-a_1},~~~~~~~~~~~~~~~~~~~\mbox{ otherwise. }
	\end{array}
	\right.
\end{eqnarray*}\\
	\begin{eqnarray*}
	\delta_{1S3}=\left\{
	\begin{array}{ll}
		\min\{c_{01},\beta_2(1+ Z_1+(n_1-a_1+1)Z_2+(n_2-a_2+1)Z_3)\}V_1,~Z_i>0, i=1,2,3\\\\
		c_{01}V_1,~~~~~~~~~~~~~~~~~~~\mbox{ otherwise. }
	\end{array}
	\right.
\end{eqnarray*}
where, $\beta_2=\frac{1}{b-a+2}$ respectively.

For symmetric loss function $L_3(t)$ using the Theorem \ref{Th1sigma1}, \ref{th2sigma1} and \ref{th3sigma1} the improved estimators  are
\begin{eqnarray*}
	\delta_{1S1}=\left\{
	\begin{array}{ll}
		\left\{\frac{1+ Z_1}{\sqrt{b-a-1}\sqrt{b-a}}\right\}V_1, \text{if}\ 0<Z_1<\frac{\sqrt{(b-a-1)(b-a)}}{\sqrt{(b_1-a_1)(b_1-a_1-1)}}-1\\
		\frac{1}{\sqrt{(b_1-a_1)(b_1-a_1-1)}}V_1,~~~~~~~~~~~~~~~~~~~\mbox{ otherwise. }
	\end{array}
	\right.
\end{eqnarray*}

\begin{eqnarray*}
	\delta_{1S2}=\left\{
	\begin{array}{ll}
		\min\left\{\frac{(1+ Z_1+(n_1-a_1+1)Z_2)}{\sqrt{(b-a+1)(b-a)}},	\frac{1}{\sqrt{(b_1-a_1)(b_1-a_1-1)}}\right\}V_1,~Z_i>0,~i=1,2\\\\
		\frac{V_1}{\sqrt{(b_1-a1)(b_1-a_1-1)}},~~~~~~~~~~~~~~~~~~~\mbox{ otherwise. }
	\end{array}
	\right.
\end{eqnarray*}

\begin{eqnarray*}
	\delta_{1S3}=\left\{
	\begin{array}{ll}
		\min\{c_{01},\beta_2(1+ Z_1+(n_1-a_1+1)Z_2+(n_2-a_2+1)Z_3)\}V_1,~Z_i>0, i=1,2,3\\\\
		c_{01}V_1,~~~~~~~~~~~~~~~~~~~\mbox{ otherwise. }
	\end{array}
	\right.
\end{eqnarray*}
where, $\beta_2=\frac{1}{\sqrt{(b-a+2)(b-a+1)}}$.

	\begin{remark}
		From \cite{fernandez2002computing}  we have the unrestricted MLE of $\sigma_1$  $\delta_{mle}=\frac{V_1}{b_1-a_1+1}$. 
	From \cite{misra2002smooth}, 	the restricted maximum likelihood estimator for $\sigma_1$ with order restriction $\sigma_{1}\le \sigma_{2}$ can be obtain as 
	$$\delta_{Rmle}=
	\min \left\lbrace \frac{V_1}{b_1-a_1+1},\frac{V_1+V_2}{b-a+2}\right\rbrace=\min\left\lbrace\frac{1}{b_1-a_1+1},\frac{1+Z_1}{b-a+2}\right\rbrace V_1=\phi_{Rmle}(Z_1)V_1.$$
Using Theorem \ref{Th1sigma1}, we get the estimator $\delta_{\phi_1^r} = \phi_1^r(Z_1)V_1$, where $\phi_1^r(Z_1) = \min\{\phi_{Rmle}(Z_1), \beta(1 + Z_1)\}$, dominates $\delta_{Rmle}$ for estimating $\sigma_{1}$, with order restriction $\sigma_{1}\le\sigma_{2}$, under a general scale invariant loss function $L(t)$.
	\end{remark}

\subsection{Class of improved estimators of $\sigma_1$}\label{subsec2.3}
In the previous section we have proposed several improved estimators. Now we will derive a class of improved estimators. The Joint density of $T_1=\frac{V_1}{\sigma_1}$ and $Z_1$ is $$g_\eta(z_1,t_1)\propto\eta^{b_2-a_2}e^{-(\eta z_1 +1)t_{1}}z_{1}^{b_2-a_2-1}t_{1}^{b-a-1},   t_1>0,z_1>0,0<\eta\le1.$$ 
Define 
 $G_{\eta}(x,t_1)=\int_{0}^{x}g_{\eta}(u,t_1)du$ and $G_{1}(x,t_1)=\int_{0}^{x}g_{1}(u,t_1)du$.
Now we give sufficient conditions under which the estimator $\delta_{\phi_{1}}(\underline{X_{a}},\underline{V})$ dominate the BAEE.

\begin{theorem}
The risk of estimator $\delta_{\phi_{1}}(\underline{X_a},\underline{V})$ given in (\ref{imp1}) is nowhere larger than $\delta_{01}$ under a bowl shaped loss function $L(t)$ provided $\phi_{1}(u_1)$ satisfies the following conditions:
\begin{enumerate}
	\item[(i)]$\phi_1(u_1)$ is non-decreasing and $\lim_{u_1\rightarrow0}\phi_1(u_1) = c_{01},$
	\item[(ii)] $\int_{0}^{\infty}L'(\phi_{1}(u_1)t_1)t_1G_1(x,t_1)dt_1\ge0$.
\end{enumerate}
\end{theorem}
\textbf{Proof}: The proof is similar to the Theorem 2.8 of \cite{patra2021componentwise}. So we omit it. 

\noindent In the following corollaries we obtain class of improved estimators for special loss function. 
\begin{corollary}
	Consider the quadratic loss function $L_1(t)=(t-1)^2$. The risk of estimator $\delta_{\phi_{1}}(\underline{X_{a}},\underline{V})$ given in (\ref{imp1}) is nowhere larger than 
$\delta_{01}^1$ provided $\phi_{1}(u_1)$ satisfies the following conditions.
\begin{enumerate}
	\item [(i)] $\phi_{1}(u_1)$ is non decreasing and $\lim \limits_{u_1\rightarrow0}\phi_{1}(u_1)=\frac{1}{b_1-a_1+1}$
	\item [(ii)]$\phi_{1}(u_1)\ge\phi_{1}^{01}(u_1)$ where 
	$$\phi_{1}^{01}(u_1)=\frac{1}{b-a+1}\frac{\displaystyle \int_{0}^{1}y^{b_2-a_2-1}(yu_1+1)^{-(b-a+1)}dy}{\displaystyle\int_{0}^{1}y^{b_2-a_2-1}(yu_1+1)^{-(b-a+2)}dy}$$
\end{enumerate}
\end{corollary}
\begin{corollary}
Under the entropy loss function $L_2(t)=t-\ln t-1$, the risk of estimator $\delta_{\phi_{1}}(\underline{X_{a}},\underline{V})$ given in (\ref{imp1}) is nowhere larger than 
	$\delta_{01}^2$ provided $\phi_{1}(u_1)$ satisfies the following conditions.
	\begin{enumerate}
		\item [(i)] $\phi_{1}(u_1)$ is non decreasing and $\lim \limits_{u_1\rightarrow0}\phi_{1}(u_1)=\frac{1}{b_1-a_1}$
		\item [(ii)]$\phi_{1}(u_1)\ge\phi_{1}^{02}(u_1)$, where 
			$$\phi_{1}^{02}(u_1)=\frac{1}{b-a}\frac{\displaystyle \int_{0}^{1}y^{b_2-a_2-1}(yu_1+1)^{-(b-a)}dy}{\displaystyle\int_{0}^{1}y^{b_2-a_2-1}(yu_1+1)^{-(b-a+1)}dy}$$
	\end{enumerate}
\end{corollary}
\begin{corollary}
For symmetric loss function $L_3(t)=t+\frac{1}{t}-2$, the risk of estimator $\delta_{\phi_{1}}(\underline{X_{a}},\underline{V})$ given in (\ref{imp1}) is nowhere larger than 
$\delta_{01}^3$ provided $\phi_{1}(u_1)$ satisfies the following conditions.
\begin{enumerate}
	\item [(i)] $\phi_{1}(u_1)$ is non decreasing and $\lim \limits_{u_1\rightarrow0}\phi_{1}(u_1)=\frac{1}{\sqrt{(b_1-a_1)(b_1-a_1-1)}}$
	\item [(ii)]$\phi_{1}(u_1)\ge\phi_{1}^{03}(u_1)$ where 
	\begin{align*}
\phi_{1}^{03}(u_1)=\left( \frac{1}{(b-a-1)(b-a)}\frac{\displaystyle \int_{0}^{1}y^{b_2-a_2-1}(yu_1+1)^{-(b-a-1)}dy}{\displaystyle\int_{0}^{1}y^{b_2-a_2-1}(yu_1+1)^{-(b-a+1)}dy}\right)^{\frac{1}{2}}
\end{align*}
\end{enumerate}
\end{corollary}
\subsection{Maruyama type estimators}
Consider the following estimators for $\alpha>1$
\begin{equation*}\label{marL1}
	\delta_{\alpha,1}=\phi_{\alpha,1}(U_1)V_1~ \mbox{ with }~
	\phi_{\alpha,1}(u_1)=\frac{1}{b-a+1}\frac{\displaystyle \int_{0}^{1}y^{\alpha(b_2-a_2-1)}(yu_1+1)^{-\alpha(b-a+1)}dy}{\displaystyle\int_{0}^{1}y^{\alpha(b_2-a_2-1)}(yu_1+1)^{-\alpha(b-a+2)}dy}
\end{equation*}
\begin{equation*}\label{marL2}
	\delta_{\alpha,2}=\phi_{\alpha,2}(U_1)V_1, \mbox{ where }
	\phi_{\alpha,2}(u_1)=\frac{1}{(b-a)}\frac{\displaystyle \int_{0}^{1}y^{\alpha(b_2-a_2-1)}(yu_1+1)^{-\alpha(b-a)}dy}{\displaystyle\int_{0}^{1}y^{\alpha(b_2-a_2-1)}(yu_1+1)^{-\alpha(b-a+1)}dy}
\end{equation*}
and 
\begin{equation*}\label{marL3}
	\delta_{\alpha,3}=\phi_{\alpha,3}(U_1)V_1, \mbox{ where }
	\phi_{\alpha,3}(u_1)=\left( \frac{1}{(b-a+1)(b-a)}\frac{\displaystyle \int_{0}^{1}y^{\alpha(b_2-a_2-1)}(yu_1+1)^{-\alpha(b-a-1)}dy}{\displaystyle\int_{0}^{1}y^{\alpha(b_2-a_2-1)}(yu_1+1)^{-\alpha(b-a+1)}dy}\right)^{\frac{1}{2}}
\end{equation*}

\begin{theorem}\label{maruth}
	For estimating $\sigma_1$, risk of the estimator $\delta_{\alpha,i}$ is no where larger than that of $\delta^i_{01}$ with respect to $L_i(t)$, for $i=1,2,3$
\end{theorem}
\textbf{Proof:} We will give short proof for $L_1(t)$.  \cite{maruyama1998minimax} demonstrated that $\phi_{\alpha,1}(u_1)$ exhibits a non-decreasing trend and value greater than or equal to $\phi_{0,1}(u_1)$. Additionally, as $u_1$ approaches infinity, $\phi_{\alpha,1}(u_1)$ tends to  limit of $1/(b_1-a_1+1)$, which is very easy to show by using the transformation $y = \frac{t(1 + u_1)}{(1 + tu_1)}$ in $\phi_{\alpha,1}(u_1)$. Consequently, this follows from Theorem 2.2(i) in \cite{bobotas2017estimation}. Similarly we can proof the theorem for $L_2(t)$ and $L_3(t)$.
\begin{remark}\label{remBayes}\rm
	The boundary estimators $\delta_{\phi_{1}^{01}},\delta_{\phi_{1}^{02}},\delta_{\phi_{1}^{03}}$ are generalized Bayes for the loss functions $L_1(t),L_2(t),L_3(t)$ respectively. To show this we consider  prior as 
	$$\Pi(\theta)=\frac{1}{\sigma_{1}\sigma_{2}}, \ \  \mu_1<X_{a_1}, \ \ \mu_2<X_{a_2}, \ \ \sigma_1,\sigma_2>0, \ \ \sigma_{1}\le \sigma_{2}$$ 
	Posterior distribution can be obtained as
	$$\pi^*(\theta|\underline{X_a},\underline{V}) \propto \frac{1}{\sigma_1^{b_1-a_1+1}}\frac{1}{\sigma_2^{b_2-a_2+1}}e^{-\frac{v_1}{\sigma_1}}e^{-\frac{v_2}{\sigma_{2}}}$$
After some simplification we get the generalized Bayes for squared error loss is 
	$$\delta_B^{01}=\frac{\int_{0}^{\infty}\int_{0}^{\infty}\frac{1}{\sigma_{1}^{b_1-a_1+2}}\frac{1}{\sigma_{2}^{b_2-a_2+1}} e^{-\frac{v_1}{\sigma_1}}e^{-\frac{v_2}{\sigma_{2}}}d\sigma_{2}d\sigma_1}{\int_{0}^{\infty}\int_{0}^{\infty}\frac{1}{\sigma_{1}^{b_1-a_1+3}}\frac{1}{\sigma_{2}^{b_2-a_2+1}} e^{-\frac{v_1}{\sigma_1}}e^{-\frac{v_2}{\sigma_{2}}}d\sigma_{2}d\sigma_1}$$
Making the transformation as  $u_1=\frac{v_1}{\sigma_1}$ and $x=\frac{v_2\sigma_{2}}{v_1\sigma_1}$, we get 
$$\delta_B^{01}=\frac{\displaystyle\int_{0}^{\infty}\int_{0}^{\frac{v_2}{v_1}}e^{-(x+1)t_1}x^{b_2-a_2-1}t_1^{b-a}dxdt_1}{\displaystyle\int_{0}^{\infty}\int_{0}^{\frac{v_2}{v_1}}e^{-(x+1)t_1}x^{b_2-a_2-1}t_1^{b-a+1}dxdt_1}V_1$$
Similarly, we have obtained generalized Bayes estimators for entropy loss and symmetric loss functions respectively
$$\delta_B^{02}=\frac{\displaystyle\int_{0}^{\infty}\int_{0}^{\frac{v_2}{v_1}}e^{-(x+1)t_1}x^{b_2-a_2-1}t_1^{b-a-1}dxdt_1}{\displaystyle\int_{0}^{\infty}\int_{0}^{\frac{v_2}{v_1}}e^{-(x+1)t_1}x^{b_2-a_2-1}t_1^{b-a}dxdt_1}V_1$$
$$\delta_B^{03}=\left( \frac{\displaystyle\int_{0}^{\infty}\int_{0}^{u_1}e^{-(x+1)t_1}x^{b_2-a_2-1}t_1^{b-a-2}dxdt_1}{\displaystyle\int_{0}^{\infty}\int_{0}^{u_1}e^{-(x+1)t_1}x^{b_2-a_2-1}t_1^{b-a}dxdt_1}\right)^{\frac{1}{2}}V_1$$
	\end{remark}
\subsection{Strawderman type improved estimator}
In this subsection we have obtained Strawderman-type improved estimator for $\sigma_1$ for quadratic loss and entropy loss function. 
\begin{theorem}\label{man1}
The estimator 
$$\delta^1_{\phi}=(1-\phi(Z_1))c_{01}V_1,$$ has smaller risk than $\delta^1_{01}$ with respect to the quadratic loss function $L_1(t)$  provided, for some $\epsilon>0$ the function $\phi(z_1)$ satisfy the following conditions:
\begin{itemize}
	\item [(i)] $(1+z_1)^{\epsilon}\phi(z_1)$ is decreasing in $(0,\infty).$
\item [(ii)] $0 \le \phi(z_1) \le r^*(c_{01})= \min\left\lbrace r(c_{01}),r_1(c_{01})\right\rbrace $, where $r(c_{01})= \frac{(2(c_{01}((b_1-a_1)+2)-1)}{c_{01}(b_1-a_1+2)}$ and  $$r_1(c_{01})= \frac{\Gamma(b_1-a_1+1)(\Gamma(b-a+\epsilon+1)}{\Gamma(b-a+1)(\Gamma(b_1-a_1+\epsilon+1)}-\frac{2\Gamma(b_1-a_1)\Gamma(b-a+\epsilon+1)}{c_{01}(b-a+2)\Gamma(b_1-a_1+\epsilon+1)\Gamma(b-a)}$$
\end{itemize}
\end{theorem}
\noindent\textbf{Proof:} To prove this theorem we invoke the argument similar to \cite{maruyama2006new}, and \cite{bobotas2019improved}. Let $\psi(z_1)=(1+z_1)^{\epsilon}\phi(z_1)$ and by the assumption we have $\phi(z_1) \le \lim \limits_{z_1\rightarrow 0}\psi(w)=\lim \limits_{z_1\rightarrow 0}\phi(z_1)=r$. The risk difference between $\delta^1_{01}$ and $\delta_{\phi}$  is 
\begin{align*}
RD\left(\theta;\delta^1_{01}, \delta_{\phi}\right)&=E\left( c_{01}\frac{ V_1}{\sigma_1}-1\right) ^2-E\left(c_{01}\{1-\phi(Z_1)\}\frac{V_1}{\sigma_{1}}-1\right)^2  \\	
&\ge E\left[c_{01} \frac{\psi(Z_1)}{(1+Z_1)^{\epsilon}}\left\lbrace 2c_{01}Y^2-c_{01}\frac{r}{(1+Z_1)^{\epsilon}}Y^2-2Y\right\rbrace  \right] 
\end{align*}
as we have $V_1 \sim \Gamma(b_1-a_1,\sigma_1)$ and  $\sigma_1 \le \sigma_2$ and $V_2 \sim \Gamma(b_2-a_2,\sigma_2)$ this implies $V_2 \sim \sigma_{1}\Gamma(b_2-a_2,\frac{\sigma_2}{\sigma_1})$ this implies $\frac{V_2}{\sigma_1}\sim \Gamma(b_2-a_2,\frac{\sigma_2}{\sigma_1})$ 
from Corollary $1$  in \cite{fearnhead2004exact} and \cite{bobotas2019improved}, we can write $\frac{V_2}{\sigma_1}$ follows  $G(b_2-a_2+L,1)$ where $L$ follows $NB(b_2-a_2,1-\frac{\sigma_{1}}{\sigma_{2}})$. Therefore,
\begin{align*}
RD&\ge E_L\left[ E\left[c_{01} \frac{\psi(Z_1)}{(1+Z_1)^{\epsilon}}\left\lbrace 2c_{01}Y^2-c_{01}\frac{r}{(1+Z_1)^{\epsilon}}Y^2-2Y\right \rbrace \bigg|L=l\right]\right]
\end{align*}
We proceed by computing inner conditional expectation. Given $L=l$ and setting $Z_l =Z_1 |L=l$ we have that $YZ_l$ follows $\Gamma(b_2-a_2+L,1)$. 
\begin{align*}
RD&\ge E\left\lbrace E\left[ c_{01} \frac{\psi(W)}{(1+Z_1)^\epsilon}\frac{(b-a+l)(b-a+l+1)}{(Z_l+1)^2} \left\lbrace2c_{01}-\frac{rc_{01}}{(1+Z_l)^\epsilon}-2\frac{1+ Z_l}{b_1-a_1+l+1} \right \rbrace |L=l\right]\right\rbrace
\end{align*}
Transforming $U_l=(1+z_l)^{-1}$ where $U_l$ has $Beta(b_1-a_1,b_2-a_2)$ distribution. So finally we get
\begin{equation*}
\ge E\left\lbrace E \left[c_{01}^2 \psi\left( \frac{1-U_l}{U_l}\right) U_l^{\epsilon+2}(b-a+l)(b-a+l+1)\lbrace2-rU_l^{\epsilon}-\frac{2}{c_{01}U_l(b-a+l+1)}\rbrace \big|L=l\right]\right\rbrace
\end{equation*} 
The expression $f(u)=2U_l-rU_l^{\epsilon+1}-\frac{2}{c_{01}(b-a+l+1)}$ changes sign at least once in $[0,1]$ provided 

\begin{equation}
	r \le \frac{(2(c_{01}((b-a)+l+1)-1)}{c_{01}(b-a+l+1)}
\end{equation}
Using lemma A.1 in \cite{bobotas2011improved}
 $$ E\left\lbrace E \left[c_{01}^2 \psi\left( \frac{1-U_l}{U_l}\right) U_l^{\epsilon+1}(b-a+l)(b-a+l+1)\lbrace2U_l-rU_l^{\epsilon+1}-\frac{2}{c_{01}(b-a+l+1)}\rbrace\big|L=l\right]\right\rbrace  $$ 
 $$\ge c_{01}^2 \psi\left( \frac{1-U_{l0}}{U_{l0}}\right) U_{l0}^{\epsilon+1}(b-a+l)(b-a+l+1) E_L\left\lbrace E \left[\lbrace2U_l-rU_l^{\epsilon+1}-\frac{2}{c_{01}(b-a+l+1)}\rbrace\big|L=l\right]\right\rbrace $$
 Using the argument similar to \cite{bobotas2019improved} the right hand side expectation will be non-negative provided
 $$r \le \frac{\Gamma(b_1-a_1+1)(\Gamma(b-a+\epsilon+1)}{\Gamma(b-a+1)(\Gamma(b_1-a_1+\epsilon+1)}-\frac{2\Gamma(b_1-a_1)\Gamma(b-a+\epsilon+1)}{c_{01}(b-a+l+1)\Gamma(b_1-a_1+\epsilon+1)\Gamma(b-a)}
 =r_1(c_{01},l)$$
 $r_1(c_{01},l)$ is non negative and non decreasing in $l$ and taking $r_1(c_{01})= r_1(c_{01},1)$ completes the proof. 
 
Now will derive the Strawderman-type improved estimator under the entropy loss function. We prove the following lemma using the argument similar to \cite{pal1995improved}. 
 \begin{Lemma}\label{Le1}
  $f(u)=\alpha_2(b-a+l)u^{\epsilon+1}+\frac{1}{r}\log(1-ru^{\epsilon})$ has exactly one sign change from negative to positive provided $r\le\frac{1}{1+\epsilon}$   \end{Lemma}
\textbf{Proof:}  $f(u_l)=\alpha_2(b-a+l)u^{\epsilon+1}+\frac{1}{r}\log(1-ru^\epsilon)$, $u_l$ varries from $0$ to $1$.  Since $f(0)=0$ and $f(1)= \alpha_2(b-a+l)+\frac{1}{r}\ln(1-r)$. For one sign change $f(u)$  is increasing if $u-ru^{\epsilon+1}\ge \frac{\epsilon}{(\epsilon+1)\alpha_2(b-a+l)}$ and decreasing if $u-ru^{\epsilon+1}\le \frac{\epsilon}{(\epsilon+1)\alpha_2(b-a+l)}$. Let $w(u)=u-ru^{\epsilon+1}$ and $w(u)$ increases if $r \le \frac{1}{\epsilon+1}$. As we increase $u$, $w(u)$  will  increases provied $r \le \frac{1}{\epsilon+1}$.  So $f(u)$ will increase and eventually  it cuts the horizontal axis as $f(1)>0$ and $u \rightarrow 0$ \\
\begin{theorem}
The estimator 
$$\delta^2_{\phi}=(1-\phi(Z_1))c_{01}V_1,$$ dominates $\delta^2_{01}$ with respect to the entropy loss function $L_2(t)$ provided, 
for some $\epsilon >0$  the function $\phi(z_1)$ satisfy the following conditions:
\begin{enumerate}
\item[(i)] $(1+z_1)^\epsilon \phi(z_1)$ is decreasing
in $(0, \infty)$.
\item[(ii)] $0 \le \phi(z_1) \le B^*(\epsilon)$  with $B^*(\epsilon)=\min\{r^*,\frac{1}{1+\alpha_2},\frac{B(\epsilon)}{b-a}\}$, where  $r^*\in (0,1)$ is unique solution of the equation $c_{01}(b-a)+\frac{1}{r}\ln(1-r)=0$ and   $B(\epsilon)=\frac{\Gamma(b_1-a_1+\epsilon)\Gamma(b_1-a_1+\epsilon+1+b_2-a_2)}{\Gamma(b-a+\epsilon)\Gamma(b_1-a_1+\epsilon+1)}$.
\end{enumerate}
\end{theorem}
\textbf{Proof:}  To prove this theorem we have used the techniques of \cite{maruyama2006new}, \cite{bobotas2011improved} and \cite{bobotas2019improved}. The risk difference of $\delta^2_{01}$ and $\delta^2_{\phi}$ is 
\begin{align*}
	RD& = E\left(c_{01} \frac{V_1}{\sigma_1}-\ln(\frac{c_{01} V_1}{\sigma_1})-1\right)-E\left(c_{02}\{1-\phi(z_1)\}\frac{V_1}{\sigma_1}-\ln(c_{02}\{1-\phi(z_1)\}\frac{V_1}{\sigma_1})-1\right)\\
 &\ge E\left[\psi(Z)\left\lbrace\frac{1}{(1+Z_1)^{\epsilon}}c_{01}Y_1+\frac{1}{r}\ln\left(1-\frac{r}{(1+Z_1)^{\epsilon}}\right) \right\rbrace \right] 
\end{align*}
Again using similar argument as in Theorem \ref{man1}. $L \sim NB(b_2-a_2,1-\frac{\sigma_1}{\sigma_2})$, Now taking conditional on $L$
$$RD \ge E\left\{E\left[\psi(Z_1)\left\lbrace\frac{1}{(1+Z_1)^{\epsilon}}c_{01}Y_1+\frac{1}{r}\ln\left(1-\frac{r}{(1+Z_1)^{\epsilon}}\right) \right\rbrace \big|L\right] \right\}$$
Suppose $Z_l=Z_1|L=l$ then we have $YZ_l$ follows $\Gamma(b_2-a_2+L,1)$. So the risk difference becomes 
\begin{align*}
	RD&\ge E\left\lbrace E\left[\psi(Z_l)\left\lbrace\frac{c_{01}}{(1+Z_l)^{\epsilon}}E(Y_1|Z_l)+\frac{1}{r}\ln\left(1-\frac{r}{(1+Z_l)^{\epsilon}}\right) \right\rbrace \big|L=l\right]\right\rbrace
\end{align*}
Proceeding by computing the inner expectations.
\begin{align*}
E\left\lbrace E\left[\psi(Z_l)\left\lbrace\frac{c_{01}}{(1+Z_l)^{\epsilon}}\frac{b-a+l}{1+Z_l}+\frac{1}{r}\ln\left(1-\frac{r}{(1+Z_l)^{\epsilon}}\right) \right\rbrace \big|L=l\right]\right\rbrace
\end{align*}
Transforming  $\frac{1}{1+z_l}=u_l$
we get
$$E\left\lbrace E\left[\psi(\frac{1-u_l}{u_l})\left\lbrace{c_{01}}{(u_l)^{\epsilon+1}}(b-a+l)+\frac{1}{r}\ln\left(1-r{u_l^{\epsilon}}\right) \right\rbrace \big|L=l\right]\right\rbrace$$
Let $K_l(u_l)={c_{01}(b-a+l)}{(u_l)^{\epsilon+1}}+\frac{1}{r}\ln\left(1-{r}{u_l^{\epsilon}}\right)$. Using Lemma \ref{Le1} we can conclude $K_l(u_l)$ changes sign once from negative to positve as $u_l$ going to $0$ to $1$ provided $r\le \frac{1}{1+\epsilon}$ and $r$ should be in $(0,r^*)$, where $r^*$ is the root of $c_{01}(b-a)+\frac{1}{r}\ln(1-r)=0$. Since for $K_l(1)>0$ we must have $c_{01}(b-a)+\frac{1}{r}\ln(1-r)>0$ and $c_{01}(b-a) >1$ because $\frac{1}{r}\ln(1-r)$ is decreasing in $r$ and bounded above by $\lim_{r\rightarrow 0} \frac{1}{r}\ln(1-r)=-1$. Let $r^* \in (0,1)$  is unique solution of $c_{01}(b-a)+\frac{1}{r}\ln(1-r)=0$. So for $c_{01}(b-a)+\frac{1}{r}\ln(1-r)>0$ we should have $ r \in (0,r^*)$.
\begin{align*}
E[k_l(U_l)|L=l]&=E[c_{01}(b-a+l)U_l^{\epsilon+1}+\frac{1}{r}\log(1-rU_l^{\epsilon})|L=l] 
\end{align*}
Using the argument given in \cite{pal1995improved} we have obtained
$E[k_l(U_l)|L=l]\ge0$ provided $(1-r)\ge B(\epsilon) C_l$
where, $B(\epsilon)=\frac{\Gamma(b_1-a_1+\epsilon)\Gamma(b_1-a_1+\epsilon+1+b_2-a_2)}{\Gamma(b-a+\epsilon)\Gamma(b_1-a_1+\epsilon+1)}$ and $C_l=\frac{1}{b-a+l}$.
Since $C_l$ is decreasing function of $l$ and $C_l\le\frac{1}{b-a}$
this implies  $E[k_l(U_l)|L=l]\ge0$ provided $r\le \frac{B}{b-a}$.
\section{Improved estimation of $\sigma_2$}\label{sect3}
\noindent In this section, we are considering the problem of estimation of scale parameter 
$\sigma_2$ under the order restriction $\sigma_{2}\le\sigma_2$. Where location parameters $\mu_{i}s$ are unknown. We will derive estimators which dominate the BAEE of $\sigma_2$. To derive the dominating estimator we consider the class of estimators of the form
\begin{equation}\label{Baee*}
	C_1^*=\left\lbrace 	\delta_{\phi_{1}^*}=V_2\phi_1^{*}(Z^*):Z^*=\frac{V_1}{V_2},  \text{and} \ \phi_1^* \  \text{is positive  measurable function}\right\rbrace 
\end{equation}
Using the argument similar to Theorem \ref{Th1sigma1} we can proof the following theorem. 
\begin{theorem}\label{th1sigma2}
	Let $\beta^*$ be the unique solutions of the \begin{equation}
		E[L'(\beta^*U)]=0
	\end{equation}
	where $U$ follows $Gamma(b-a+1,1)$ distribution. Define a function $\phi_{2,1}=max\{\phi_1^*,c(1,z^*)\}$ with $c(1,z)=\beta^*(1+z^*).$ Then the estimator $\delta_{\phi_{2,1}}$ has uniformly smaller risk than $\delta_{\phi_1}^*$ under a scale invariant loss function $L(.)$, provided $P_\eta(\phi_{1}^*(Z^*)<c(1,Z^*))\neq0$ for some $\eta>0$.
\end{theorem}
\begin{corollary}\label{ST21}
	The risk of the estimator $\delta_{2S1}=\max\{c_{02},\beta^*(1+z^*)\}V_2,$ is nowhere larger than $\delta_{02}$ under a bowl shaped loss function $L(t)$ provided $\beta^*\ge c_{02}$.
\end{corollary}
\textbf{Proof:} Taking $\phi_1^*(Z^*)=c_{02}$ we have $c(\eta,z^*)>\phi_{2,1}^*(z^*)>c_{02}$ on a set of positive probability as we have $\beta^*>c_{02}$. Now,  $R_1^*(\theta,c) $ is strictly bowl shaped. So $R_1^*(\theta,c)$ is strictly decreasing in $(c(\eta,z^*),c)$. Hence, $R_1^*(\theta,\delta_{\phi_{2S1}})<R^*(\theta,\delta_{02})$ and $z^*>0$. Hence we get the result.~~~~~~~~~~ $\blacksquare$
Consider the another class of estimator 
\begin{equation}
	C_2^*=\left\lbrace 	\delta_{\phi_{2}^*}=\phi_2^*(Z_1^*)V_2, Z_1^{*}=\frac{X_{2a}}{V_2} \right\rbrace 
\end{equation}
\begin{theorem}
	Suppose that a random variable $W$ has $Gamma(b_2-a_2+2,1)$ distribution and let $\beta_1^*$ is unique solution of $$EL'(W\beta_1^*)=0.$$ Define a function $\phi_{22}=\min\{\phi_{2}(Z_1^*),\beta_1^*((n_2-a_2+1)Z_1^*+1)\}$ for $Z_1^*>0$ then the estimator
	\begin{eqnarray*}
		\delta_{\phi_{22}}=\left\{
		\begin{array}{ll}
			\phi_{22}(Z_1^*)V_2,~Z_1^*>0\\\\
			\phi_{2}^*(Z_1^*)V_2,~~~~~~~~~\mbox{ otherwise. }
		\end{array}
		\right.
	\end{eqnarray*}
	dominates the estimator $\delta_{\phi_2^*}$ with respect to a scale invariant loss function provided $P(\phi_{2}^*(Z_1^*)>\beta_1^*(1+ (n_2-a_2+1)Z_1^*) \ne 0$ for $Z_1^*>0$. 
\end{theorem}
\textbf{Proof}: Proof is similar to Theorem \ref{th2sigma1}.
\begin{corollary}\label{th2sigma2}
	The estimator \begin{eqnarray}
		\delta_{2S2}=\left\{
		\begin{array}{ll}
			\min \left\{c_{02},\beta_{1}^*(1+(n_2-a_2+1)Z_1^*)\right\}V_2,~Z_1^*>0\\\\
			c_{02}V_2,~~~~~~~~~~~~~~~~~~~\mbox{ otherwise. }
		\end{array}
		\right.
	\end{eqnarray}
	dominates $\delta_{02}$ under scale invariant loss function $L(t)$ provided $\beta_{1}^*<c_{02}$.
\end{corollary}

\subsection{Double shrinkage improved estimators }
Now we will propose double shrinkage estimator involving the estimators $\delta_{2S1}$ and $\delta_{2S2}$ following the techniques of \cite{iliopoulos1999improving}. Consider the estimator $\delta_{2S1}$ and $\delta_{2S2}$ defined in corollary \ref{ST21} and \ref{th2sigma2} respectively. Let $ \phi_{21},\phi_{22}$ and $L(t)$ are almost differentiable. $ \phi_{21}(z^*)$ is non increasing and $\phi_{22}(z_1^*)$ is non decreasing. Also 
$\lim _{z^*\rightarrow\infty }\phi_{21}(z^*)=\lim _{z_1^*\rightarrow\infty} \phi_{22}(z_1^*)=c_{02}$. Then we have the following result.
\begin{theorem}
The estimator $\delta_{2S3}=\{\phi_{21}+\phi_{22}-c_{02}\}V_2$ dominates $\delta_{2S1}$ as well as $\delta_{2S2}$ provided $L'(t)$ is non decreasing.
\end{theorem}
\textbf{Proof:} Proof is similar to the Theorem 3.5 in \cite{iliopoulos1999improving}. 
\subsection{Improved estimators of $\sigma_2$ for special loss functions} 
In this subsection we will derive improved estimators for various special loss functions as an application of Theorem  \ref{th1sigma2}, \ref{th2sigma2} . Using the Theorem \ref{th1sigma2}, \ref{th2sigma2}   improved estimators for quadratic loss function are
\begin{eqnarray*}
	\delta_{2S1}=\left\{
	\begin{array}{ll}
		\frac{V_2}{b-a+1}(1+ Z^*),~ \text{if}\ \ Z^* > \frac{b-a}{b_2-a_2+1}-1\\
		\frac{1}{b_2-a_2+1}V_2,~~~~~~~~~~~~~~~~~~~\mbox{ otherwise. }
	\end{array}
	\right.
\end{eqnarray*}

\begin{eqnarray*}
	\delta_{2S2}=\left\{
	\begin{array}{ll}
		\min\left\{\frac{(1+(n_2-a_2+1)Z_1^*)}{b_2-a_2+2},\frac{1}{b_2-a_2+1}\right\}V_2,~Z_1^*>0,\\\\
		\frac{1}{b_2-a_2+1}V_2,~~~~~~~~~~~~~~~~~~~\mbox{ otherwise. }
	\end{array}
	\right.
\end{eqnarray*}


By the Theorems \ref{th1sigma2}, \ref{th2sigma2}  improved estimators for entropy loss function are
\begin{eqnarray*}
	\delta_{2S1}=\left\{
	\begin{array}{ll}
		\frac{1}{b-a}(1+ Z^*)V_2, \text{if}\ \   Z^*>\frac{b-a}{b_2-a_2}-1\\
		\frac{1}{b_2-a_2}V_2,~~~~~~~~~~~~~~~~~~~\mbox{ otherwise. }
	\end{array}
	\right.
\end{eqnarray*} \\
\begin{eqnarray*}
	\delta_{2S2}=\left\{
	\begin{array}{ll}
		\min\left\{\frac{(1+(n_2-a_2+1)Z_1^*)}{b_2-a_2+1},	\frac{1}{b_2-a_2}\right\}V_2,~Z_1^*>0\\\\
		\frac{V_2}{b_2-a_2},~~~~~~~~~~~~~~~~~~~\mbox{ otherwise. }
	\end{array}
	\right.
\end{eqnarray*}
Using the Theorem \ref{th1sigma2}, \ref{th2sigma2}  improved estimators for symmetric loss function $L_3(t)$ are
\begin{eqnarray*}
	\delta_{2S1}=\left\{
	\begin{array}{ll}
		\left\{\frac{1+ Z^*}{\sqrt{b-a-1}\sqrt{b-a}}\right\}V_2, \text{if}\ Z^*>\frac{\sqrt{(b-a-1)(b-a)}}{\sqrt{(b_2-a_2)(b_2-a_2-1)}}-1\\
		\frac{1}{\sqrt{(b_2-a_2)(b_2-a_2-1)}}V_2,~~~~~~~~~~~~~~~~~~~\mbox{ otherwise. }
	\end{array}
	\right.
\end{eqnarray*}

\begin{eqnarray*}
	\delta_{2S2}=\left\{
	\begin{array}{ll}
		\min\left\{\frac{(1+(n_2-a_2+1)Z_1^*)}{\sqrt{(b_2-a_2+1)(b_2-a_2)}},	\frac{1}{\sqrt{(b_2-a_2)(b_2-a_2-1)}}\right\}V_2,~Z_1^*>0\\\\
		\frac{V_2}{\sqrt{(b_2-a_2)(b_2-a_2-1)}},~~~~~~~~~~~~~~~~~~~\mbox{ otherwise. }
	\end{array}
	\right.
\end{eqnarray*}
%

\begin{remark}
From \cite{fernandez2002computing}  we have the unrestricted MLE of $\sigma_2$  is $\delta_{Rmle}^*=\frac{V_2}{b_2-a_2+1}$. 
From \cite{misra2002smooth}	The restricted maximum likelihood estimator for $\sigma_2$ with order restriction $\sigma_{1}\le \sigma_{2}$ can be obtain as 
$$\delta_{Rmle}^*=
\max \left\lbrace \frac{V_2}{b_2-a_2+1},\frac{V_1+V_2}{b-a+2}\right\rbrace=\max\left\lbrace\frac{1}{b_2-a_2+1},\frac{1+Z^*}{b-a+2}\right\rbrace V_2=\phi_{Rmle}^*(Z^*)V_2.$$
Using Theorem \ref{th1sigma2}, we get the estimator $\delta_{\phi_2^r}^* = \phi_2^r(Z^*)V_2$, where $\phi_2^r(Z^*) = \max\{\phi_{RML2}(Z^*), \beta^*(1 + Z^*)\}$, dominates $\delta_{Rmle}^*$ for estimating $\sigma_{2}$, with order restriction $\sigma_{1}\le\sigma_{2}$, under a general scale invariant loss function $L(t)$.
\end{remark}

\subsection{A class of improved estimators of $\sigma_2$ }
In this section, we have derived a class of improved estimators for $\sigma_2$ using IERD approach of \cite{kubokawa1994unified}. 
\begin{theorem}
	The risk of estimator $\delta_{\phi_{1}}(\underline{X_a},\underline{V})$ given in (\ref{Baee*}) is nowhere larger than $\delta_{02}(\underline{X_a},\underline{V})$ under a bowl shaped loss function $L(t)$ provided $\phi_{1}(u_1)$ satisfies the following conditions:
	\begin{enumerate}
		\item[(i)]$\phi_2(u_2)$ is non-decreasing and $lim_{u_2\rightarrow0}\phi_2(u_2) = c_{02},$
		\item[(ii)] $\int_{0}^{\infty}\int_{t_2u_2}^{\infty}L'(\phi_2(u_2)t_2)t_2
		\pi_1(x)\pi_2(t_2)dx dt_2 \le 0$,
	\end{enumerate}
	where $\pi_i$ is the pdf of a $Gamma( b_i-a_i$, 1), i = 1, 2. 
\end{theorem}
\begin{corollary}
		Consider the quadratic loss function $L_1(t)=(t-1)^2$. The risk of estimator $\delta_{\phi_{2}}(\underline{X_{a}},\underline{V})$given in \ref{imp1} is nowhere larger than 
	$\delta_{02}(\underline{X_a},\underline{V})$ provided $\phi_{2}(u_2)$ satisfies the following conditions.
	\begin{enumerate}
		\item [(i)] $\phi_{2}(u_2)$ is non decreasing and $\lim \limits_{u_2\rightarrow0}\phi_{2}(u_2)=\frac{1}{b_2-a_2+1}$
		\item [(ii)]$\phi_{2}(u_2)\le\phi_{2}^{01}(u_2)$ where 
		$$\phi_{2}^{01}(u_2)=\frac{\displaystyle\int_{0}^{\infty}\int_{t_2u_2}^{\infty}e^{-x}x^{b_1-a_1-1}e^{-t_2}t_2^{b_2-a_2}dxdt_2}{\displaystyle\int_{0}^{\infty}\int_{t_2u_2}^{\infty}e^{-x}x^{b_1-a_1-1}e^{-t_2}t_2^{b_2-a_2+1}dxdt_2}$$
	\end{enumerate}
	
\end{corollary}
\begin{corollary}
	Consider the entropy loss function $L_2(t)=t-\ln t-1$. The risk of estimator $\delta_{\phi_{2}}(\underline{X_{a}},\underline{V})$given in \ref{imp1} is nowhere larger than 
	$\delta_{02}(\underline{X_a},\underline{V})$ provided $\phi_{2}(u_2)$ satisfies the following conditions.
	\begin{enumerate}
		\item [(i)] $\phi_{2}(u_2)$ is non decreasing and $\lim \limits_{u_2\rightarrow0}\phi_{2}(u_2)=\frac{1}{b_2-a_2}$
		\item [(ii)]$\phi_{2}(u_2)\le\phi_{2}^{02}(u_2)$ where 
		$$\phi_{2}^{02}(u_2)=\frac{\displaystyle\int_{0}^{\infty}\int_{t_2u_2}^{\infty}e^{-x}x^{b_1-a_1-1}e^{-t_2}t_2^{b_2-a_2-1}dxdt_2}{\displaystyle\int_{0}^{\infty}\int_{t_2u_2}^{\infty}e^{-x}x^{b_1-a_1-1}e^{-t_2}t_2^{b_2-a_2}dxdt_2}$$
	\end{enumerate}
	
\end{corollary}
\begin{corollary}
	Consider the symmetric loss function $L_3(t)=t+\frac{1}{t}-2$. The risk of estimator $\delta_{\phi_{2}}(\underline{X_{a}},\underline{V})$given in \ref{imp1} is nowhere larger than 
	$\delta_{02}(\underline{X_a},\underline{V})$ provided $\phi_{2}(u_2)$ satisfies the following conditions.
	\begin{enumerate}
		\item [(i)] $\phi_{2}(u_2)$ is non decreasing and $\lim \limits_{u_2\rightarrow0}\phi_{2}(u_2)=\frac{1}{\sqrt{(b_2-a_2)(b_2-a_2-1)}}$
		\item [(ii)]$\phi_{2}(u_2)\le\phi_{2}^{03}(u_2)$ where 
		$$ \phi_{2}^{03}(u_2)=\left( \frac{\displaystyle\int_{0}^{\infty}\int_{t_2u_2}^{\infty}e^{-x}x^{b_1-a_1-1}e^{-t_2}t_2^{b_2-a_2-2}dxdt_2}{\displaystyle\int_{0}^{\infty}\int_{t_2u_2}^{\infty}e^{-x}x^{b_1-a_1-1}e^{-t_2}t_2^{b_2-a_2}dxdt_2}\right)^\frac{1}{2}$$
	\end{enumerate}
\end{corollary}
\subsection{Maruyama type improved estimator}
Consider the following estimators for $\alpha>1$
\begin{equation*}\label{marL4}
	\delta_{\alpha,1}^{*}=\phi_{\alpha,1}^*(U_2)V_2~ \mbox{ with }~
	\phi^*_{\alpha,1}(u_2)=\frac{1}{b-a+1}\frac{\displaystyle \int_{1}^{\infty}y^{\alpha(b_1-a_1-1)}(yu_2+1)^{-\alpha(b-a+1)}dy}{\displaystyle\int_{1}^{\infty}y^{\alpha(b_1-a_1-1)}(yu_2+1)^{-\alpha(b-a+2)}dy}
\end{equation*}
\begin{equation*}\label{marL5}
	\delta_{\alpha,2}^{*}=\phi_{\alpha,2}^*(U_2)V_2, \mbox{ where }
	\phi^*_{\alpha,2}(u_2)=\frac{1}{(b-a)}\frac{\displaystyle \int_{1}^{\infty}y^{\alpha(b_1-a_1-1)}(yu_2+1)^{-\alpha(b-a)}dy}{\displaystyle\int_{1}^{\infty}y^{\alpha(b_1-a_1-1)}(yu_2+1)^{-\alpha(b-a+1)}dy}
\end{equation*}
and 
\begin{equation*}\label{marL6}
	\delta_{\alpha,3}^{*}=\phi_{\alpha,3}^*(U_2)V_2, \mbox{ where }
	\phi^*_{\alpha,3}(u_2)=\left( \frac{1}{(b-a)(b-a-1)}\frac{\displaystyle \int_{1}^{\infty}y^{\alpha(b_1-a_1-1)}(yu_2+1)^{-\alpha(b-a-1)}dy}{\displaystyle\int_{1}^{\infty}y^{\alpha(b_1-a_1-1)}(yu_2+1)^{-\alpha(b-a+1)}dy}\right)^{\frac{1}{2}}
\end{equation*}

\begin{theorem}
For estimating $\sigma_2$, the risk of the estimator  $\delta_{\alpha,i}^{*}$ is no where larger than that of  $\delta^i_{02}$ with respect to $L_i(t)$ for $i=1,2,3.$
\end{theorem}
\textbf{Proof:} The idea of  the proof is close to the Theorem \ref{maruth}.

\begin{remark}
	Similarly the Remark \ref{remBayes} we can prove that the estimators $\delta_{\phi_{0,2}^1},\delta_{\phi_{0,2}^2},\delta_{\phi_{0,2}^3}$ generalized Bayes with respect to the loss functions $L_1(t),L_2(t),L_3(t)$ under the prior $$\Pi(\theta)=\frac{1}{\sigma_{1}\sigma_{2}}, \ \  \mu_1<X_{a_1}, \ \ \mu_2<X_{a_2}, \ \ \sigma_1,\sigma_2>0, \ \ \sigma_{1}\le \sigma_{2}$$
	
\end{remark}
 \section{Special sampling schemes}\label{sect4}
 \subsection{I.I.D. sample} In our model if we take $a_i=1$ and $b_i=n_i$ then we get the results for case of simple random sampling scheme.  For for simple random sampling one can see \cite{patra2021componentwise}.
 \subsection{Type-II Censoring}
 There are a number of scenarios in the field of reliability and life-testing tests when test units are either misplaced or removed from the experiment before they break. For example, in clinical studies, participants may withdraw or the experiment may cease owing to financial constraints, and in industrial experiments, units may break suddenly before the planned time. A deliberate choice to remove units before they fail is made in many circumstances, mainly to save time and lower testing costs. Censored data are those obtained from these kinds of experiments; type-II censoring is one type of censored data. The experimenter chooses to end the experiment in this method after a predetermined number of items ($r$, where $r$ is less than or equal to $n$) fail.
 In this subsection we consider estimation of $\sigma_{1}$ and 
 $\sigma_{2}$ with $\sigma_{1} \le \sigma_{2}$ under type II censoring scheme. We are considering a model in which we observe the first $n_i$ order statistics are suppose $X_{(i1)}\le X_{(i2)} \le \cdots \le X_{(i n_i)}$. from a random sample of size $N_i$, $n_i\le N_i$ follows exponential distribution $E(\mu_i,\sigma_i), i=1,2.$ The complete  sufficient statistics is $(\underline{X_{(1)}}, \underline{T})$. where,
 $\underline{X_{(1)}} = (X_{1(1)}, X_{2(1)})$ and $\underline{T} = (T_1, T_2)$. Moreover, $(X_{1(1)}, X_{2(1)}, T_1, T_2)$ are independently distributed with $X_{i(1)} \sim E(\mu_i, \sigma_{i} ) $and $T_i \sim Gamma(n_i - 1, \sigma_{i})$, i = 1, 2.
 The complete  sufficient statistics is $(\underline{X_{(1)}}, \underline{T})$. Where,
 $\underline{X_{(1)}} = (X_{1(1)}, X_{2(1)})$ and $\underline{T} = (T_1, T_2)$ also $(X_{1(1)}, X_{2(1)}, T_1, T_2)$ are independently
	distributed with $X_{i(1)} \sim E(\mu_i, \sigma_i )$ and $T_i \sim Gamma(n_i - 1, \sigma_{i}), i = 1, 2.$ Model is similar to \cite{patra2021componentwise}.
	\subsection{ Progressive Type-II censoring }
	Because of their extensive relevance, lifetime distributions evaluated under censored sample methods are frequently used in many domains, including science, engineering, social sciences, public health, and medicine. When censored data is handled by eliminating a predefined number of surviving units at the point of individual unit failure, it is classified as progressive type-II censoring. In a life-testing experiment, only $m$ (where $m < n$) are watched continually until they break. $n$ units are involved. Gradually, censoring occurs in $m$ steps.  For the $m$ units that are fully observed, these $m$ phases give failure times. 
	When the first failure (the first stage) occurs, $R_1$ of the $n-1$ remaining units are taken out at random based on the experiment. Then, when $n-2-R_1$ remaining units, $R_2$ are randomly withdrawn at the time of the second failure (the second stage). Ultimately, all of the remaining $R_m = n - m - R_1 - R_2 - \dots - R_m-1$ units are withdrawn when the $m$-th failure (the $m$-th stage) takes place. Suppose we have datasets as one set is $X_{1:n_i:m_i}^p, X_{2:n_i:m_i}^p,\dots, X_{n_i:n_i:m_i}^p$, representing the progressive type-II censored sample of lifetimes for $n_1$ independent units, following an exponential distribution $E(\mu_{1},\sigma_{1})$. The other set is $Y_{1,n_i:m_i}^p, Y_{2,n_i:m_i}^p,\dots, Y_{n_i:n_i:m_i}^p$, representing the lifetimes of another set of $n_2$ independent units following an exponential distribution $E(\mu_2 ,\sigma_{2})$ associated restriction $\sigma_{1}\le \sigma_{2}.$
	Here we are considering the problem of estimation of  parameter $\sigma_{i}$ based on these progressive type-II censored samples with order restriction $\sigma_{1}\le \sigma_{2}$ To estimate this, we define  here $(X_{1,n_1:n_2}^p,T_1)$ as complete sufficient statistic for and $(Y_{1,n_1:n_2}^p,T_2)$ are complete sufficient statistics. Where $T_1 = \sum_{j=1}^{n_2}(R_i+1)X_{j,n_1:n_2}^p-X_{1,n_1:n_2}^p$ and $T_2= \sum_{j=1}^{n_2}(R_i+1)Y_{j,n_1:n_2}^p-Y_{1,n_1:n_2}^p)$.
	It's then shown that $X_{1:n_1:n_2}^p$ follows $E(\mu_1,\frac{\sigma_1}{n_1})$, $Y_{1:n_1:n_2}^p$ follows $E(\mu_2,\frac{\sigma_2}{n
	_2})$ for  and $T_1$, $T _2$ follows a gamma distribution $\Gamma(n_1-1, \sigma_1)$, $\Gamma(n_2-1,\sigma_{2})$ respectively. This model is again similar to \cite{patra2021componentwise}
	\subsection{Record values}
	Record values are much applicable in diverse fields like hydrology, meteorology, and stock market analysis. Record values can be defined as:
Let $V_1, V_2,\dots, V_n$ be the sequence of i.i.d. random variables. We define a new sequence of random variables such that  for $k = 1$ set  $U(1) = 1$ and $U(k)$=min$\{j|j>U(k-1),V_j>V_{u(k-1)}\}$. This new sequence $\{X_k = V_{U(k)}, k \ge 1\}$, represents a sequence of maximal record statistics. In the realm of these sample sets, we are interested in estimating $\sigma_i$, with respect to a general affine-invariant loss function $L(t)$, subject to the constraint that $\sigma_{1} \le \sigma_2$. Here, $(X_1^r,Y_1^r,T_1,T_2)$ are sufficient statistics such that $T_{1}=\left(X_{(n)}^r-X_{(1)}^r \right) \sim
 \Gamma(n-1,\sigma_{1})$, $T_{2}=\left(Y_{(m)}^r-Y_{(1)}^r\right)\sim \Gamma(m-1,\sigma_{2})$. Also $X_{1}^r$ and $Y_{1}^r$ follows  exponential distribution $E(\mu_{1},\sigma_{1})$ and $E(\mu_{2},\sigma_{2})$ respectively.
 \section{Numerical comparison and data analysis}\label{sec5}
In this section, we present a numerical comparison of the proposed improved estimators. We have studied the relative risk improvement (RRI) of the proposed estimators with respect to BAEE.  The relative risk improvement of $\delta$ with respect to $\delta_{0}$ is defined as 
$$RRI= \frac{Risk(\delta)-Risk(\delta_{0})}{Risk(\delta_{0})}\times 100$$
To perform the numerical comparison, we have generated 50,000 random samples from an exponential distribution for various values of  the  $(\mu_1,\mu_2)$, as well as $(n_1,n_2)$.  We have plotted the RRI of the improved estimators for various values of  $(\mu_1,\mu_2)$, $(n_1,n_2),$ $(a_1,b_1)$ and $(a_2,b_2)$ as function of $\eta$. 

\begin{figure}[ht]
 			\begin{center}
 				\subfigure[\tiny{$(n_1,n_2)=(8,10) ,(\mu_1,\mu_2)=(0,0) $}]{\includegraphics[height=5.5cm,width=8cm]{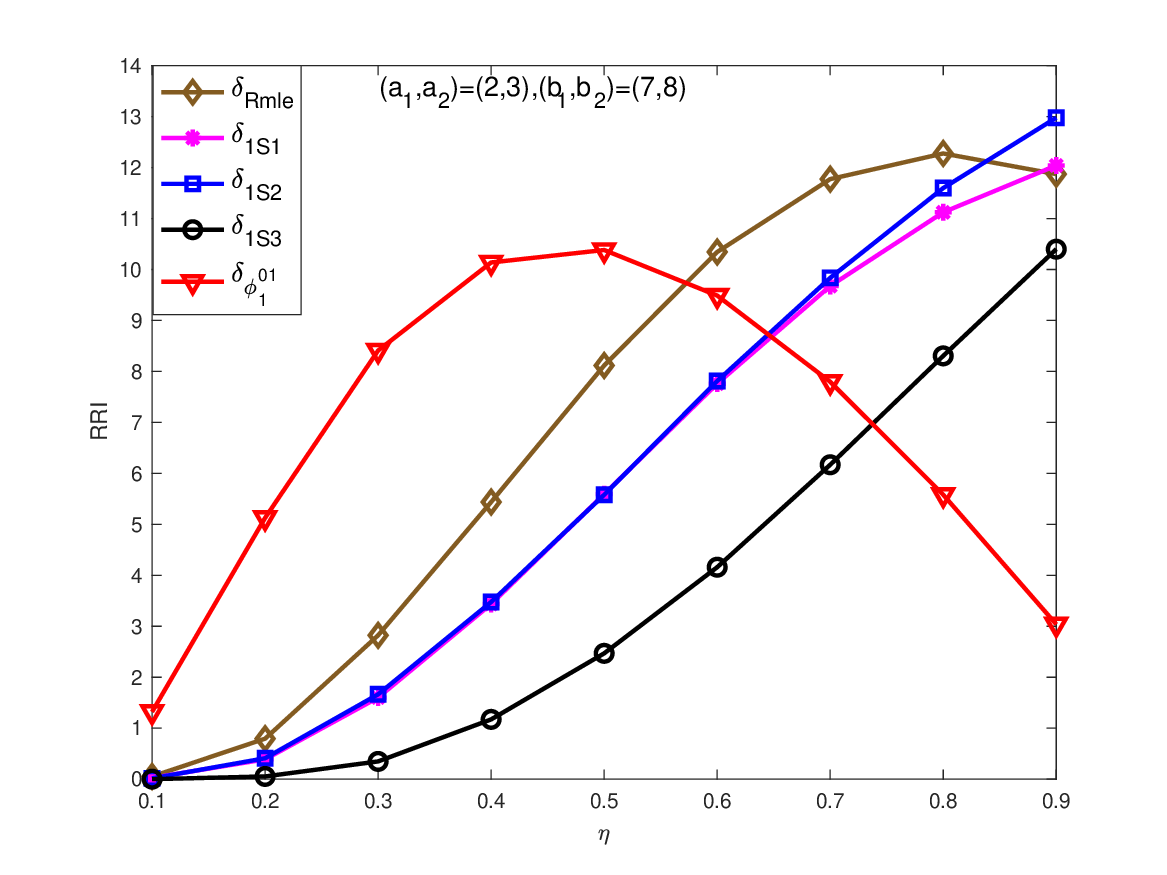}}
 				\hspace{1cm} 
 				\subfigure[\tiny{$(n_1,n_2)=(8,10) ,(\mu_1,\mu_2)=(0,0) $}]{\includegraphics[height=5.5cm,width=8cm]{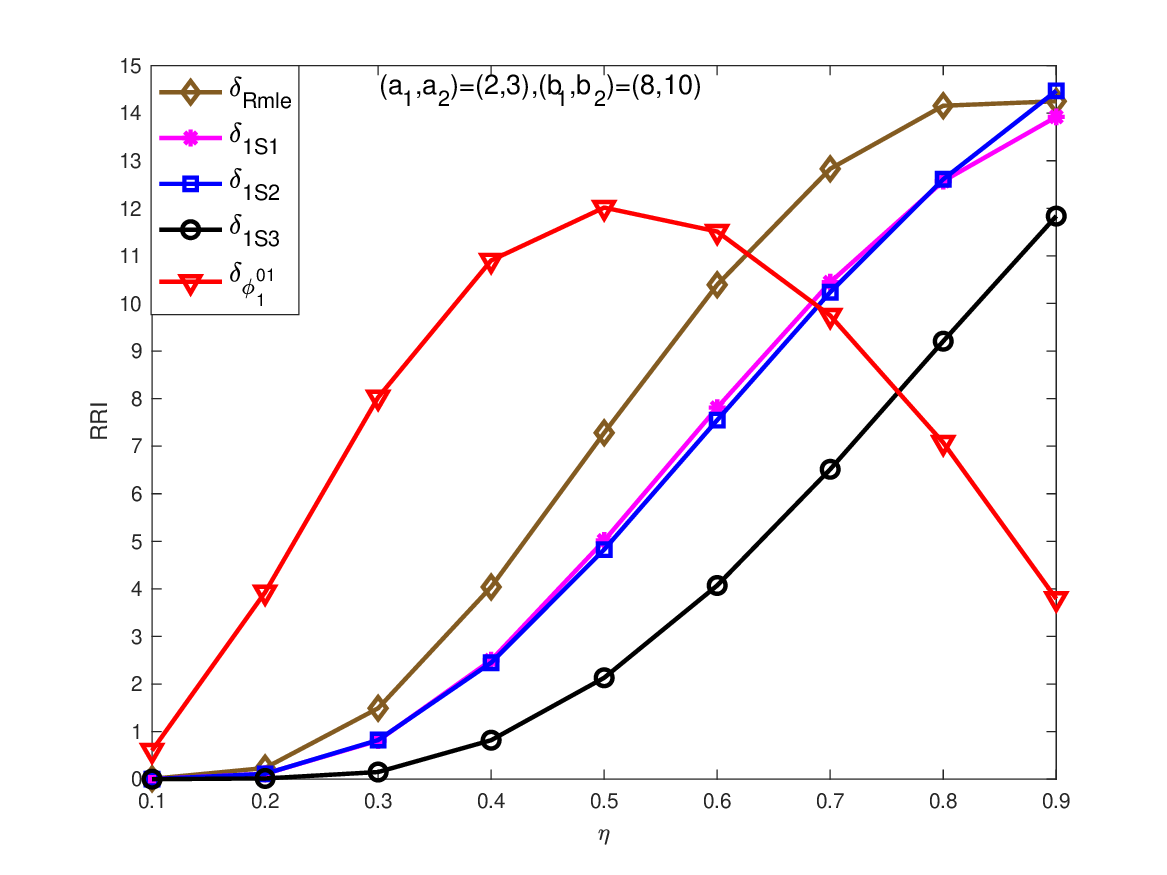}} 
 				\hspace{1cm}
 				\subfigure[\tiny{$(n_1,n_2)=(12,12) ,(\mu_1,\mu_2)=(0.2,0.3) $}]{\includegraphics[height=5.5cm,width=8cm]{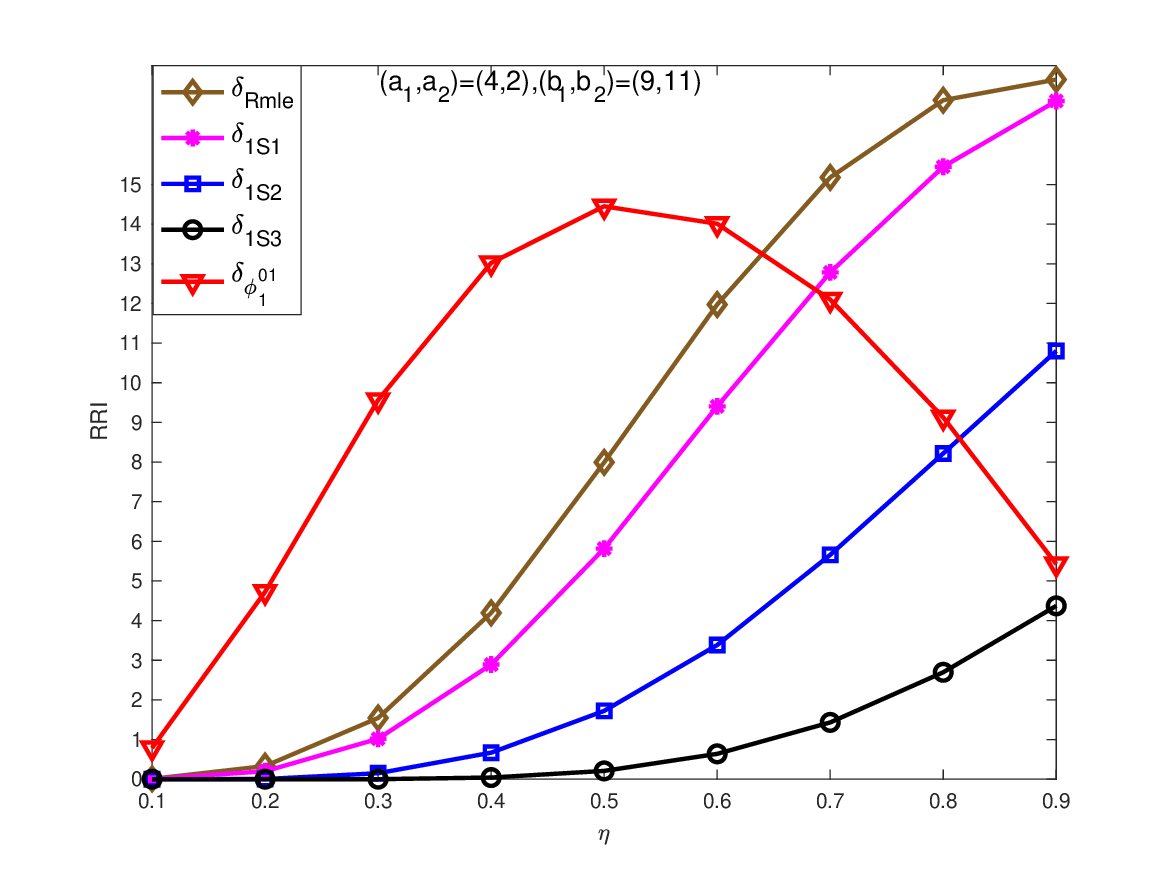}}
 				\hspace{1cm}
 				\subfigure[\tiny{$(n_1,n_2)=(12,12),(\mu_1,\mu_2)=(0.2,0.3)
 					$}]{\includegraphics[height=5.5cm,width=8cm]{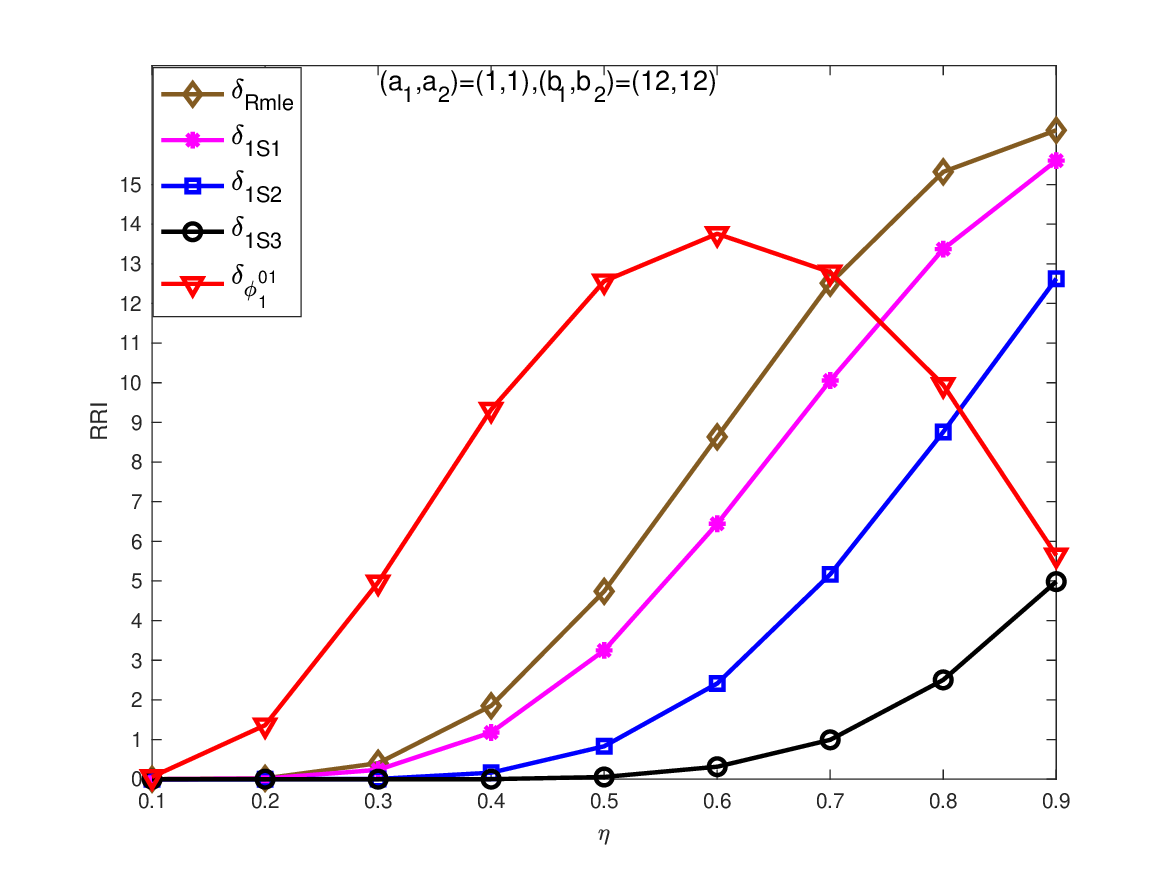}}
 				\subfigure[\tiny{$(n_1,n_2)=(14,15) ,(\mu_1,\mu_2)=(0.4,0.7) $}]{\includegraphics[height=5.5cm,width=8cm]{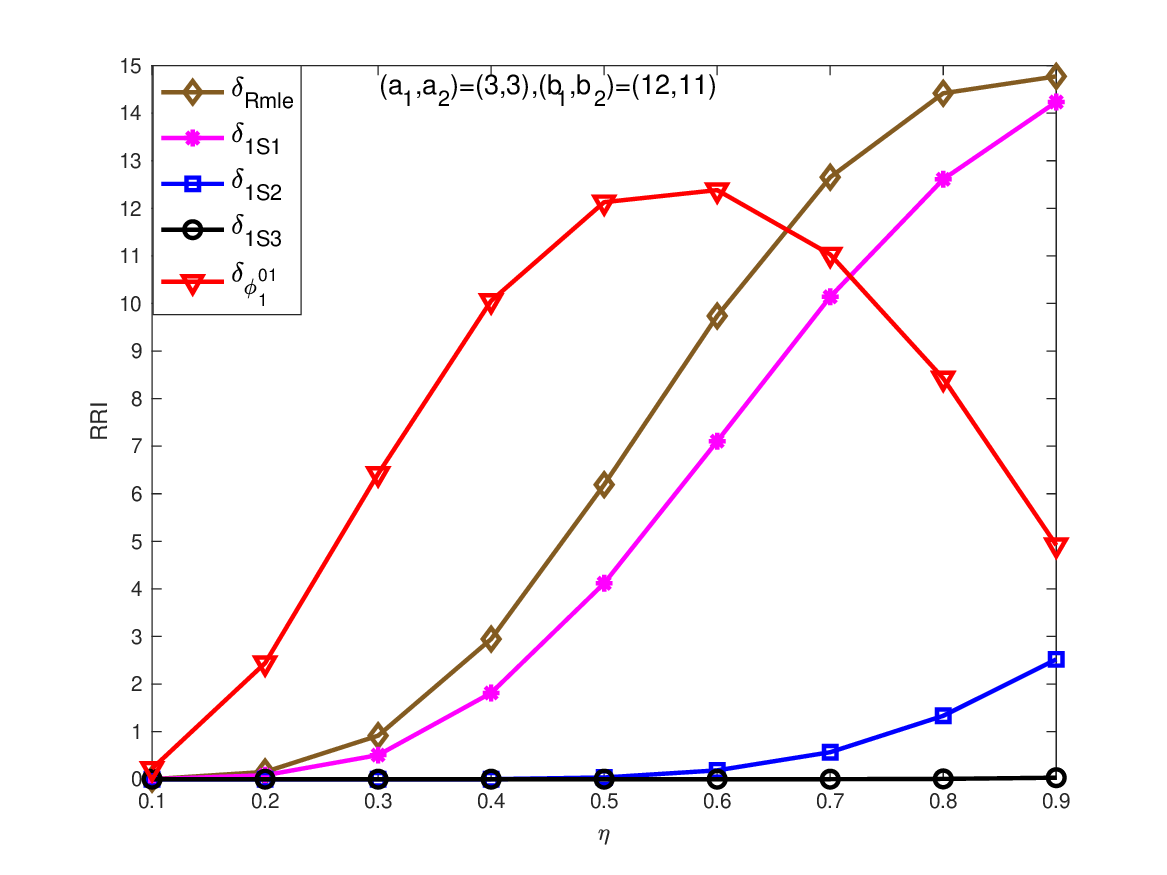}}
 				\hspace{1cm}
 				\subfigure[\tiny{$(n_1,n_2)=(14,15) ,(\mu_1,\mu_2)=(0.4,0.7) $}]{\includegraphics[height=5.5cm,width=8cm]{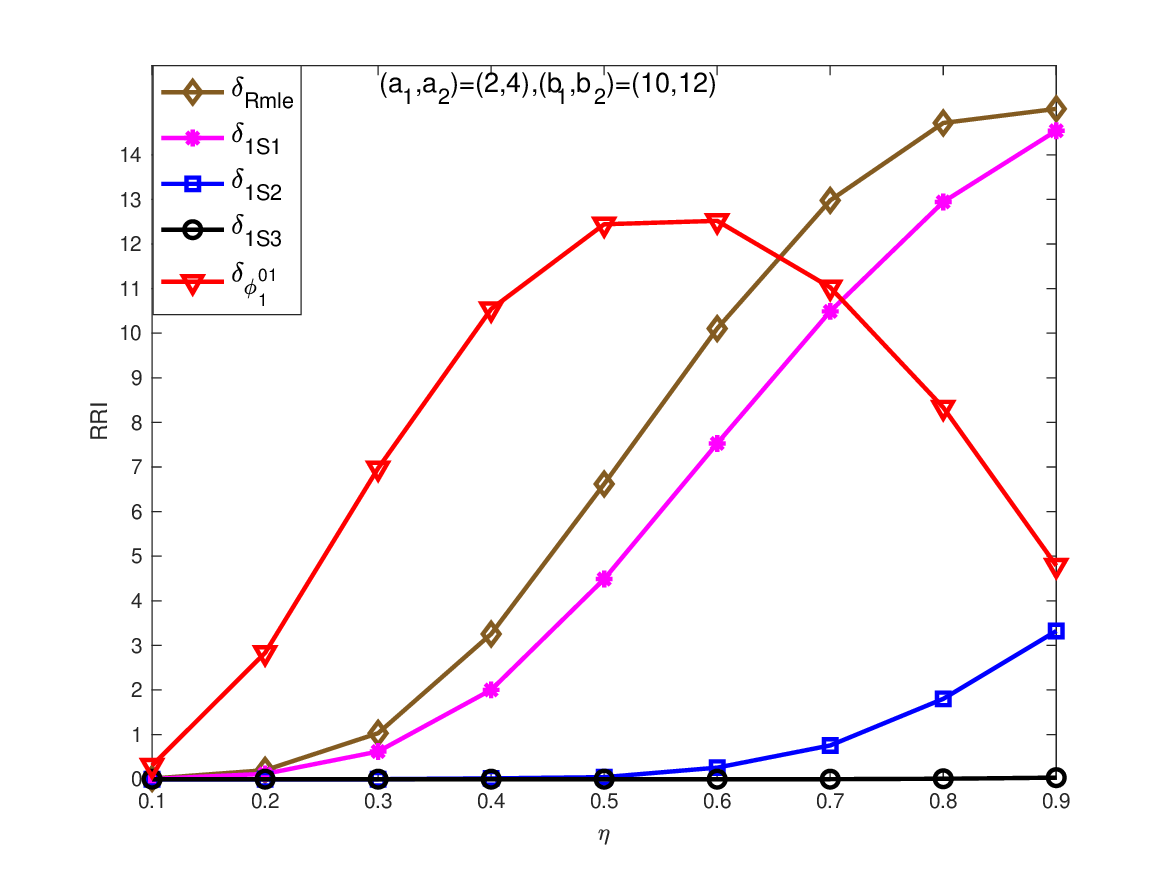}} 
 				\hspace{1cm}			
 					\end{center}
 				\caption{Relative risk improvement under quadratic loss function $L_1(t)$ for $\sigma_{1}$	\label{fig1L1}}
 			\end{figure}

\begin{figure}[ht]\label{fig1}
 				\begin{center}
 						\subfigure[\tiny{$(n_1,n_2)=(13,18) ,(\mu_1,\mu_2)=(-0.1,-0.2) $}]{\includegraphics[height=5.5cm,width=8cm]{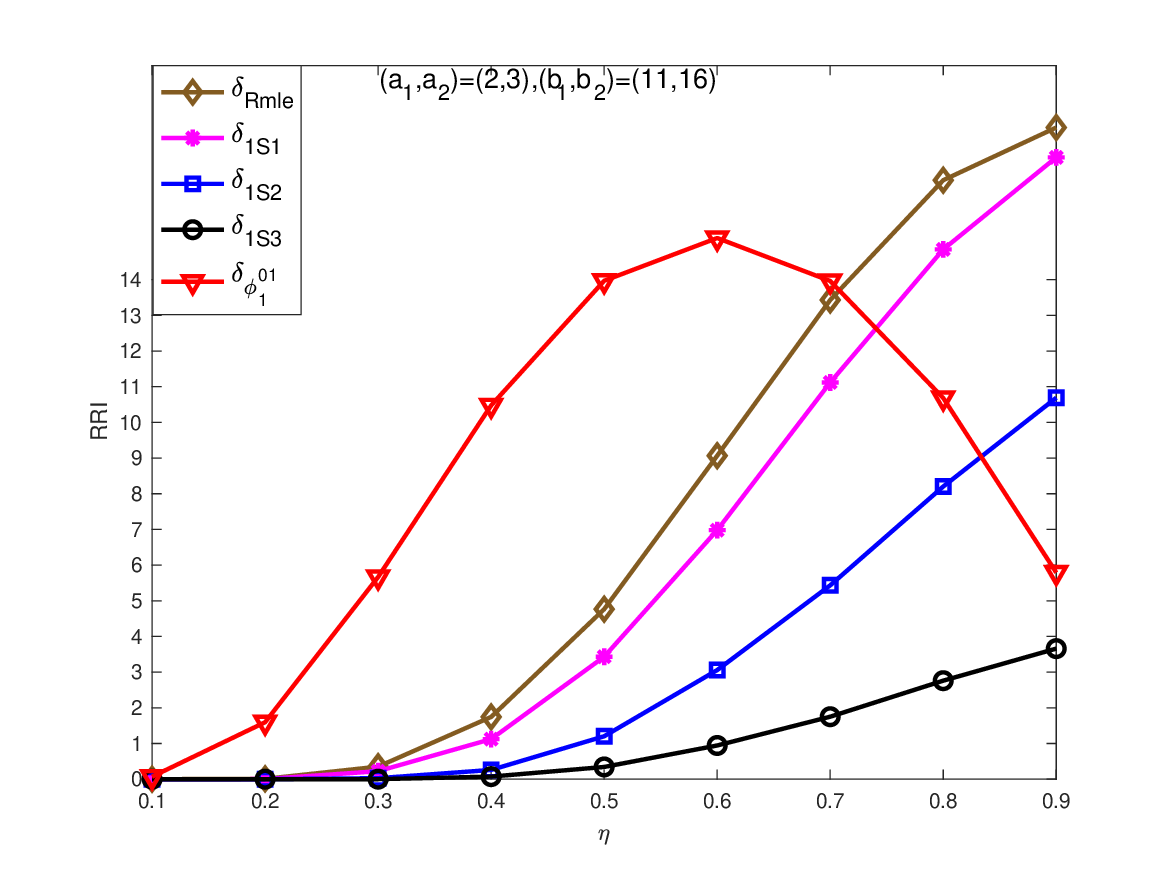}}
 					\hspace{1cm}
 					\subfigure[\tiny{$(n_1,n_2)=(13,18) ,(\mu_1,\mu_2)=(-0.1,-0.2)$}]{\includegraphics[height=5.5cm,width=8cm]{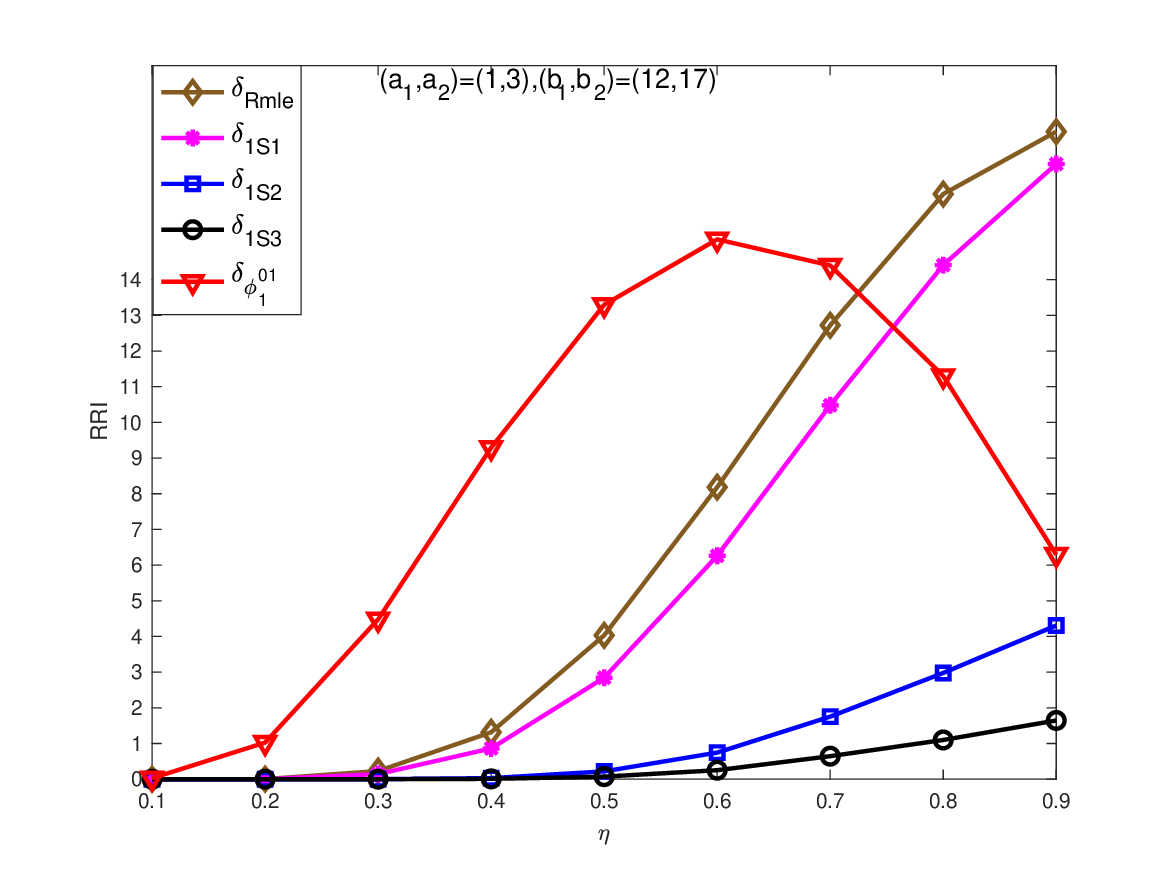}}
 				\subfigure[\tiny{$(n_1,n_2)=(7,8) ,(\mu_1,\mu_2)=(-0.4,-0.5) $}]{\includegraphics[height=5.5cm,width=8cm]{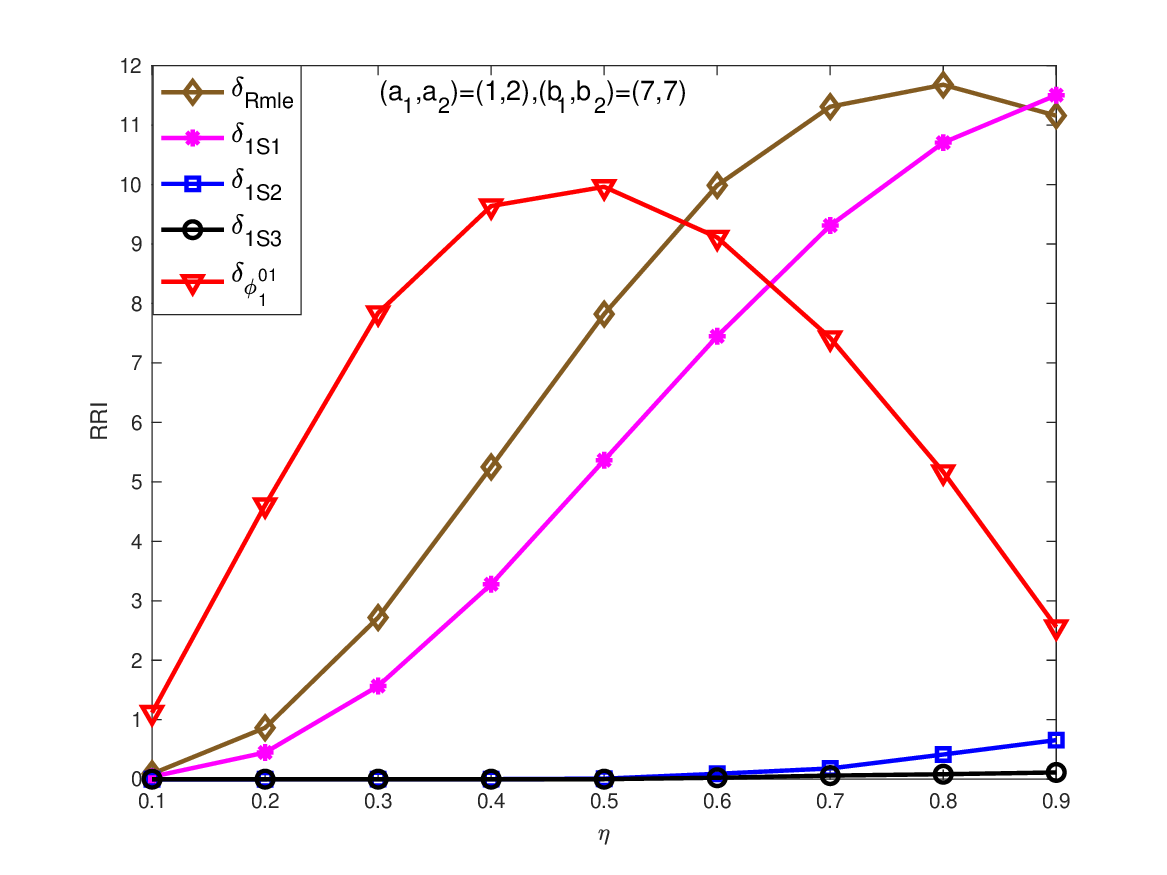}}
 				\hspace{1cm}
 				\subfigure[\tiny{$(n_1,n_2)=(7,8) ,(\mu_1,\mu_2)=(-0.4,-0.5) $}]{\includegraphics[height=5.5cm,width=8cm]{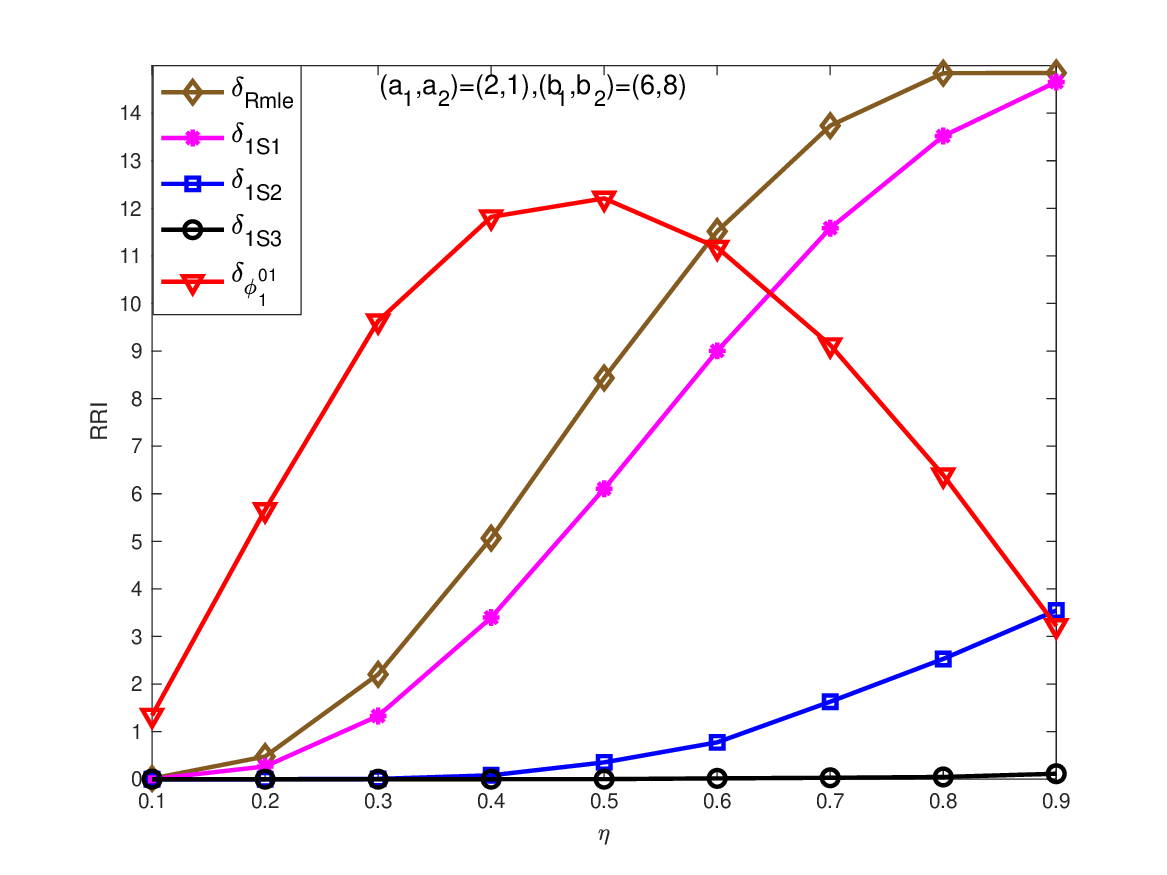}} 
 			\end{center}
 			\caption{Relative risk improvement under quadratic loss function $L_1(t)$ for $\sigma_{1}$	\label{fig2L1}}
 		\end{figure}

\begin{figure}[ht]
	\begin{center}
		\subfigure[\tiny{$(n_1,n_2)=(8,10) ,(\mu_1,\mu_2)=(0,0) $}]{\includegraphics[height=5.5cm,width=8cm]{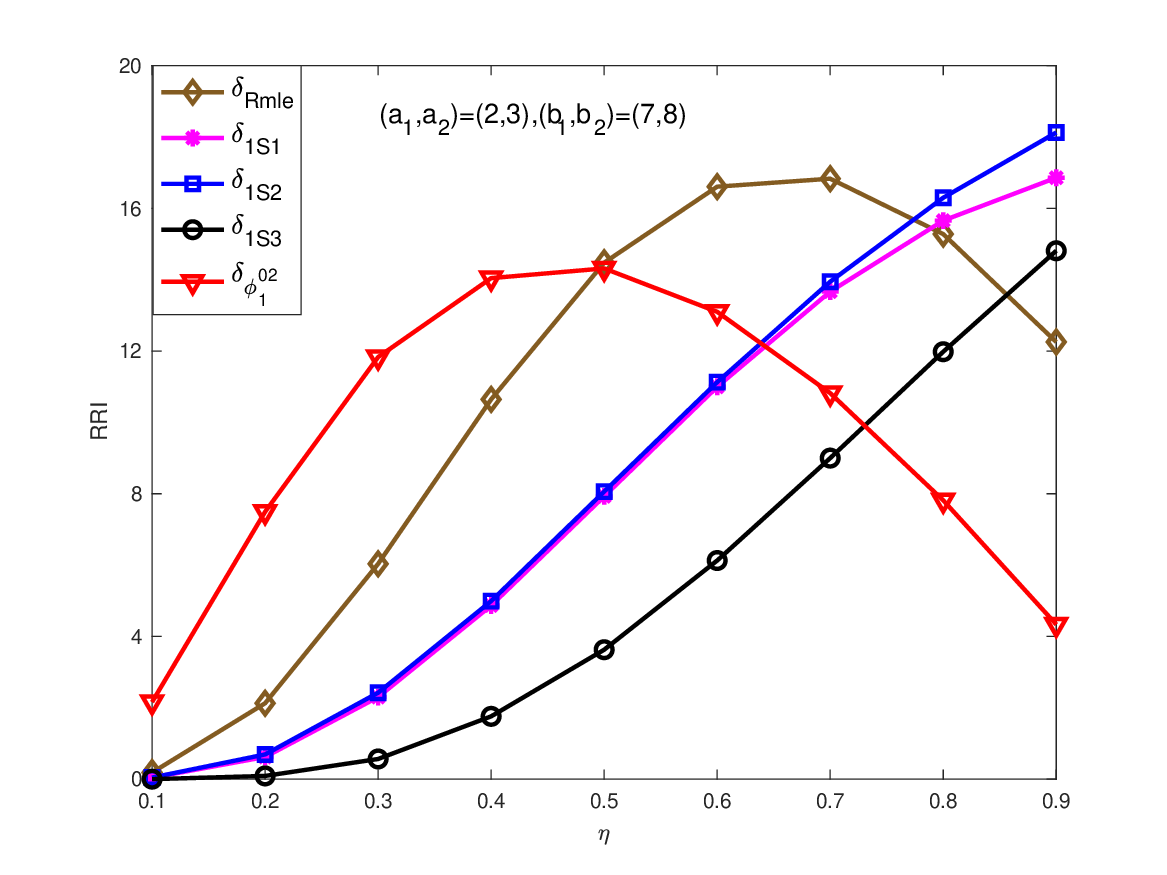}}
		\hspace{1cm} 
		\subfigure[\tiny{$(n_1,n_2)=(8,10) ,(\mu_1,\mu_2)=(0,0) $}]{\includegraphics[height=5.5cm,width=8cm]{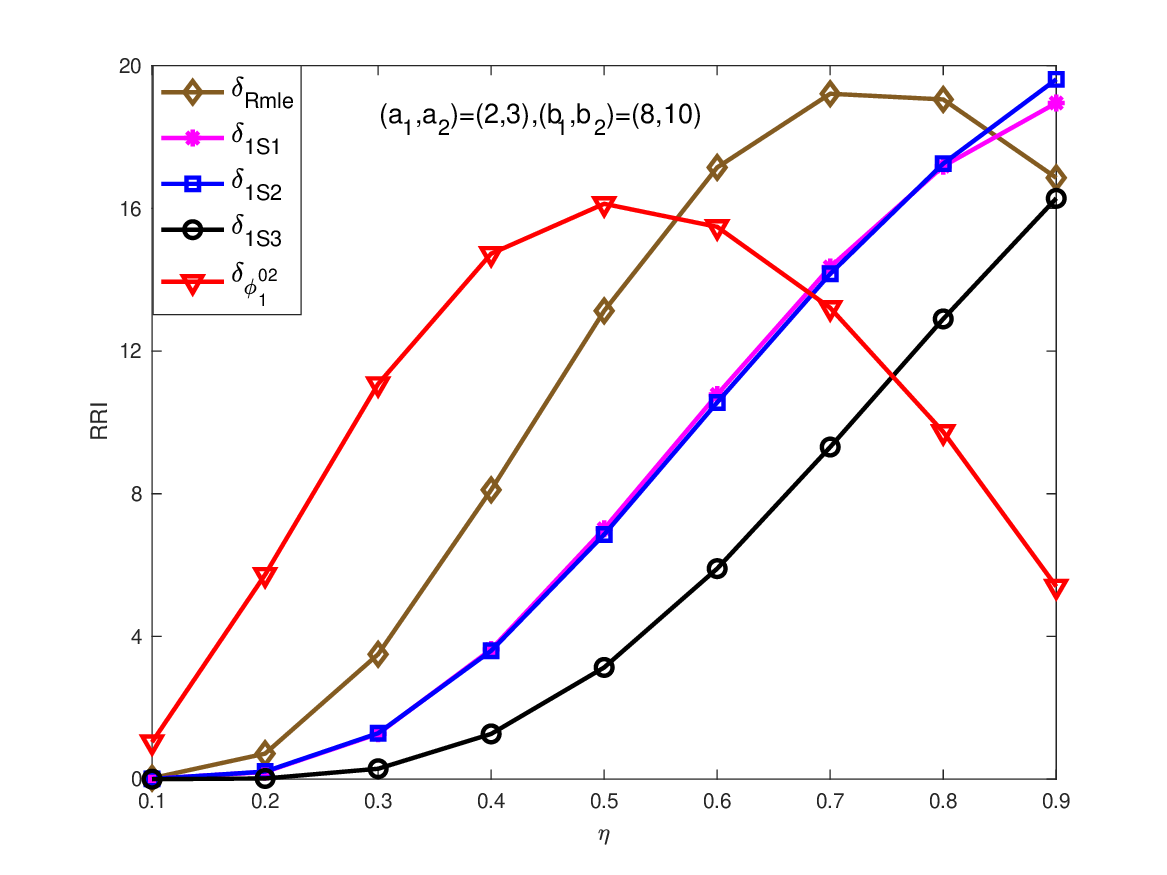}} 
		\hspace{1cm}
			\end{center}
		\caption{Relative risk improvement under entropy loss function $L_2(t)$ for $\sigma_{1}$	\label{fig1L2}}
	\end{figure}

\begin{figure}[ht]
		\begin{center}
		\subfigure[\tiny{$(n_1,n_2)=(12,12) ,(\mu_1,\mu_2)=(0.2,0.3) $}]{\includegraphics[height=5.5cm,width=8cm]{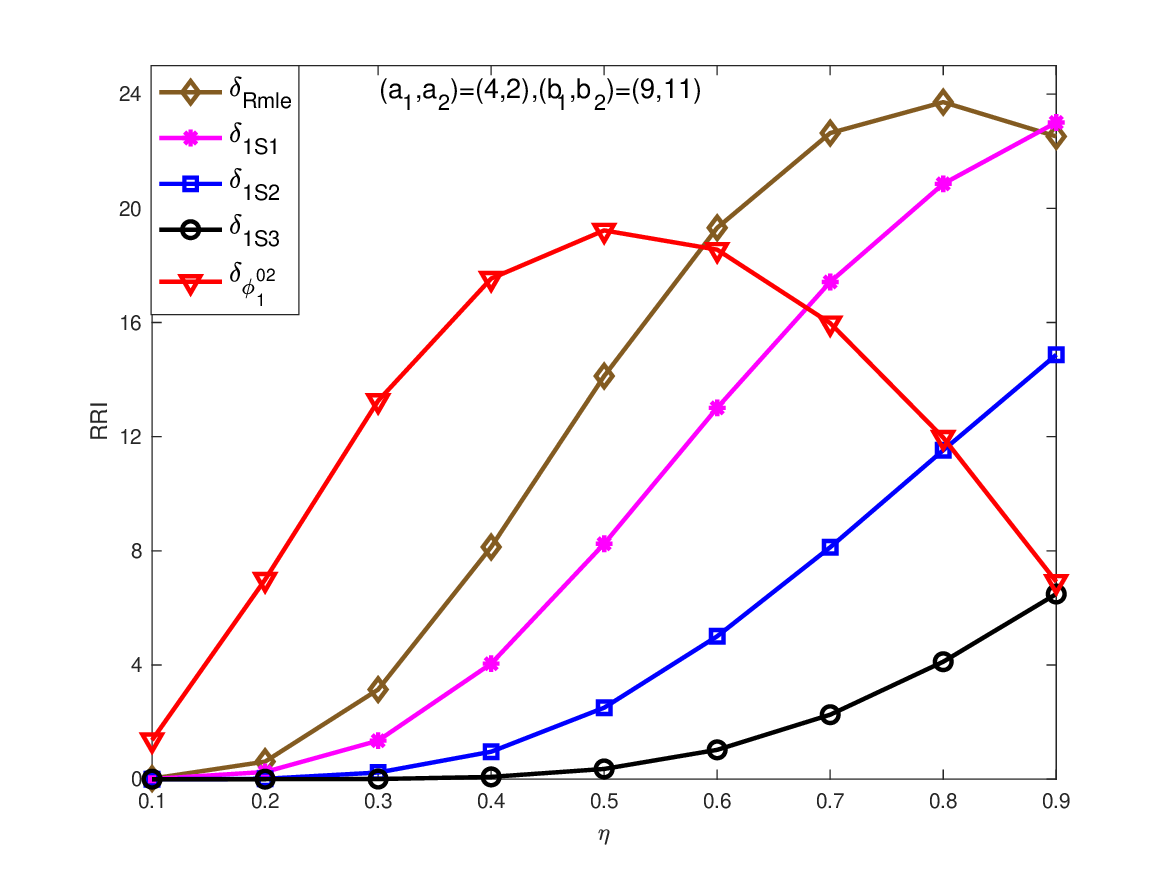}}
		\hspace{1cm}
		\subfigure[\tiny{$(n_1,n_2)=(12,12),(\mu_1,\mu_2)=(0.2,0.3)
			$}]{\includegraphics[height=5.5cm,width=7.5cm]{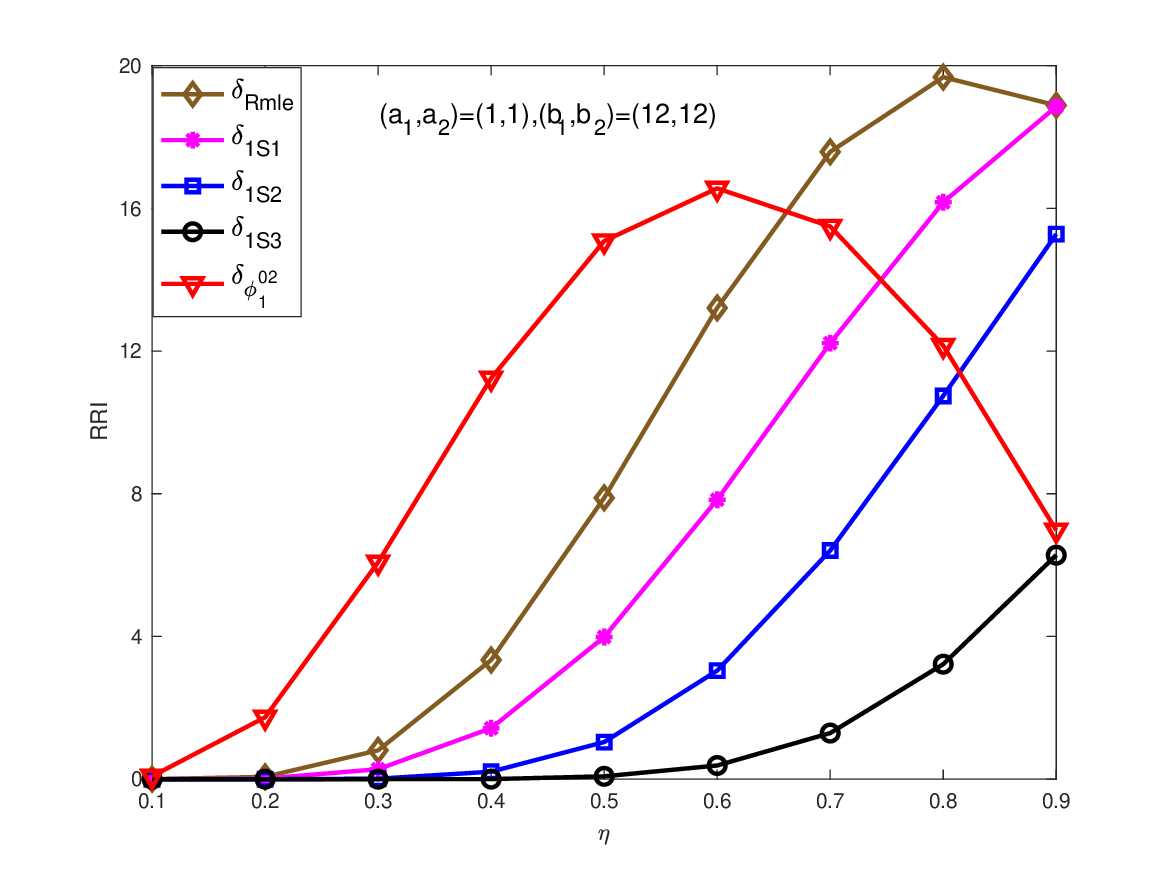}}
		\subfigure[\tiny{$(n_1,n_2)=(14,15) ,(\mu_1,\mu_2)=(0.4,0.7) $}]{\includegraphics[height=5.5cm,width=8cm]{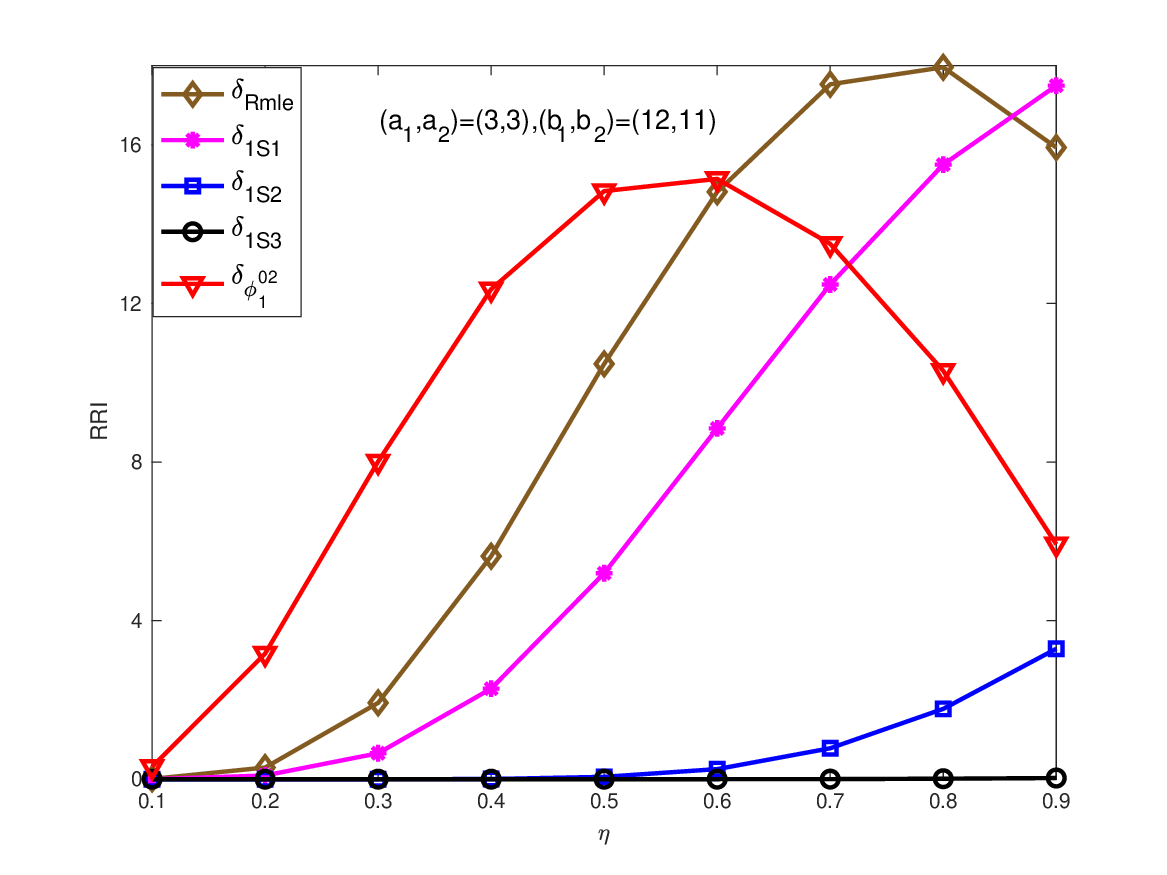}}
		\hspace{1cm}
		\subfigure[\tiny{$(n_1,n_2)=(14,15) ,(\mu_1,\mu_2)=(0.4,0.7) $}]{\includegraphics[height=5.5cm,width=8cm]{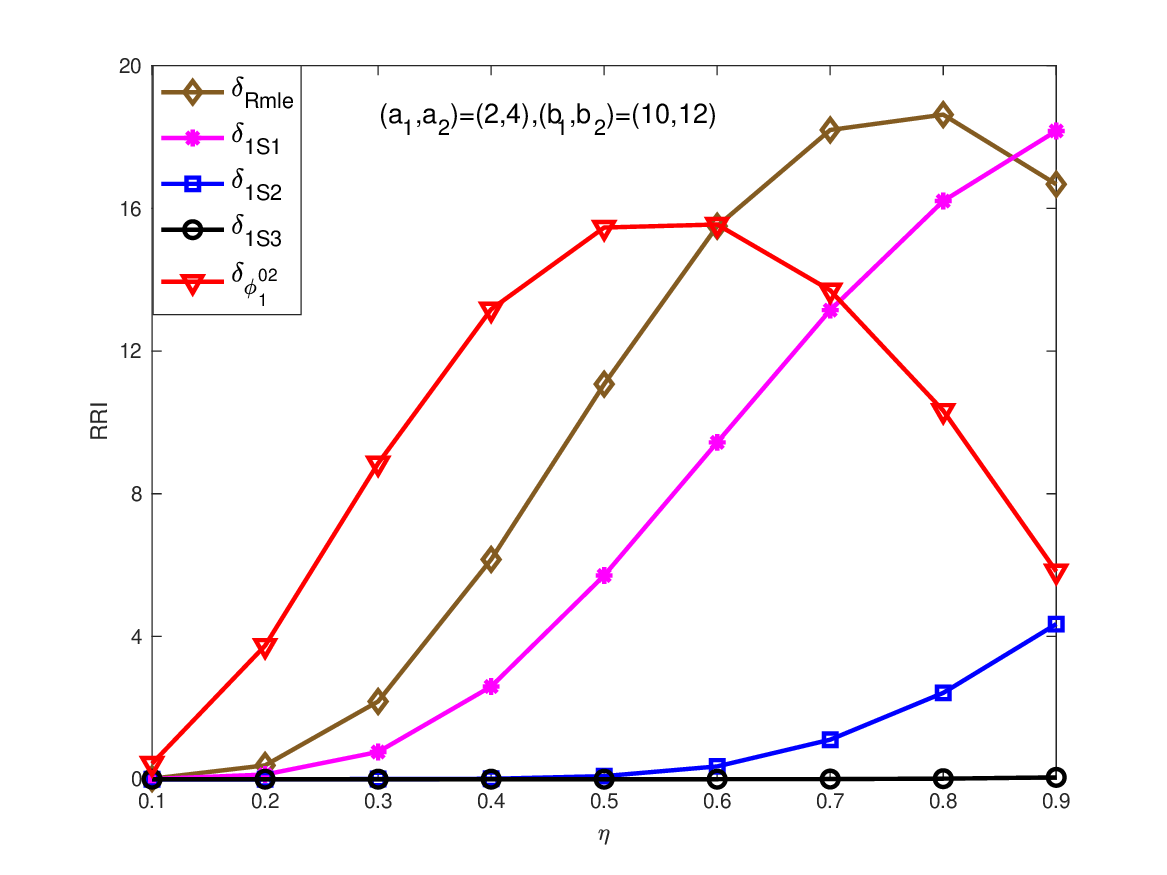}} 
		\hspace{1cm}
		\subfigure[\tiny{$(n_1,n_2)=(13,18) ,(\mu_1,\mu_2)=(-0.1,-0.2) $}]{\includegraphics[height=5.5cm,width=7.5cm]{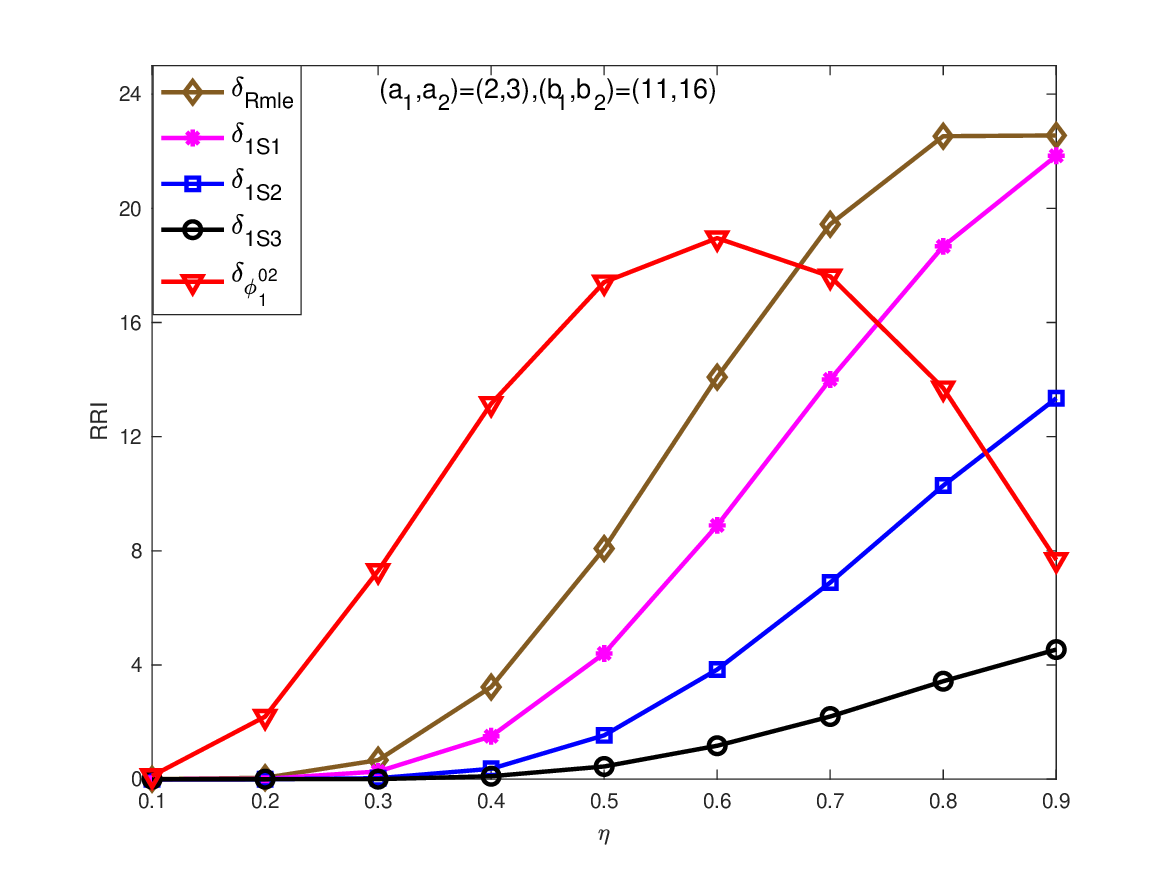}}
		\hspace{1cm}
		\subfigure[\tiny{$(n_1,n_2)=(13,18) ,(\mu_1,\mu_2)=(-0.1,-0.2)$}]{\includegraphics[height=5.5cm,width=8cm]{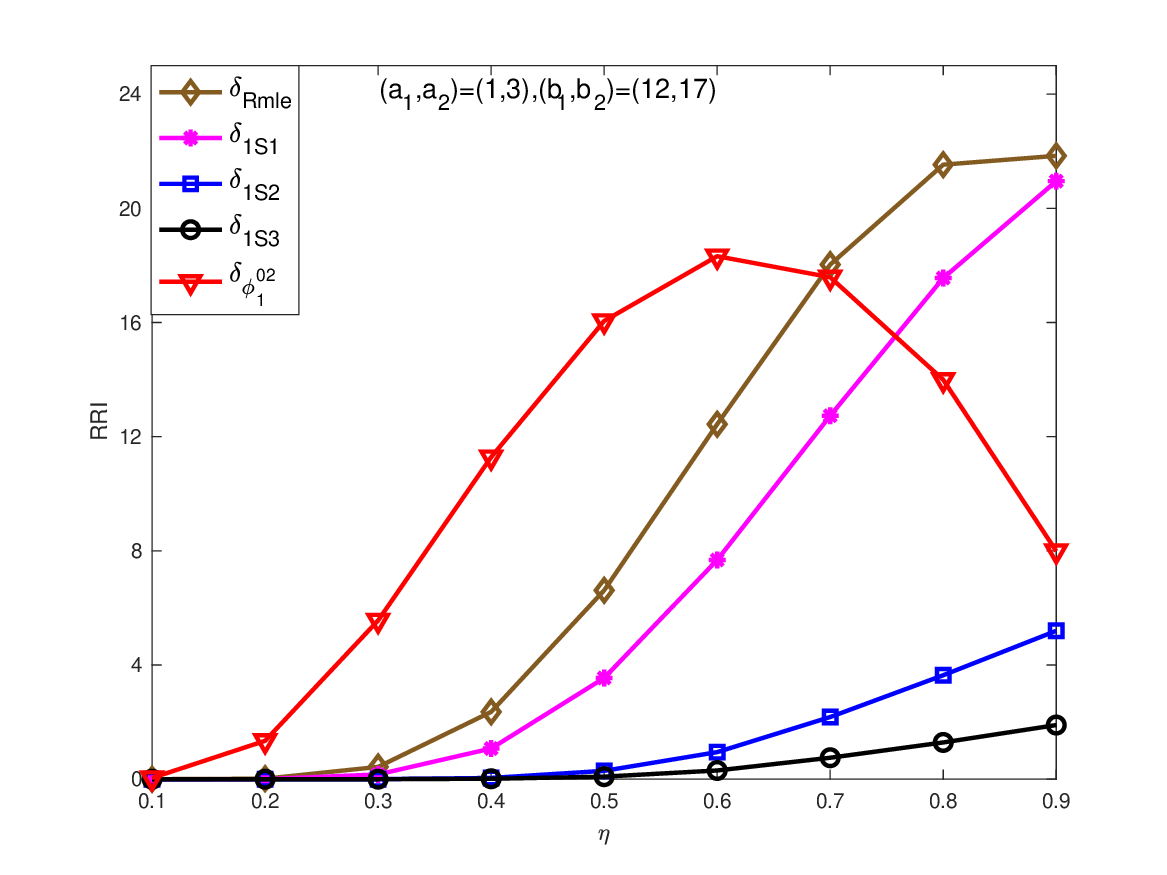}}
			\end{center}
		\caption{Relative risk improvement under entropy loss function $L_2(t)$ for $\sigma_{1}$ \label{fig2L2}}
	\end{figure}

\begin{figure}[ht]
		\begin{center}
		\subfigure[\tiny{$(n_1,n_2)=(7,8) ,(\mu_1,\mu_2)=(-0.4,-0.5) $}]{\includegraphics[height=5.5cm,width=8cm]{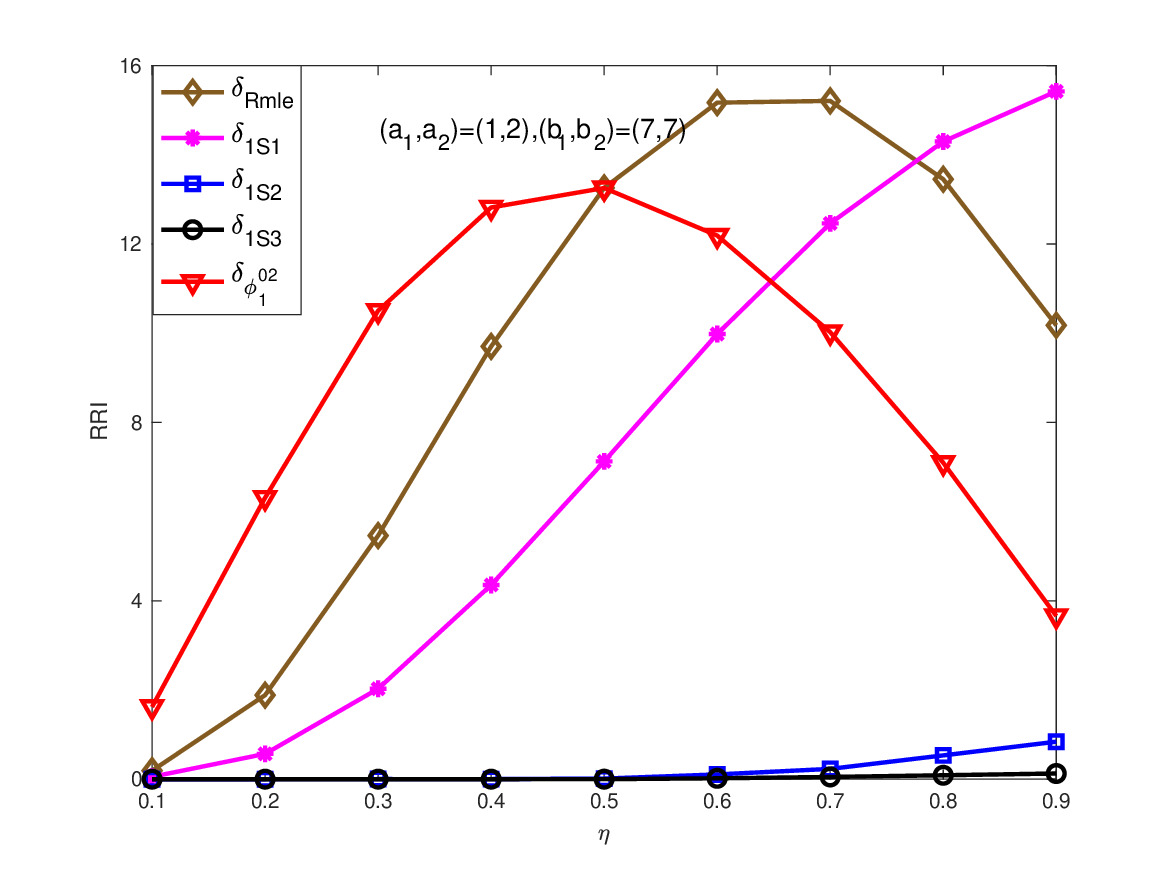}}
		\hspace{1cm}
		\subfigure[\tiny{$(n_1,n_2)=(7,8) ,(\mu_1,\mu_2)=(-0.4,-0.5) $}]{\includegraphics[height=5.5cm,width=8cm]{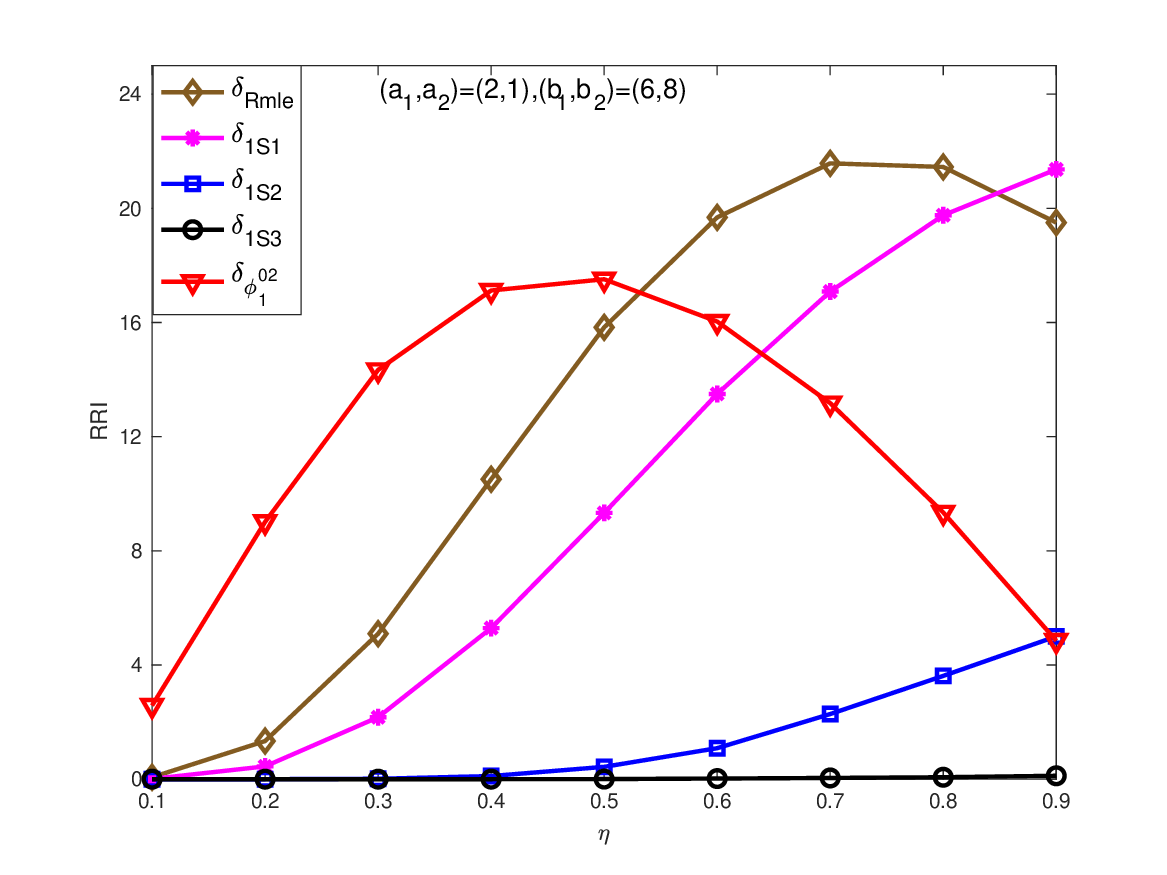}} 
	\end{center}
	\caption{Relative risk improvement under entropy loss function $L_2(t)$ for $\sigma_{1}$ \label{fig3L2}}
\end{figure}

\begin{figure}[h]
	\begin{center}
		\subfigure[\tiny{$(n_1,n_2)=(8,10) ,(\mu_1,\mu_2)=(0,0) $}]{\includegraphics[height=5.5cm,width=8cm]{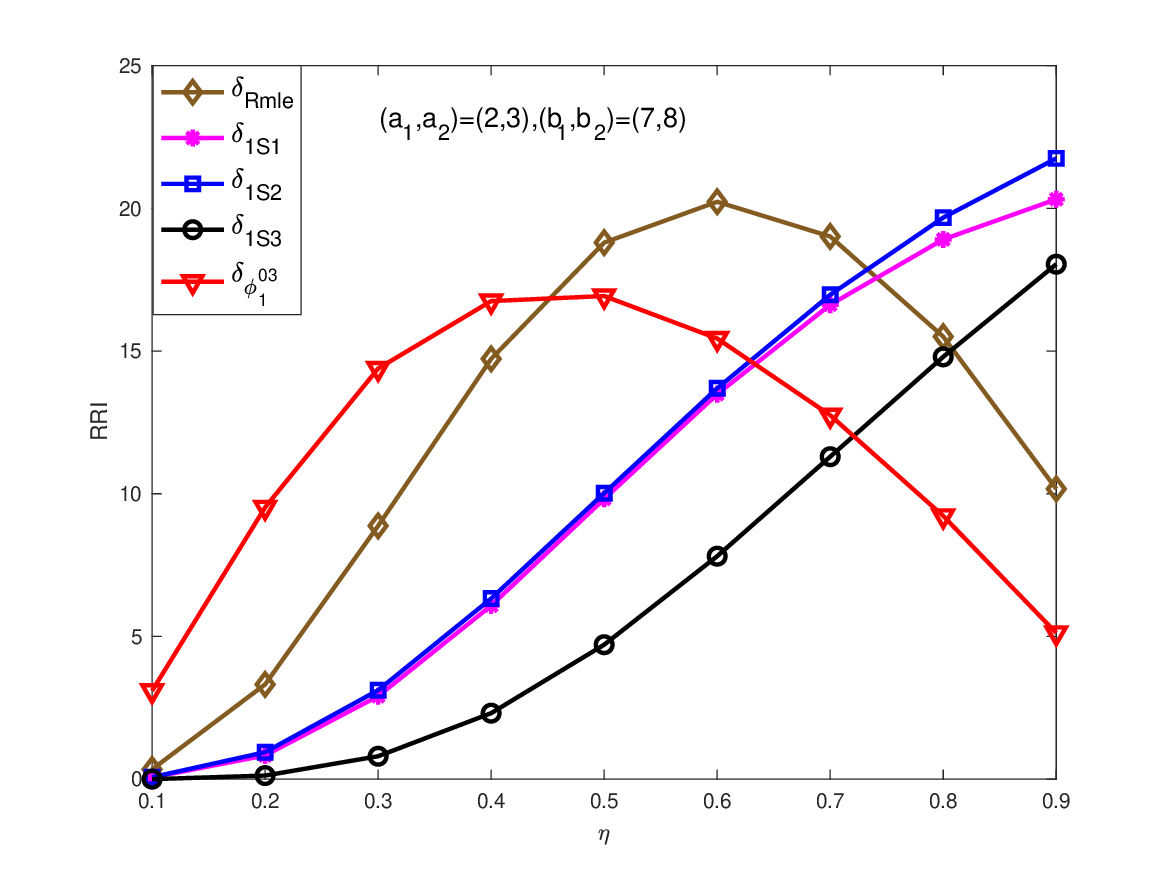}}
		\hspace{1cm} 
		\subfigure[\tiny{$(n_1,n_2)=(8,10) ,(\mu_1,\mu_2)=(0,0) $}]{\includegraphics[height=5.5cm,width=8cm]{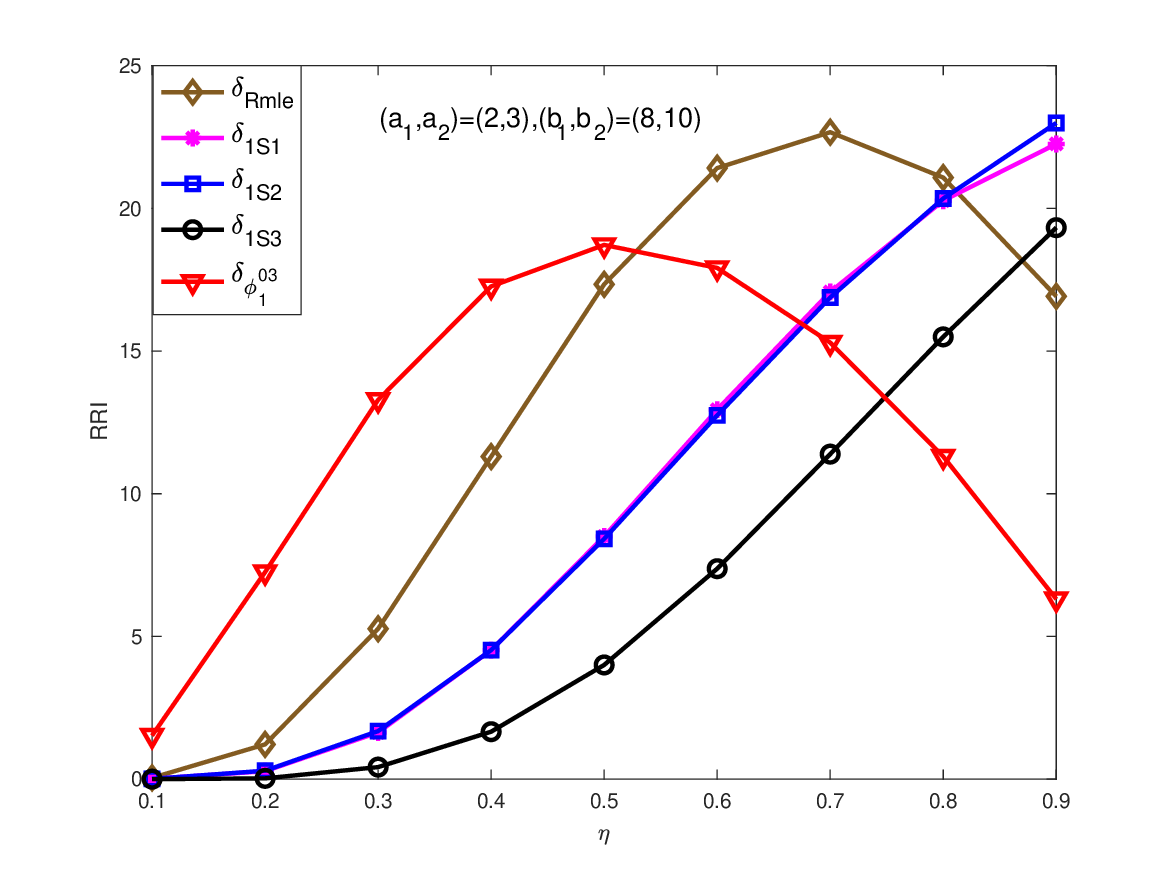}} 
		\hspace{1cm}
		\subfigure[\tiny{$(n_1,n_2)=(12,12) ,(\mu_1,\mu_2)=(0.2,0.3) $}]{\includegraphics[height=5.5cm,width=8cm]{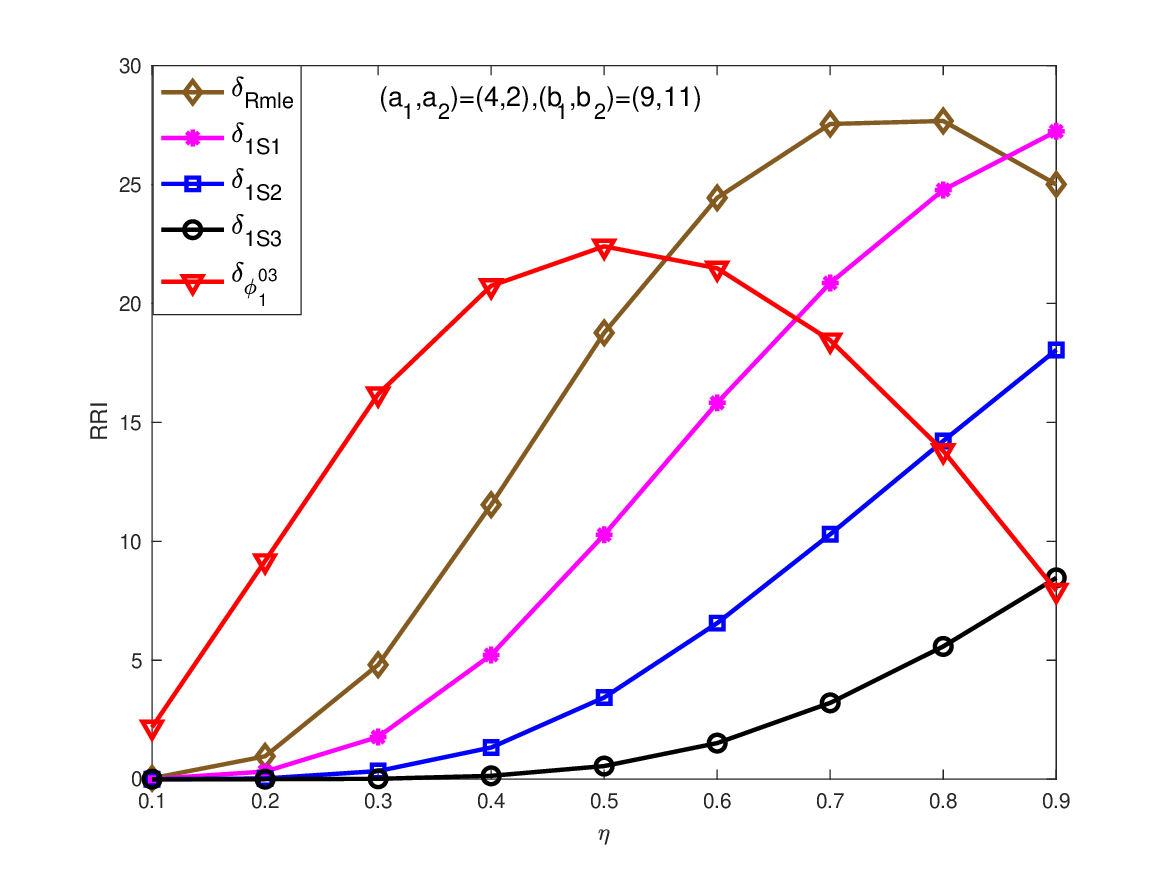}}
		\hspace{1cm}
		\subfigure[\tiny{$(n_1,n_2)=(12,12),(\mu_1,\mu_2)=(0.2,0.3)
			$}]{\includegraphics[height=5.5cm,width=8cm]{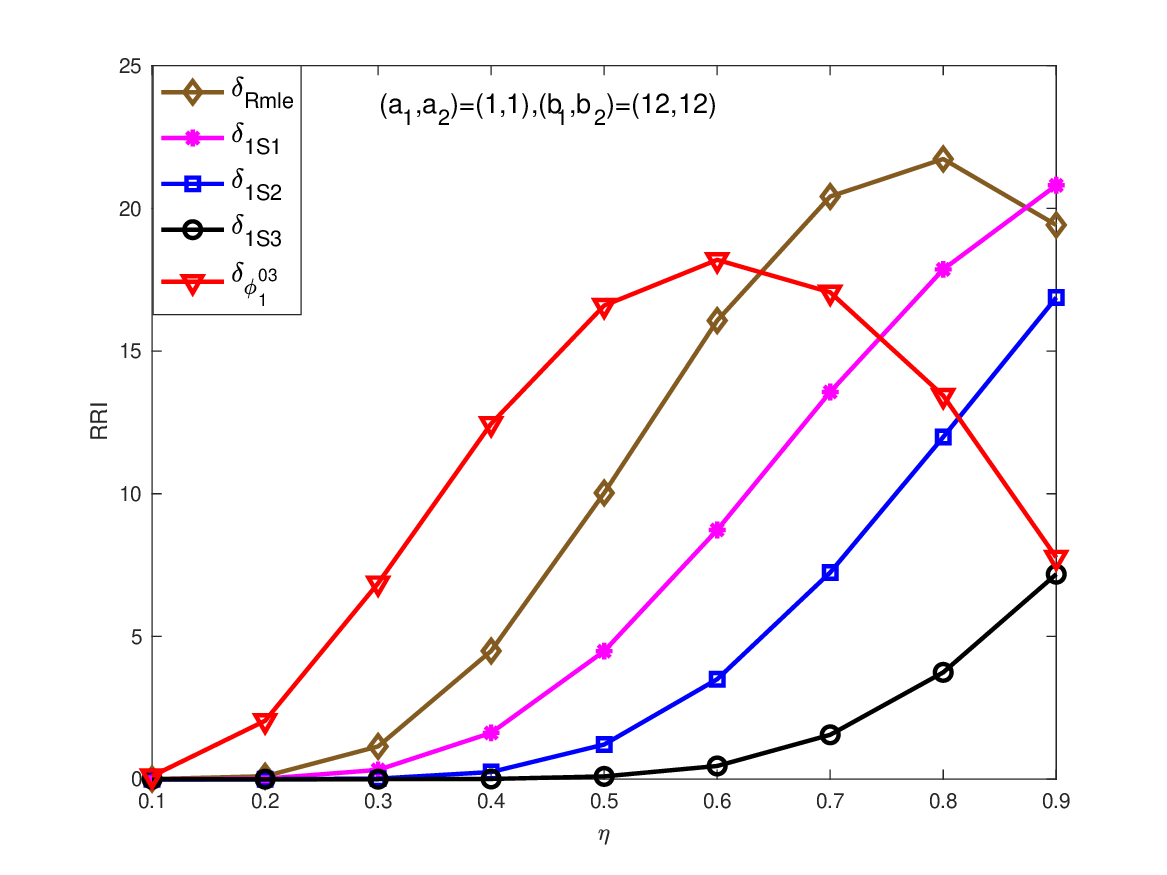}}
		\end{center}
		\caption{Relative risk improvement under symmetric loss function $L_3(t)$ for $\sigma_{1}$\label{fig1L3}}
	\end{figure}

\begin{figure}[h]
		\begin{center}
		\subfigure[\tiny{$(n_1,n_2)=(14,15) ,(\mu_1,\mu_2)=(0.4,0.7) $}]{\includegraphics[height=5.5cm,width=8cm]{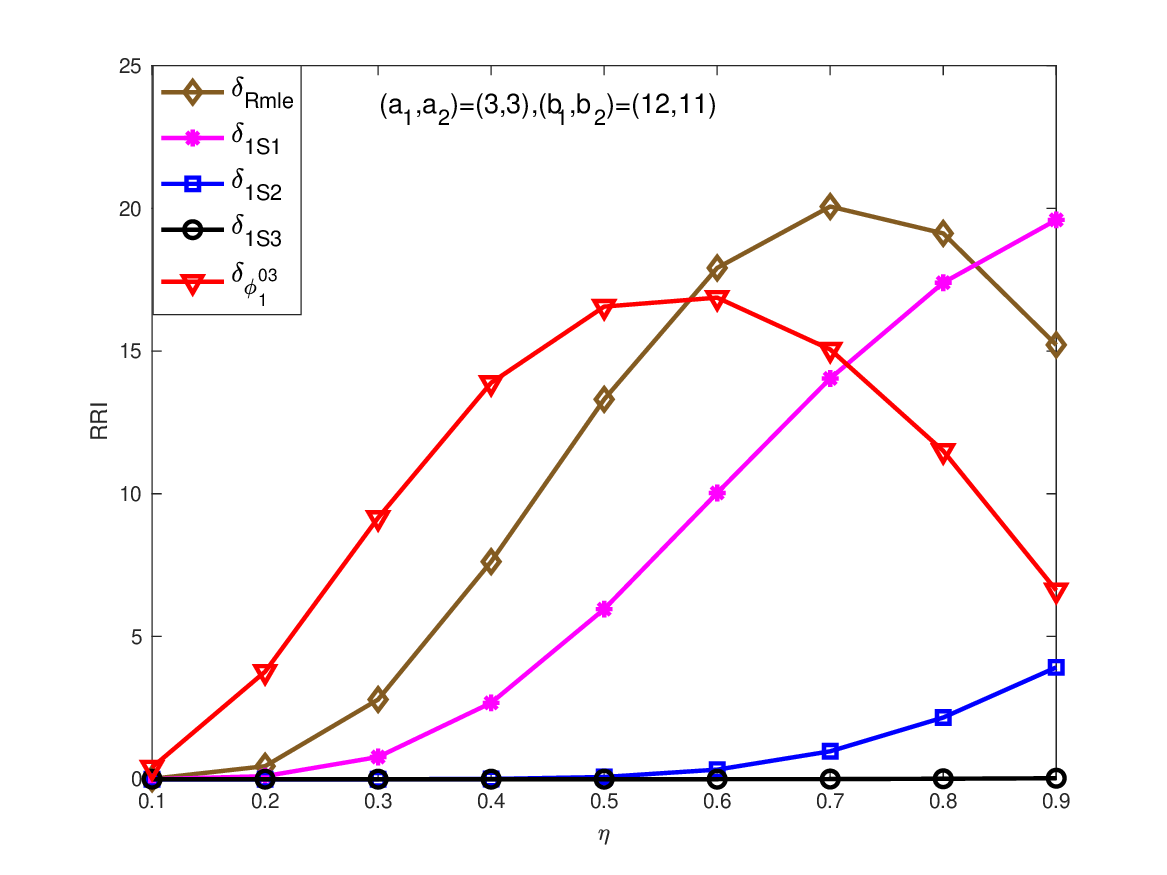}}
		\hspace{1cm}
		\subfigure[\tiny{$(n_1,n_2)=(14,15) ,(\mu_1,\mu_2)=(0.4,0.7) $}]{\includegraphics[height=5.5cm,width=8cm]{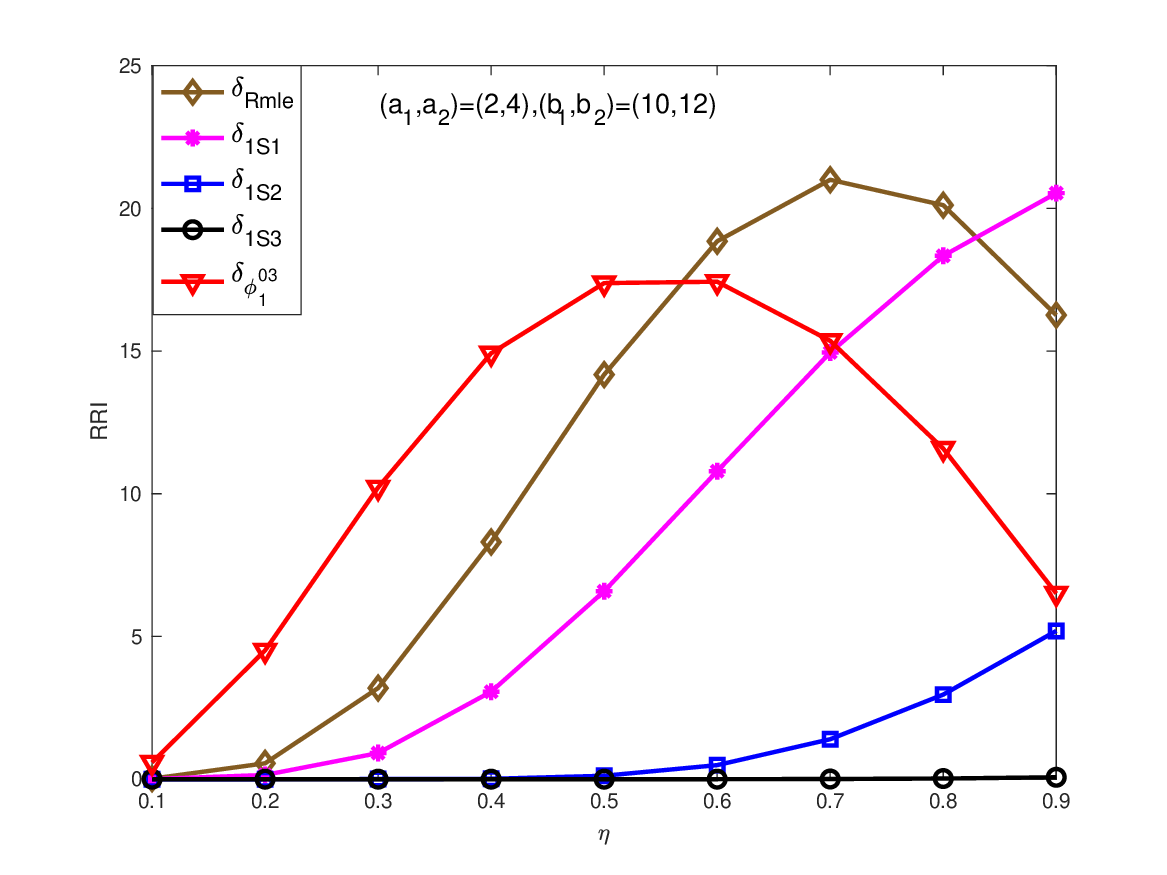}} 
		\hspace{1cm}
		\subfigure[\tiny{$(n_1,n_2)=(13,18) ,(\mu_1,\mu_2)=(-0.1,-0.2) $}]{\includegraphics[height=5.5cm,width=8cm]{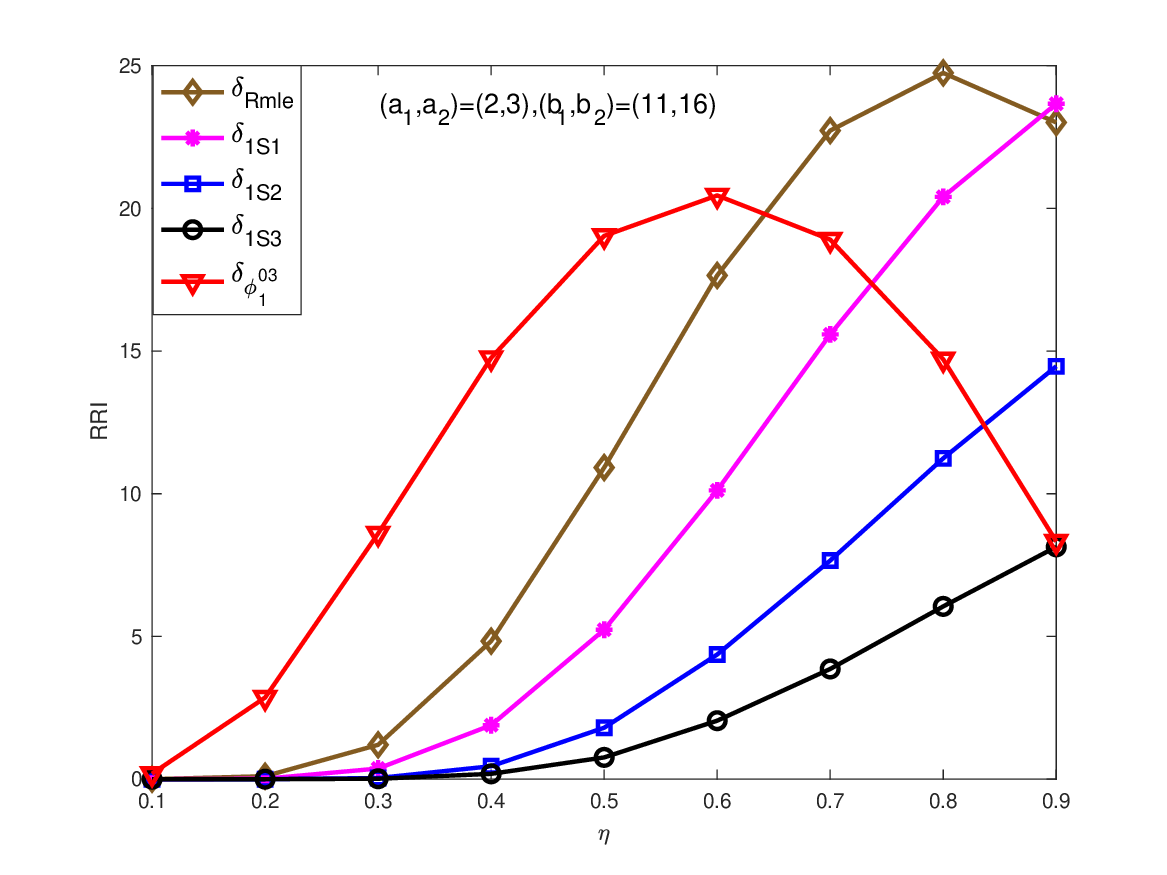}}
		\hspace{1cm}
		\subfigure[\tiny{$(n_1,n_2)=(13,18) ,(\mu_1,\mu_2)=(-0.1,-0.2)$}]{\includegraphics[height=5.5cm,width=8cm]{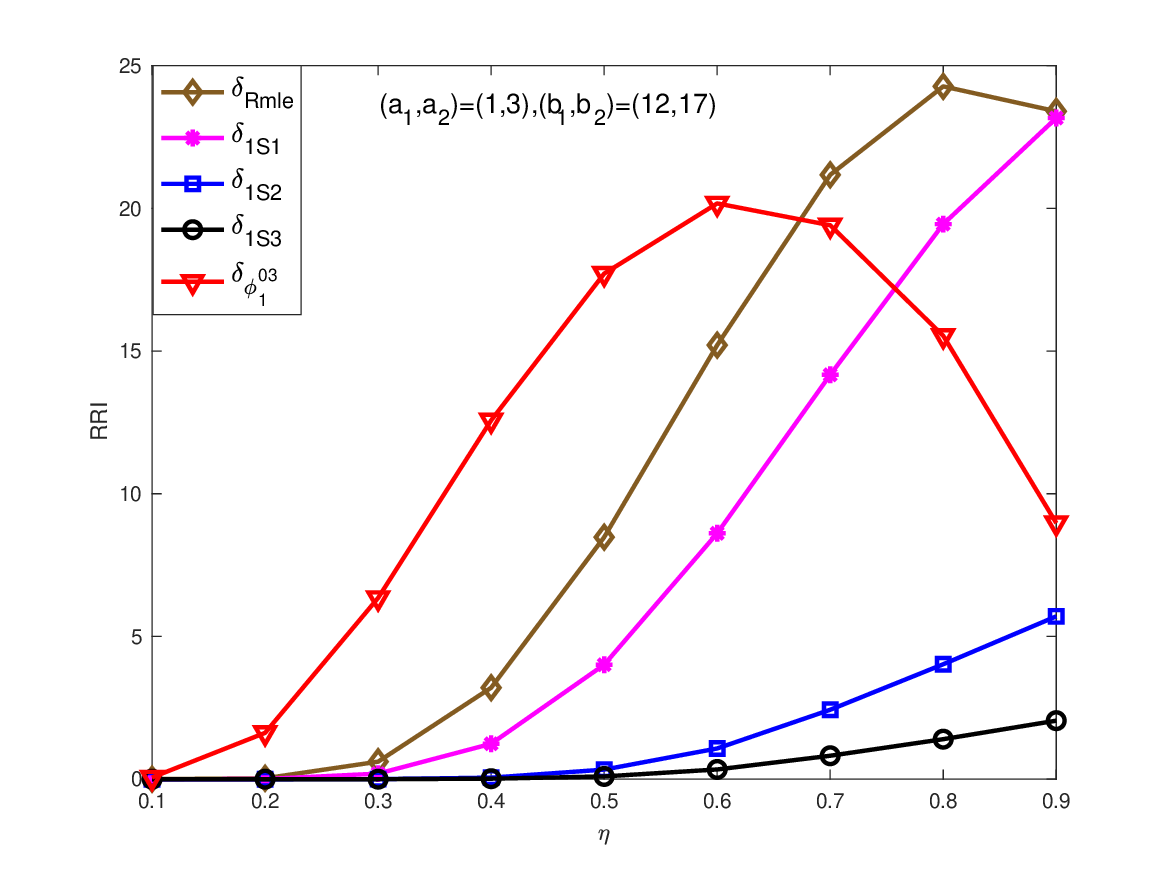}}
		\subfigure[\tiny{$(n_1,n_2)=(7,8) ,(\mu_1,\mu_2)=(-0.4,-0.5) $}]{\includegraphics[height=5.5cm,width=8cm]{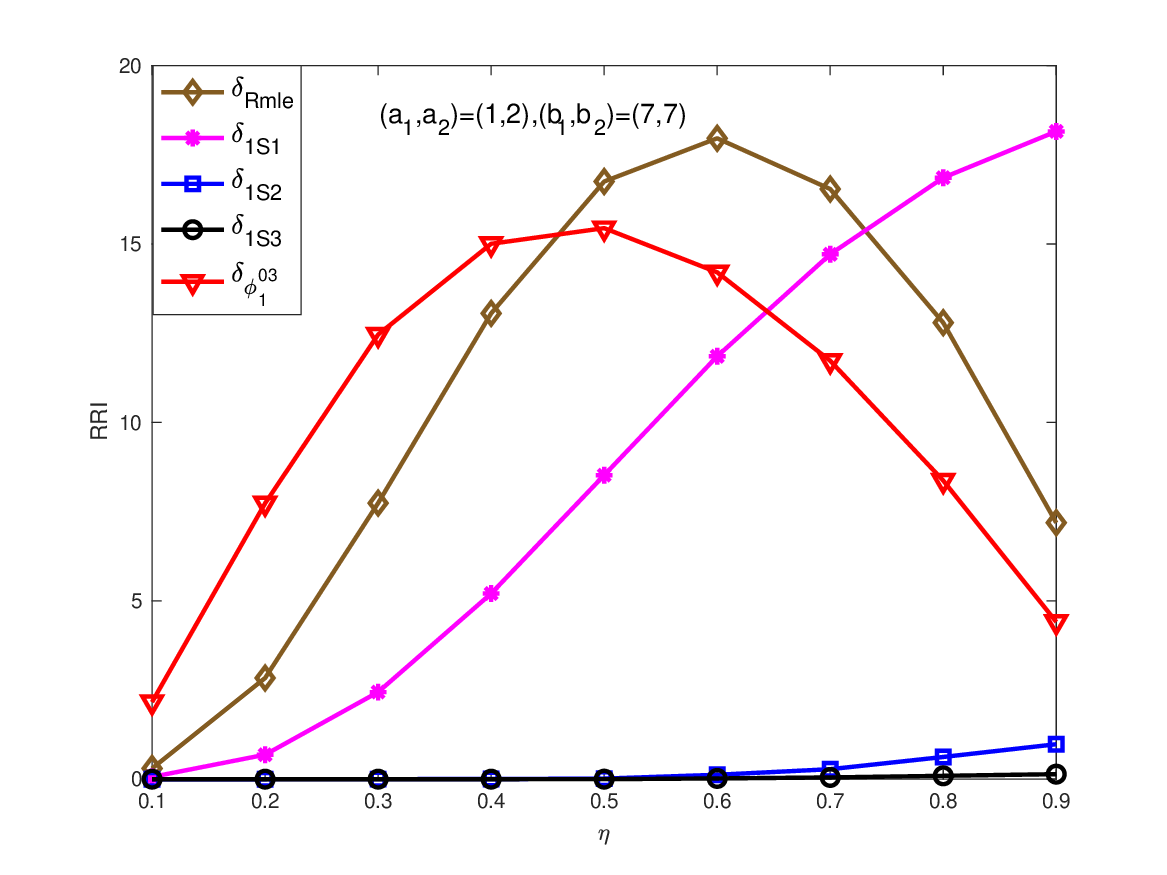}}
		\hspace{1cm}
		\subfigure[\tiny{$(n_1,n_2)=(7,8) ,(\mu_1,\mu_2)=(-0.4,-0.5) $}]{\includegraphics[height=5.5cm,width=8cm]{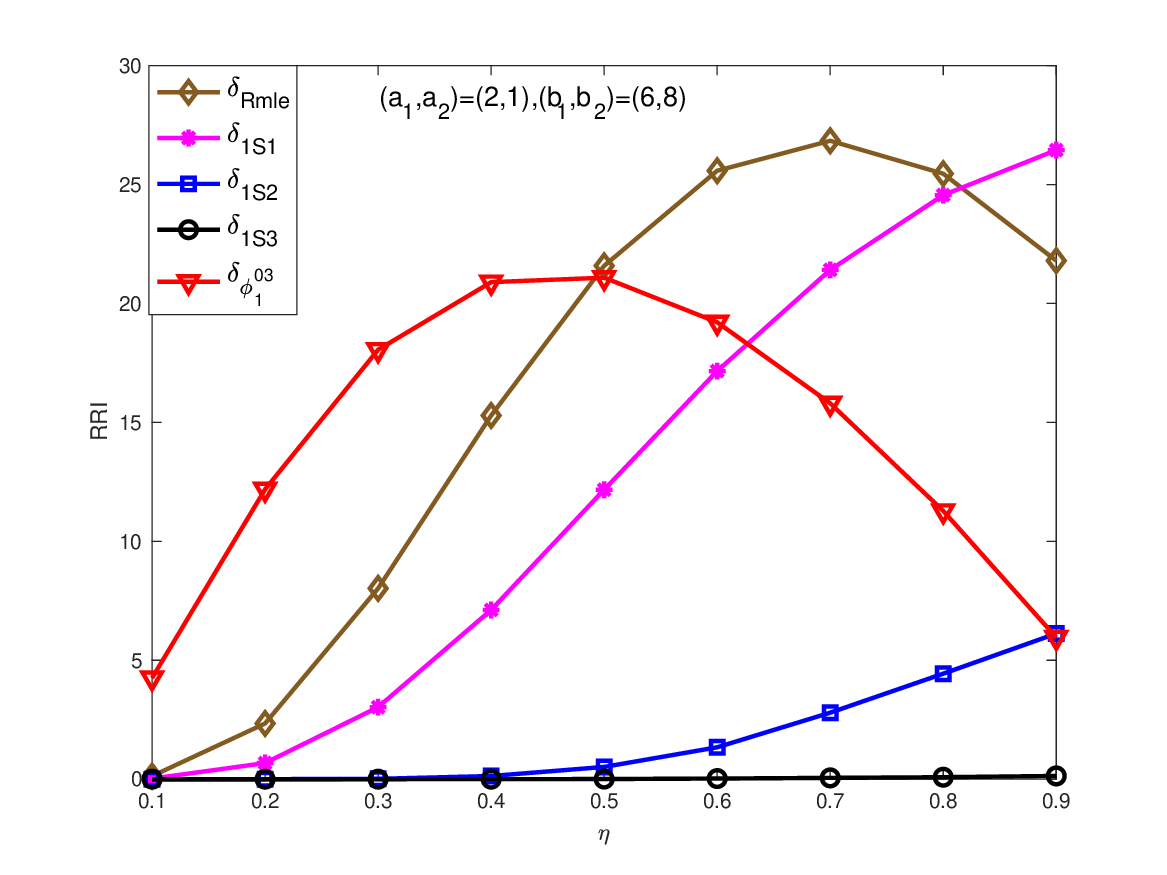}} 
	\end{center}
	\caption{Relative risk improvement under symmetric loss function $L_3(t)$ for $\sigma_{1}$\label{fig2L3}}
\end{figure}

In Figure \ref{fig1L1} and \ref{fig2L1} we have plotted the improved estimators of $\sigma_1$ with respect to quadratic loss function $L_1(t)$. The following observations are made from these figures. 
 \begin{itemize}
 	\item[(i)] RRI of $\delta_{\phi_{1}^{01}}$ and $\delta_{Rmle}$ increases and then decreases as function of $\eta$.  The RRI of $\delta_{1S1},\delta_{1S2},\delta_{1S3}$ is an increasing function of $\eta$. 
 	\item [(ii)] Risk performance of $\delta_{Rmle}$ is better than  $\delta_{1S1},\delta_{1S2},\delta_{1S3}$ for all values of $\eta$. 
 \item [(iii)]  $\delta_{\phi_{1}^{01}}$ shows  significantly better performance over $\delta_{1S1},\delta_{1S2},\delta_{1S3}$  and $\delta_{Rmle}$, for smaller  and moderate values of $\eta$ approximately $\eta \le 0.6$. The region of improvement increases as sample size increases. 
 \item [(iv)] For small sample size and $\mu_1$, $\mu_2$ is close to $0$ all the estimators perform better than $\delta_{\phi_{1}^{01}}$, when $\eta$ close to $1$ that is approximately more than $0.8$.  
 \item[(v)] The region of improvement of $\delta_{\phi_{1}^{01}}$ become larger than the region of improvement of $\delta_{1S1},\delta_{1S2},\delta_{1S3}$,  $\delta_{Rmle}$ when we increase the sample
 \item[(vi)] The estimators $\delta_{1S1},\delta_{1S2}$ has almost equal risk performance when $\mu_1,\mu_2$ are close to zero, however as $\mu_1,\mu_2$ increases  $\delta_{1S1}$ performs better than $\delta_{1S2}$.
 \end{itemize}
 In the Figure \ref{fig1L2}, \ref{fig2L2}, \ref{fig3L2} we have plotted the RRI of the improved estimators with respect ot entropy loss $L_2(t)$ and In Figures \ref{fig1L3} and \ref{fig2L3} we have plotted the RRI of improved estimators under symmetric loss $L_3(t)$. We have similar type of observations for these loss functions also.

\begin{figure}[ht]
								\begin{center}
									\subfigure[\tiny{$(n_1,n_2)=(8,10) ,(\mu_1,\mu_2)=(0,0) $}]{\includegraphics[height=5.5cm,width=8cm]{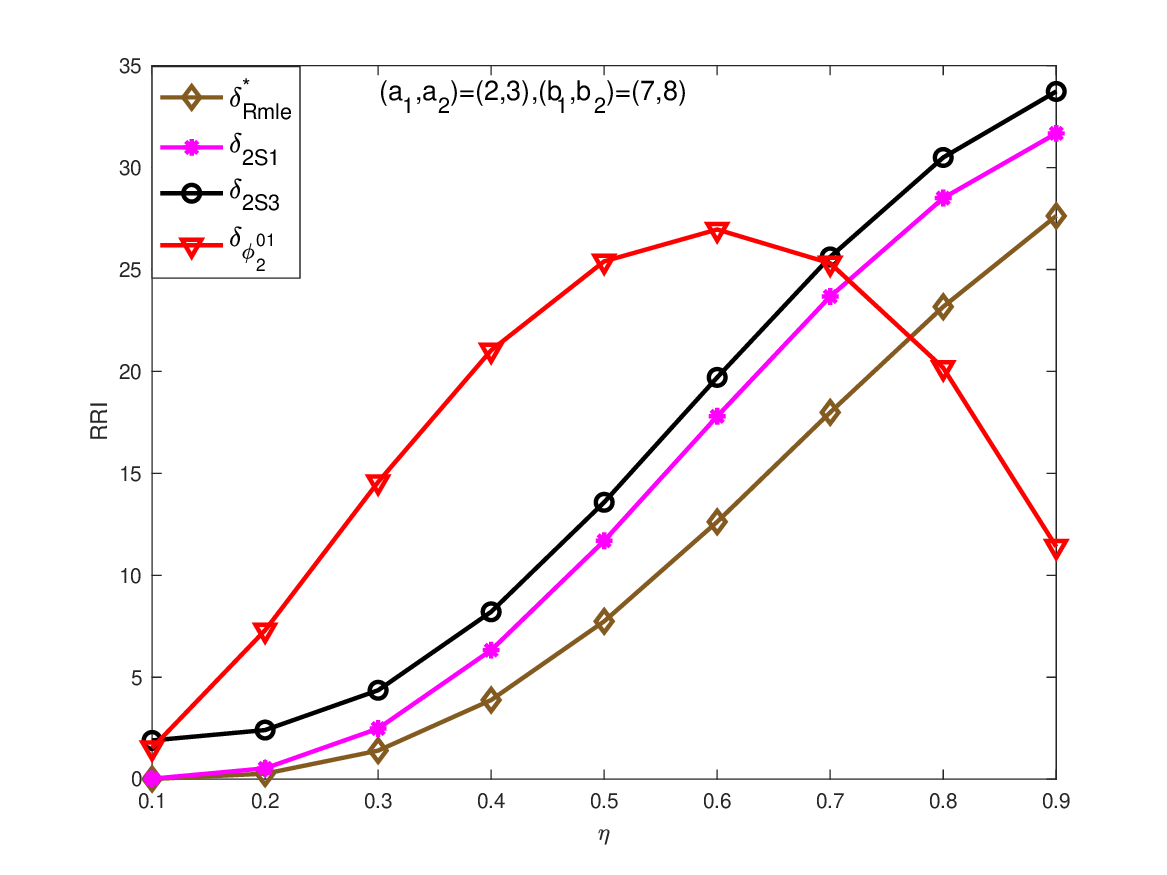}}
									\hspace{1cm} 
									\subfigure[\tiny{$(n_1,n_2)=(8,10) ,(\mu_1,\mu_2)=(0,0) $}]{\includegraphics[height=5.5cm,width=8cm]{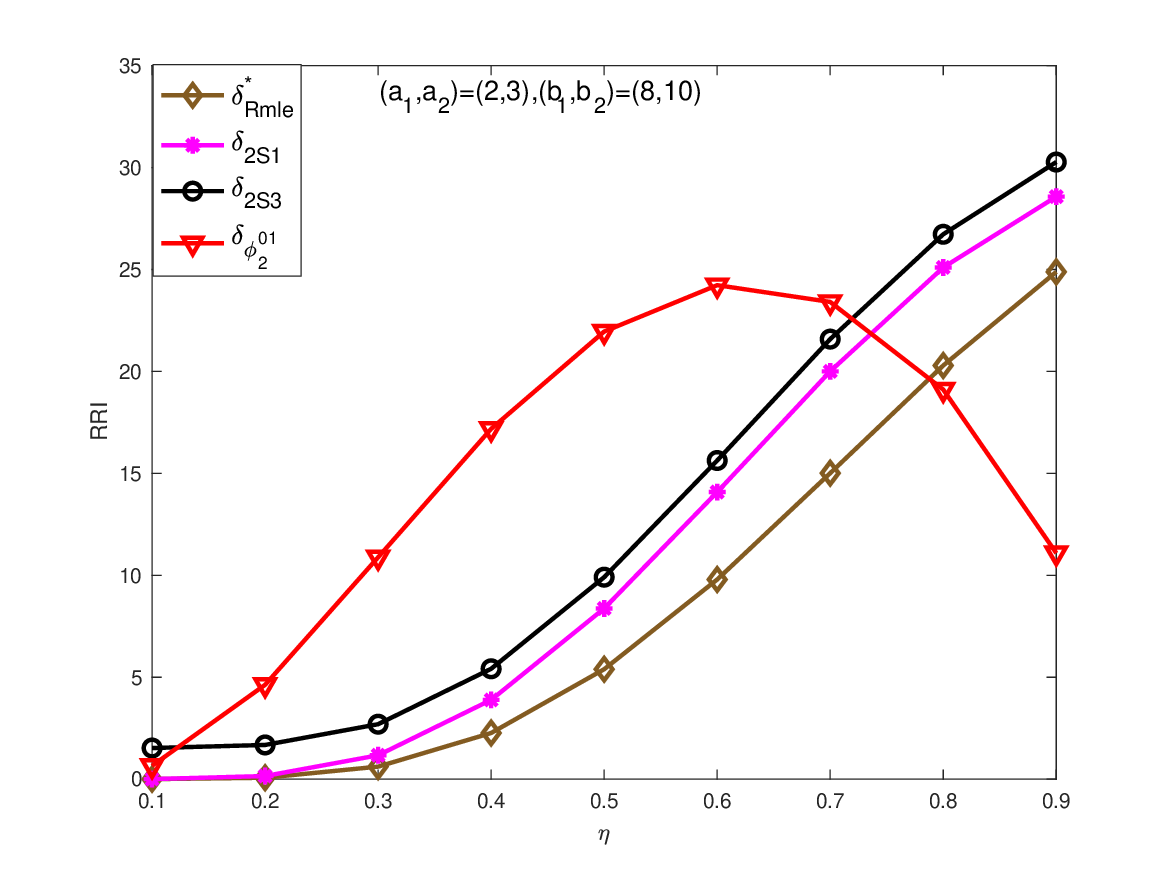}} 
									\hspace{1cm}
									\subfigure[\tiny{$(n_1,n_2)=(12,12) ,(\mu_1,\mu_2)=(0.05,0.03) $}]{\includegraphics[height=5.5cm,width=8cm]{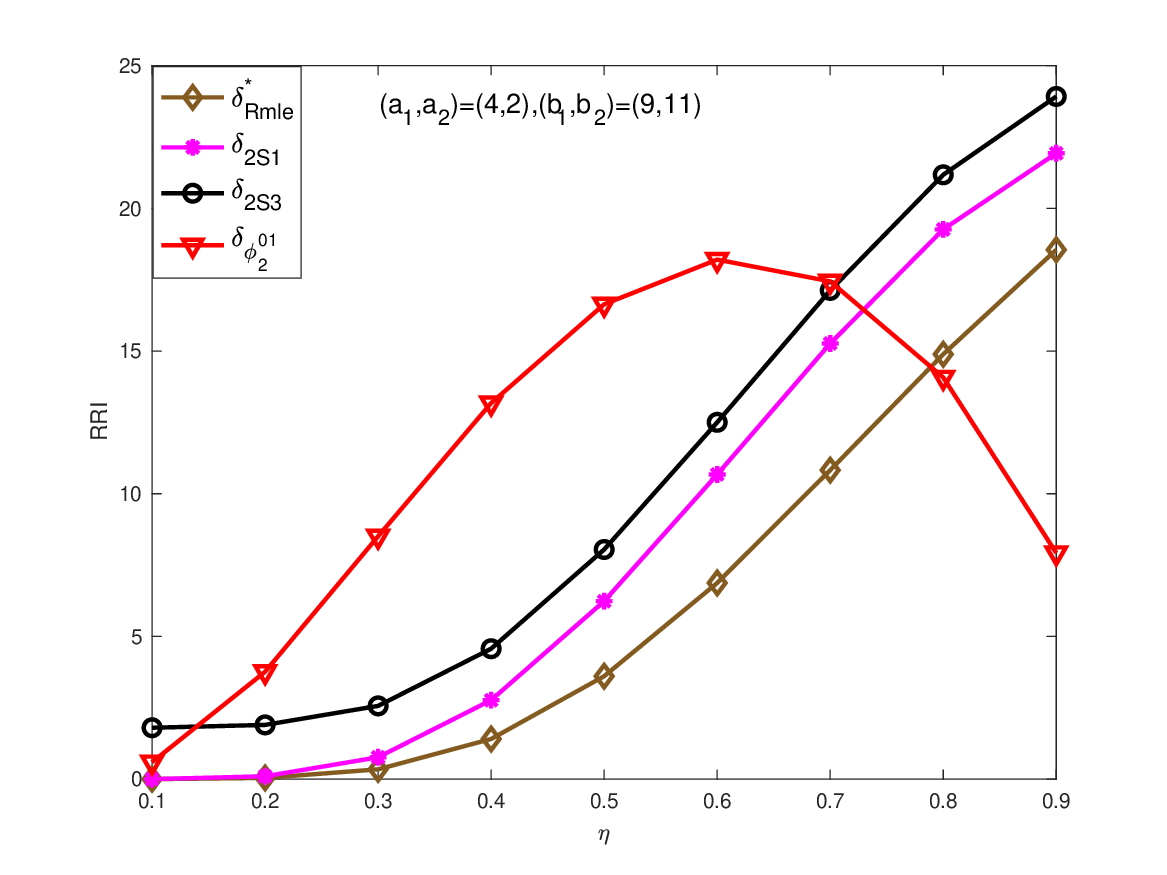}}
									\hspace{1cm}
									\subfigure[\tiny{$(n_1,n_2)=(12,12) ,(\mu_1,\mu_2)=(0.05,0.03) $}]{\includegraphics[height=5.5cm,width=8cm]{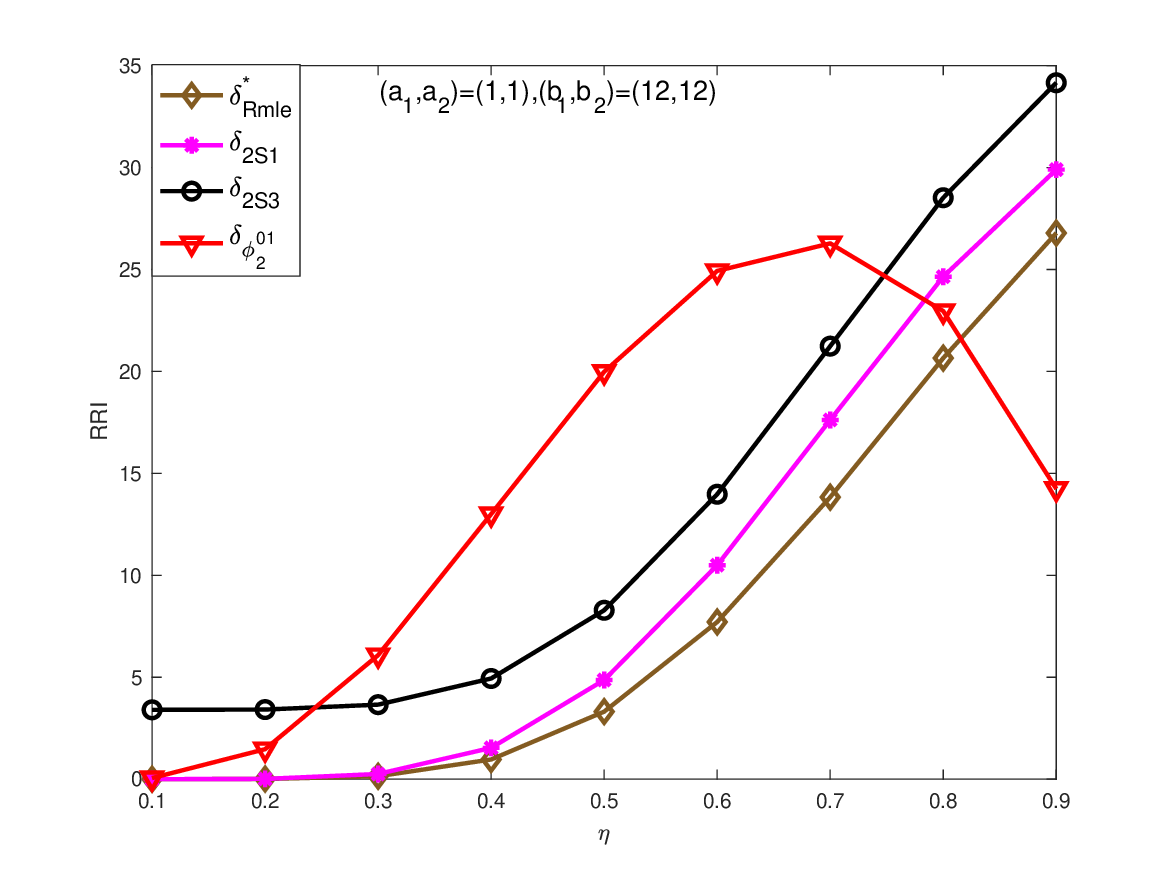}}
									\subfigure[\tiny{$(n_1,n_2)=(10,9) ,(\mu_1,\mu_2)=(0.1,0.1) $}]{\includegraphics[height=5.5cm,width=8cm]{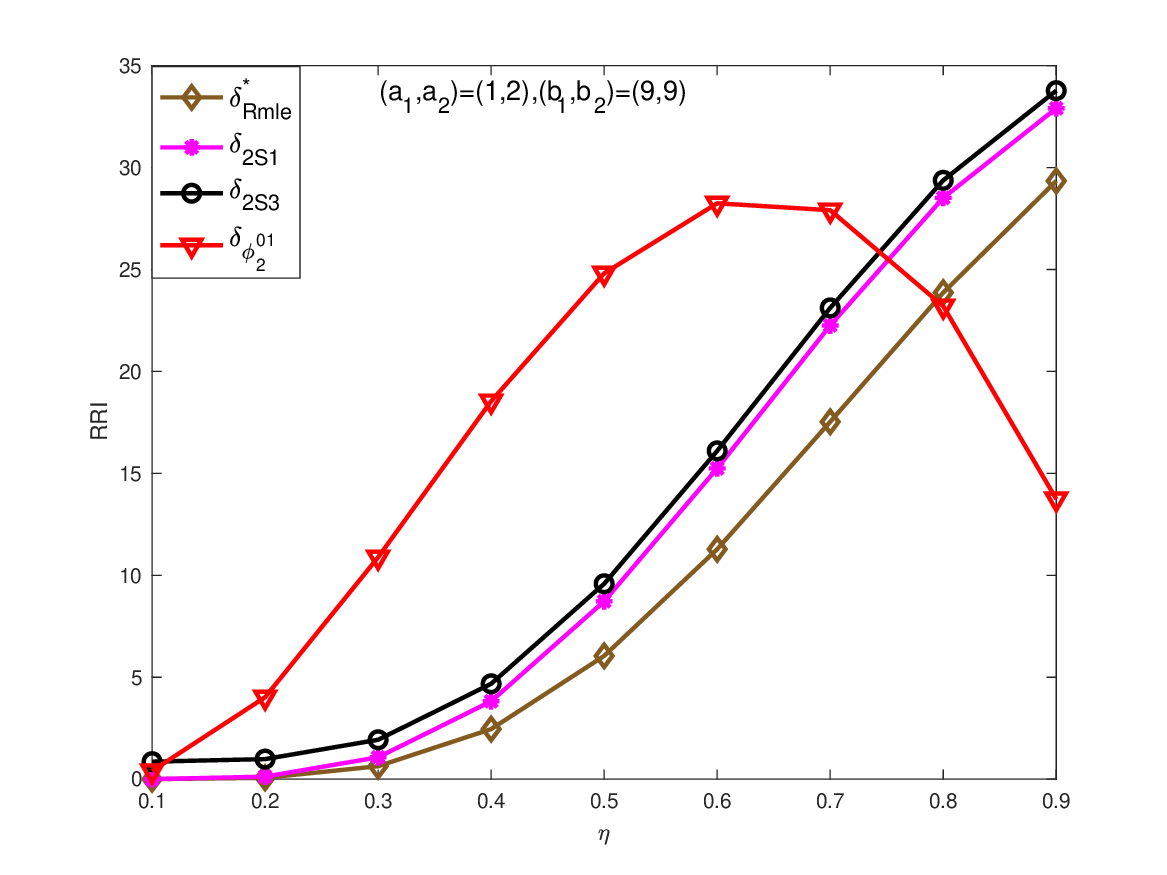}}
									\hspace{1cm}
									\subfigure[\tiny{$(n_1,n_2)=(10,9) ,(\mu_1,\mu_2)=(0.1,0.15) $}]{\includegraphics[height=5.5cm,width=8cm]{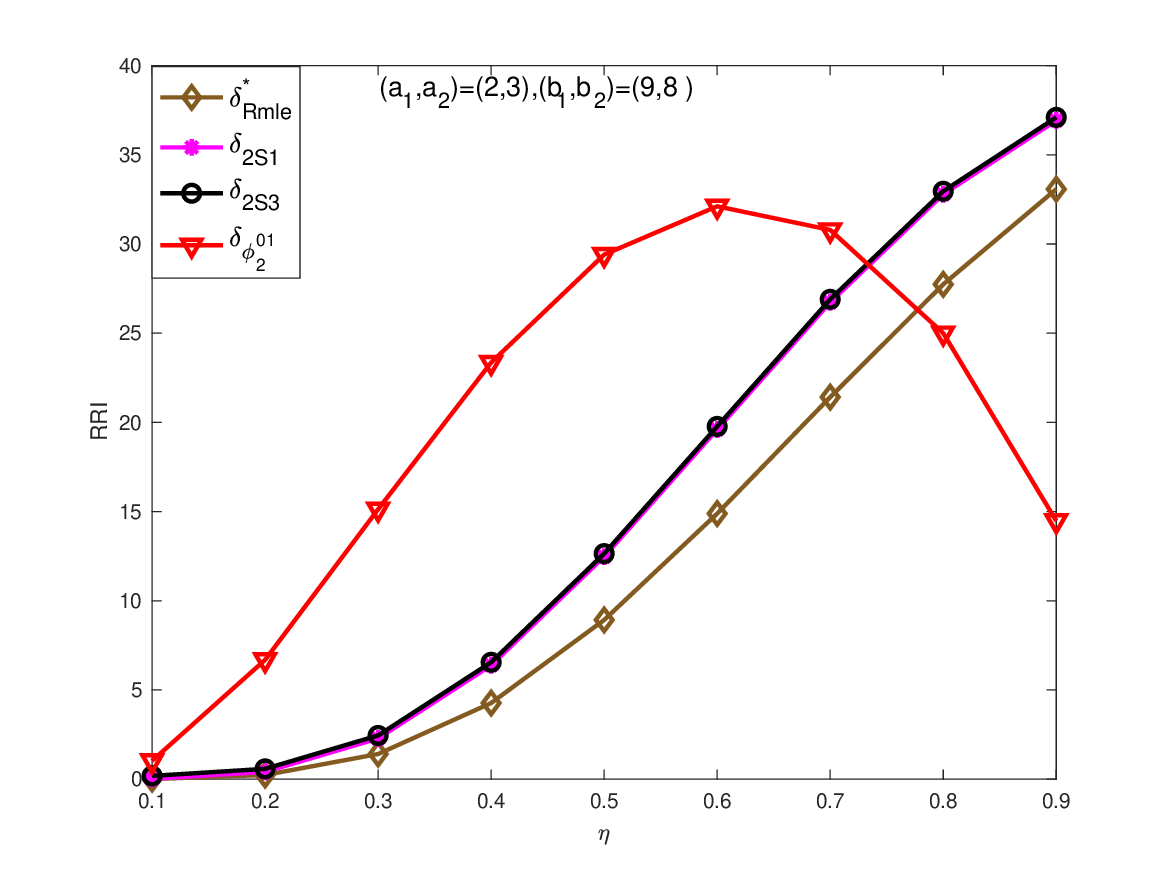}}
										\end{center}
									\caption{Relative risk improvement under quadratic loss function $L_1(t)$ for $\sigma_{2}$}\label{figs21L1}
								\end{figure}

\begin{figure}[ht]
									\begin{center}
									\subfigure[\tiny{$(n_1,n_2)=(14,15) ,(\mu_1,\mu_2)=(0.4,0.7) $}]{\includegraphics[height=5.5cm,width=8cm]{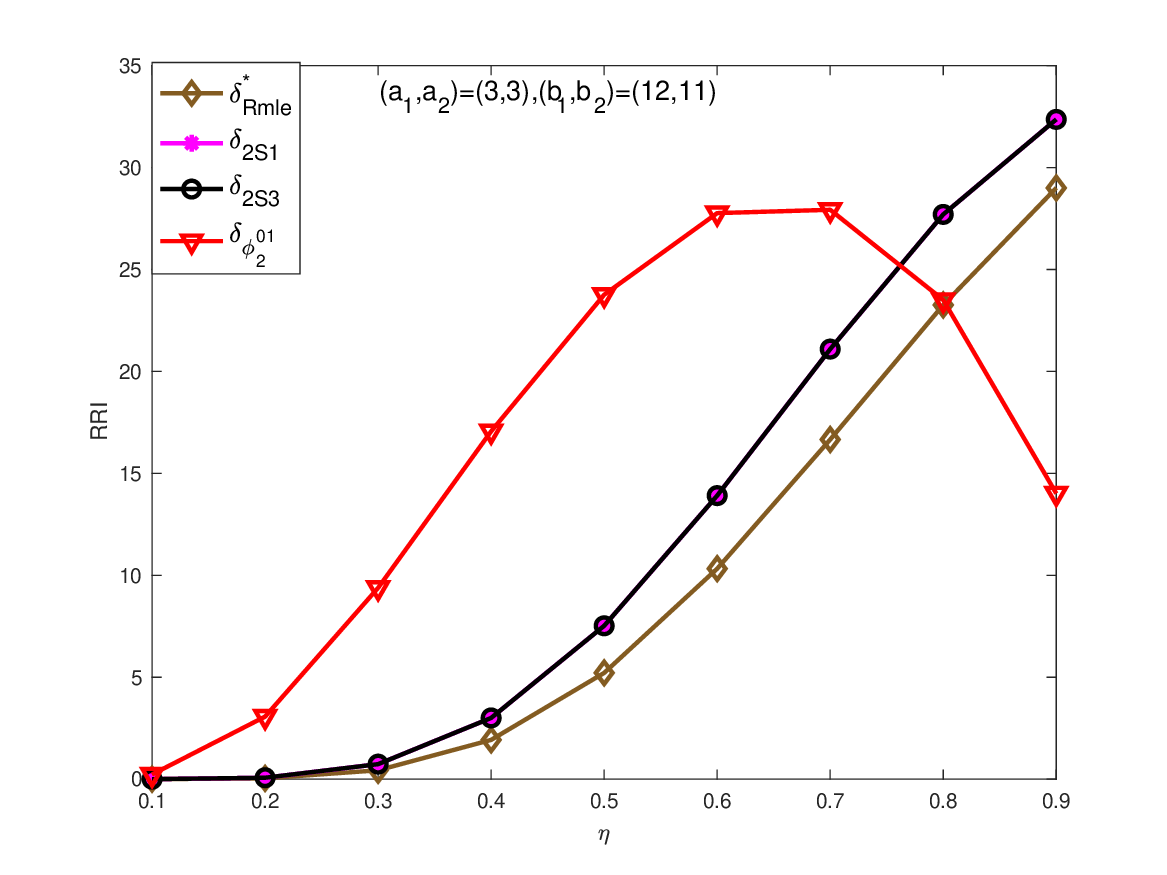}}
									\hspace{1cm}
									\subfigure[\tiny{$(n_1,n_2)=(14,15) ,(\mu_1,\mu_2)=(0.4,0.7)$}]{\includegraphics[height=5.5cm,width=8cm]{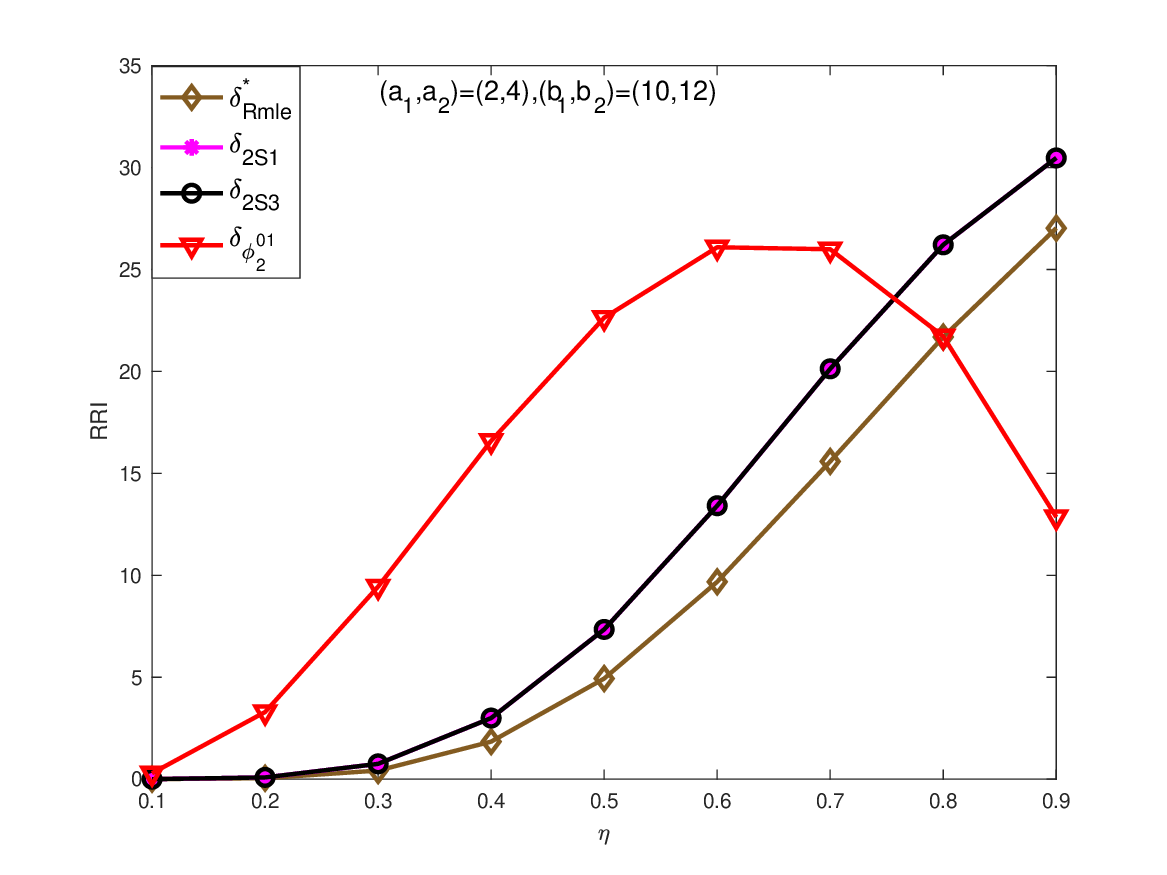}}
									
									\subfigure[\tiny{$(n_1,n_2)=(7,8) ,(\mu_1,\mu_2)=(-0.1,-0.2) $}]{\includegraphics[height=5.5cm,width=8cm]{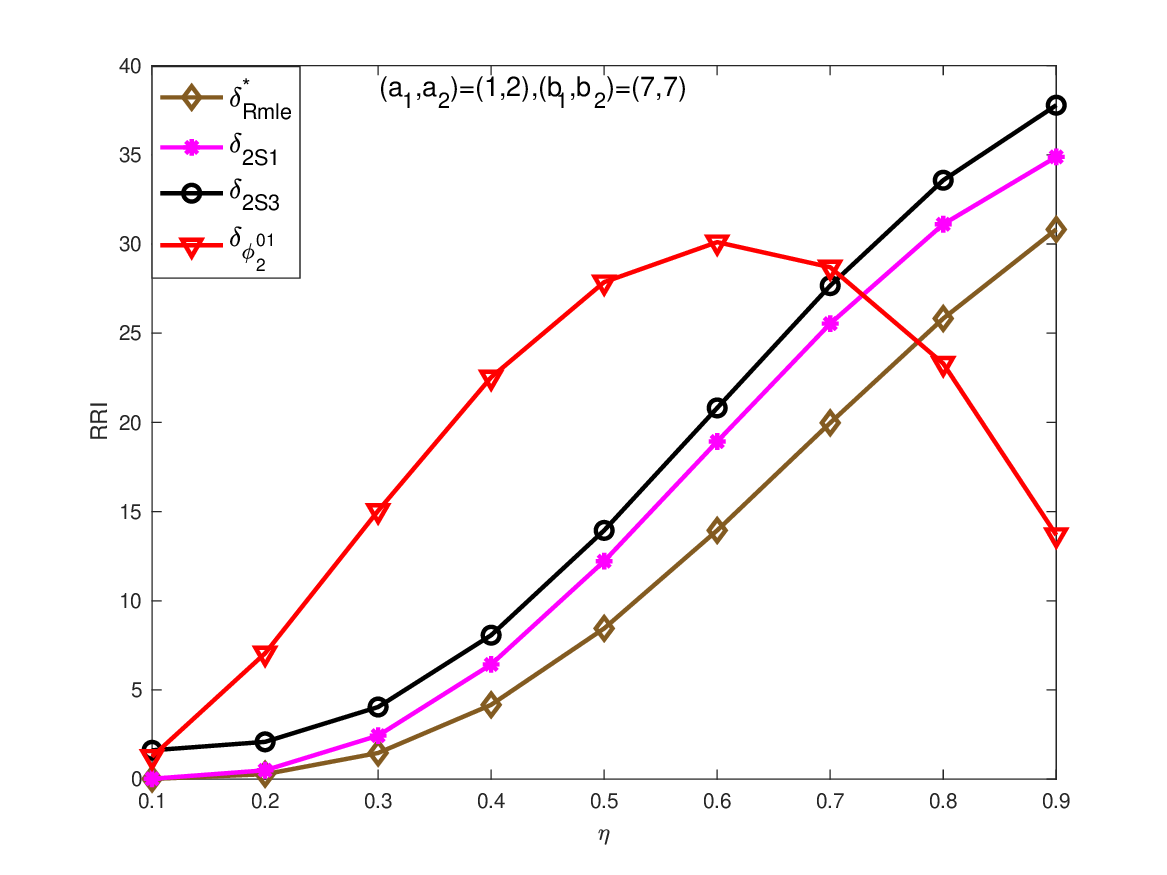}}
									\hspace{1cm}
									\subfigure[\tiny{$(n_1,n_2)=(7,8) ,(\mu_1,\mu_2)=(-0.1,-0.2) $}]{\label{fig:edge-b3}\includegraphics[height=5.5cm,width=8cm]{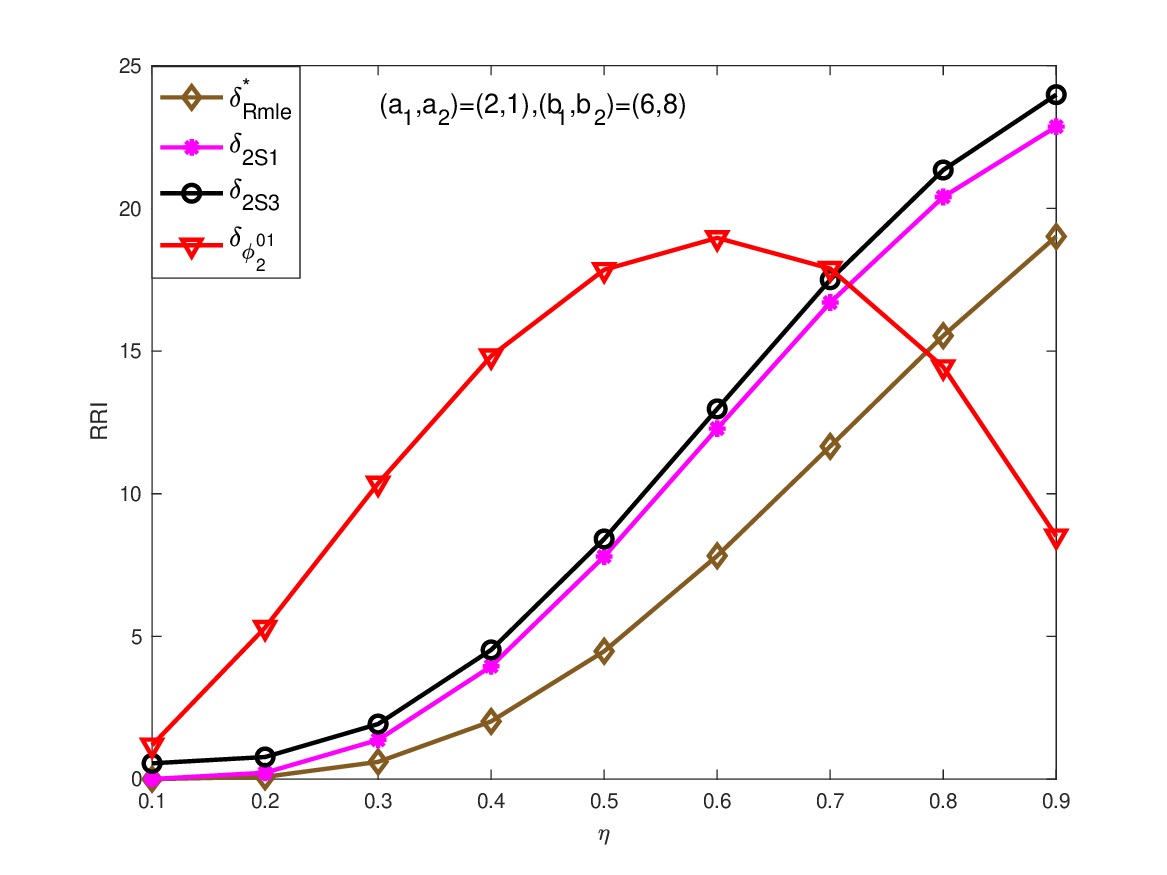}} 
								\end{center}
								\caption{Relative risk improvement under quadratic loss function $L_1(t)$ for $\sigma_{2}$}\label{figs22L1}
								\label{fig:edge}
							\end{figure}

\begin{figure}[ht]
							\begin{center}
								\subfigure[\tiny{$(n_1,n_2)=(8,10) ,(\mu_1,\mu_2)=(0,0) $}]{\includegraphics[height=5.5cm,width=8cm]{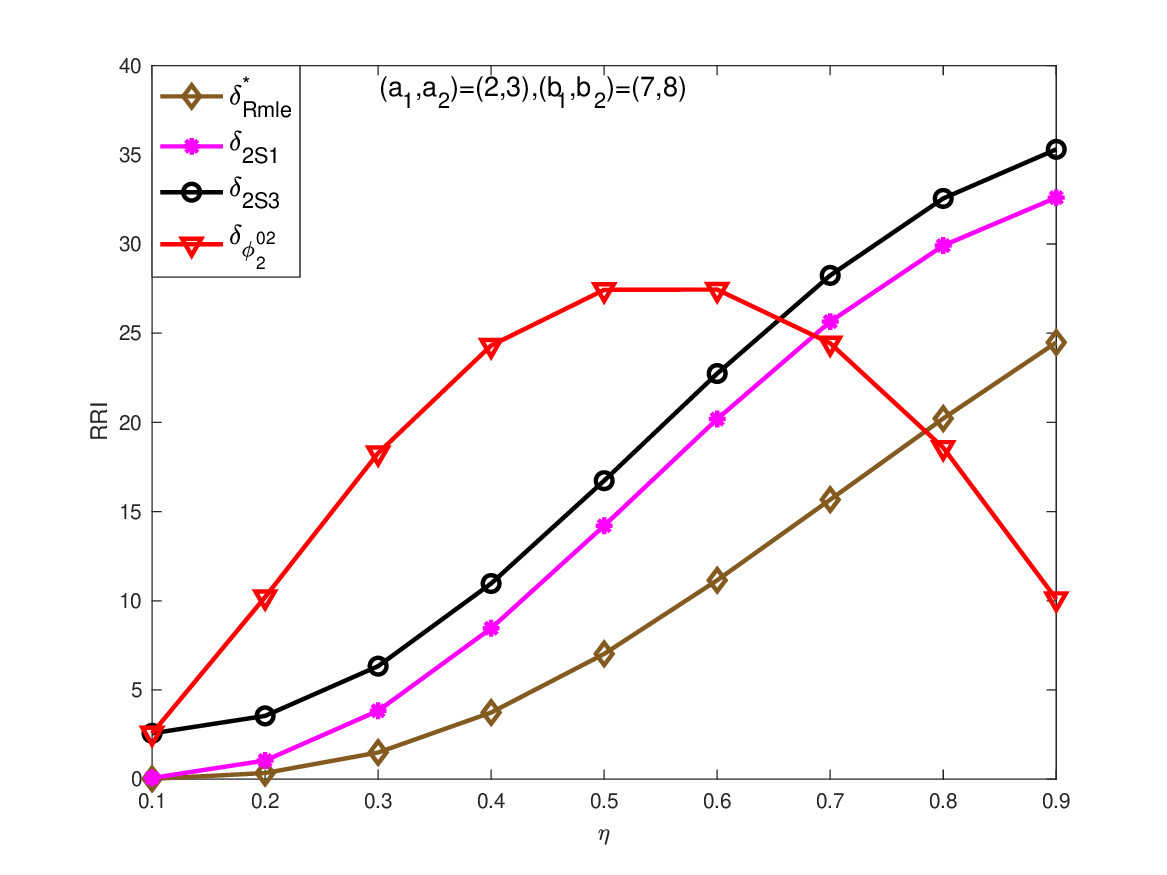}}
								\hspace{1cm} 
								\subfigure[\tiny{$(n_1,n_2)=(8,10) ,(\mu_1,\mu_2)=(0,0) $}]{\includegraphics[height=5.5cm,width=8cm]{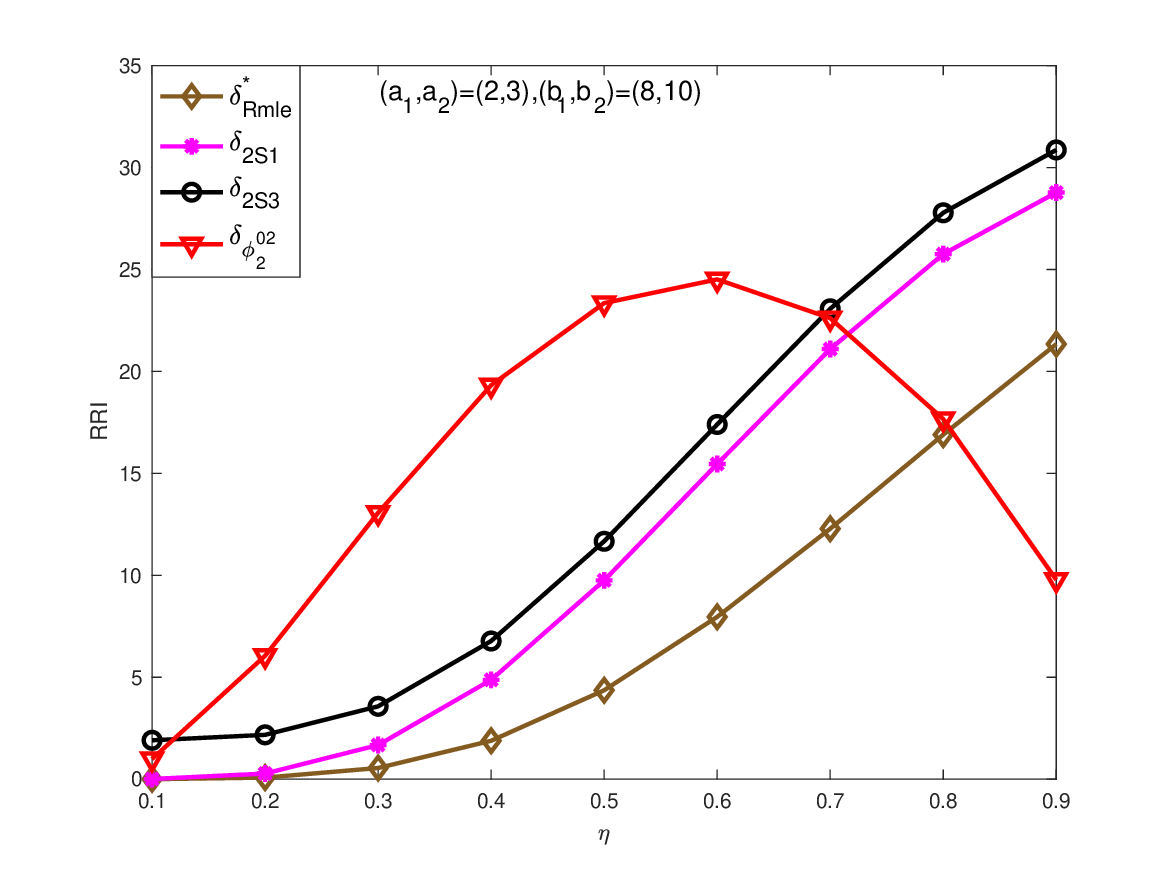}} 
								\hspace{1cm}
							\end{center}
							\caption{Relative risk improvement under entropy loss function $L_2(t)$ for $\sigma_{2}$}\label{figs21L2}
						\end{figure}

\begin{figure}[ht]
							\begin{center}
								\subfigure[\tiny{$(n_1,n_2)=(12,12) ,(\mu_1,\mu_2)=(0.05,0.03) $}]{\includegraphics[height=5.5cm,width=8cm]{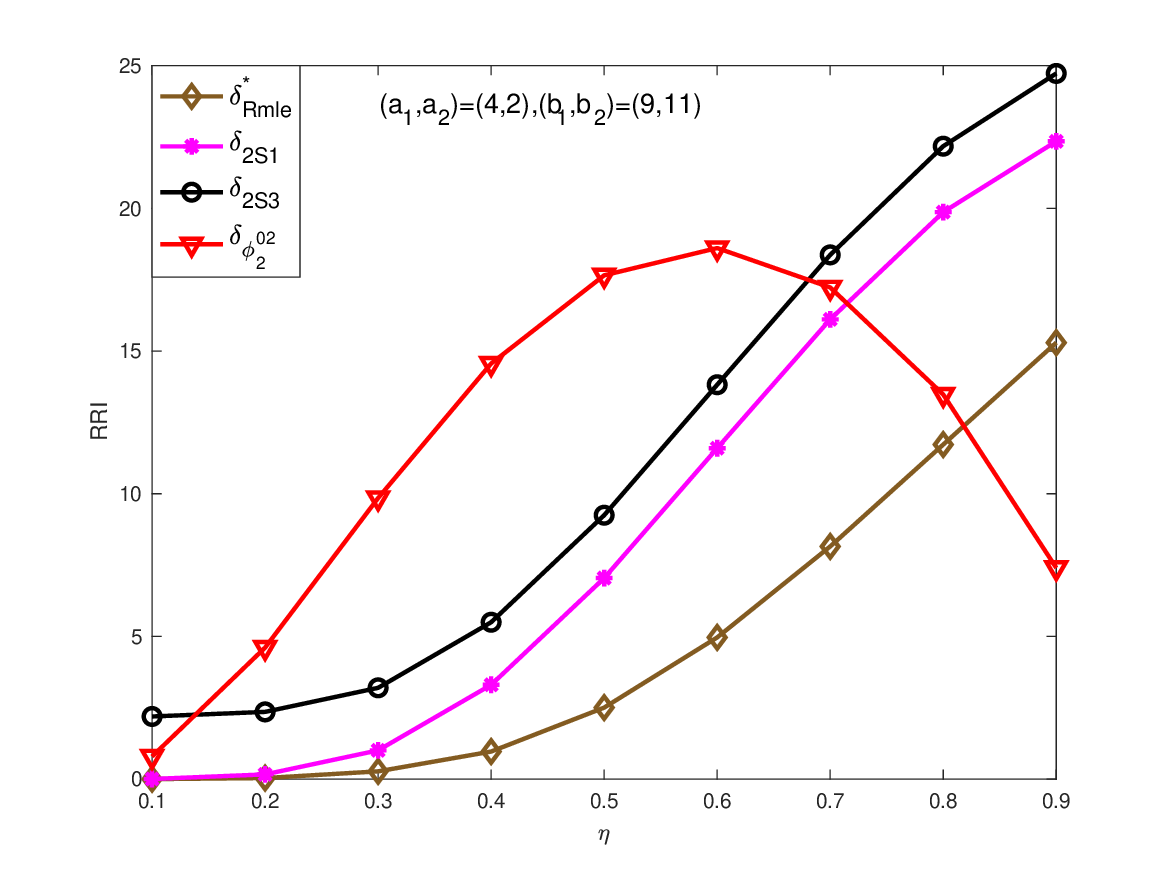}}
								\hspace{1cm}
								\subfigure[\tiny{$(n_1,n_2)=(12,12) ,(\mu_1,\mu_2)=(0.05,0.03) $}]{\includegraphics[height=5.5cm,width=8cm]{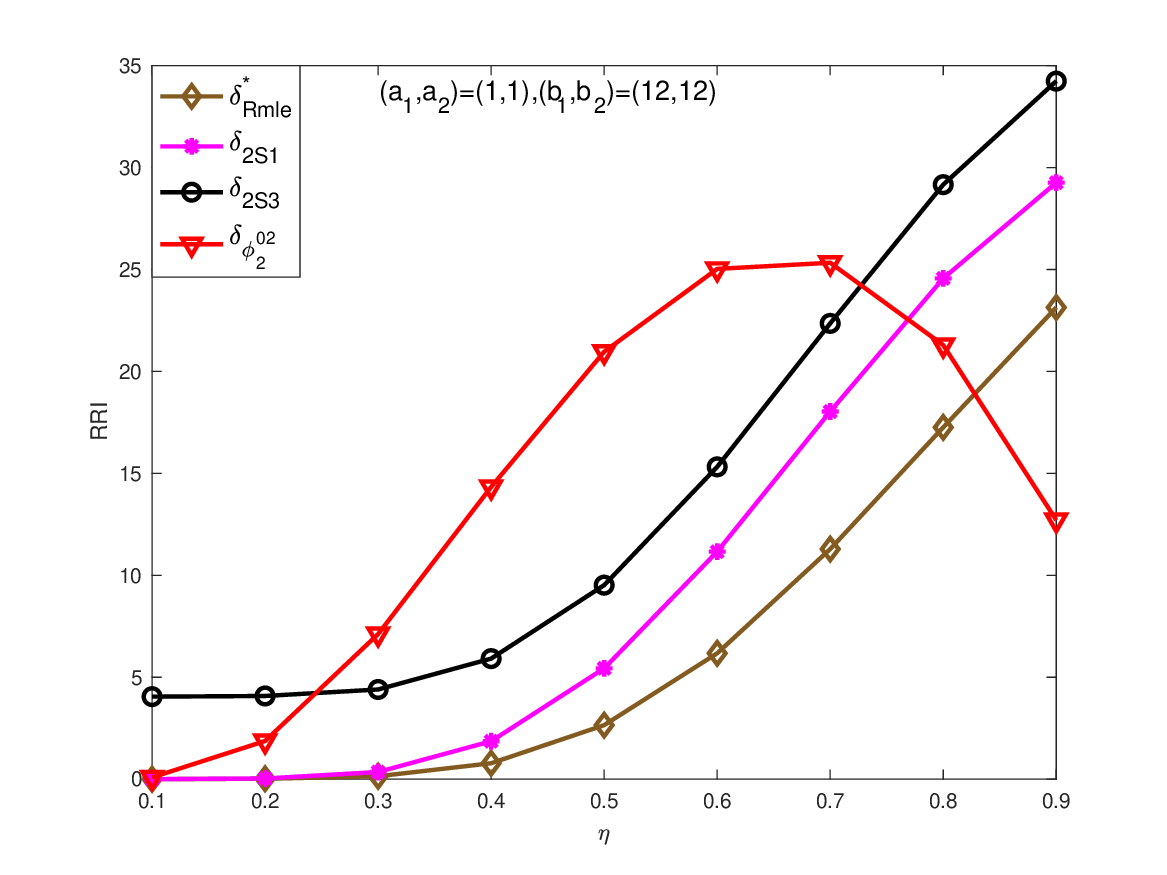}}
								\subfigure[\tiny{$(n_1,n_2)=(10,9) ,(\mu_1,\mu_2)=(0.1,0.1) $}]{\includegraphics[height=5.5cm,width=8cm]{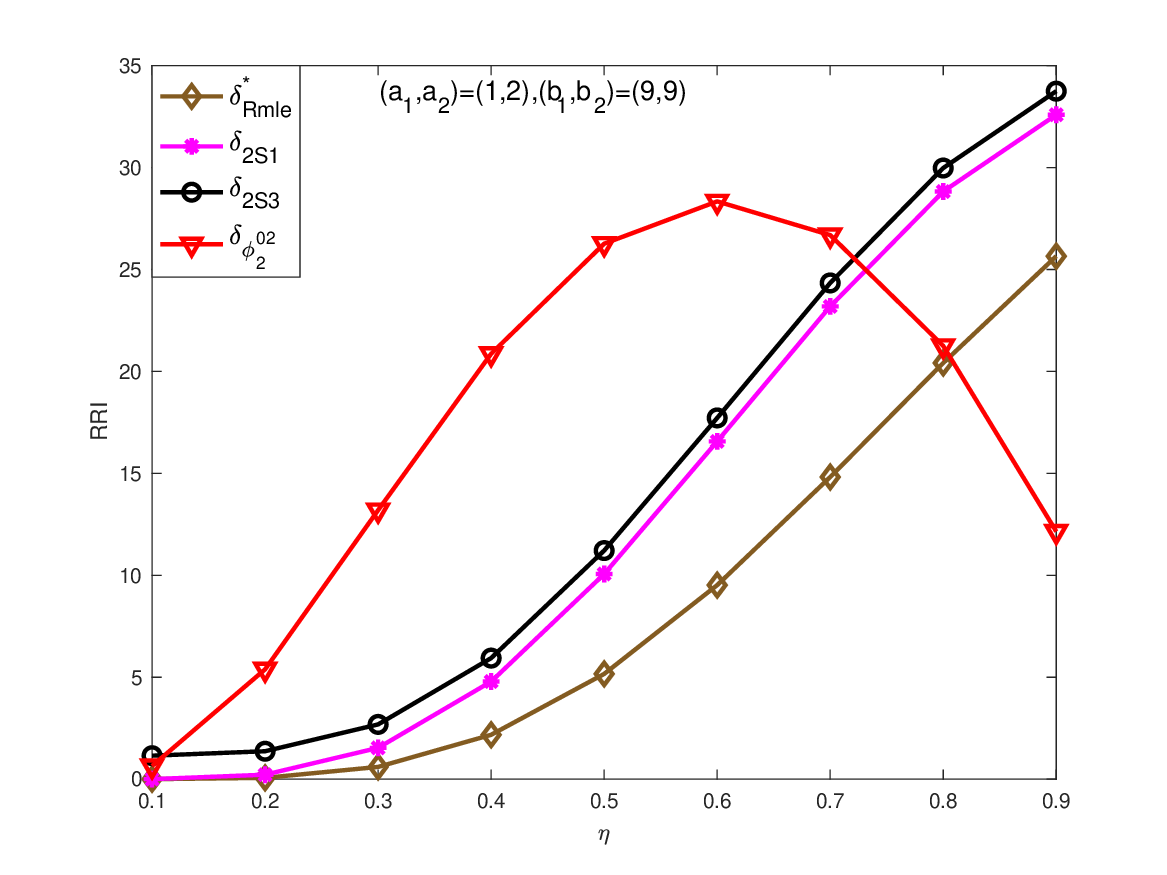}}
								\hspace{1cm}
								\subfigure[\tiny{$(n_1,n_2)=(10,9) ,(\mu_1,\mu_2)=(0.1,0.15) $}]{\includegraphics[height=5.5cm,width=8cm]{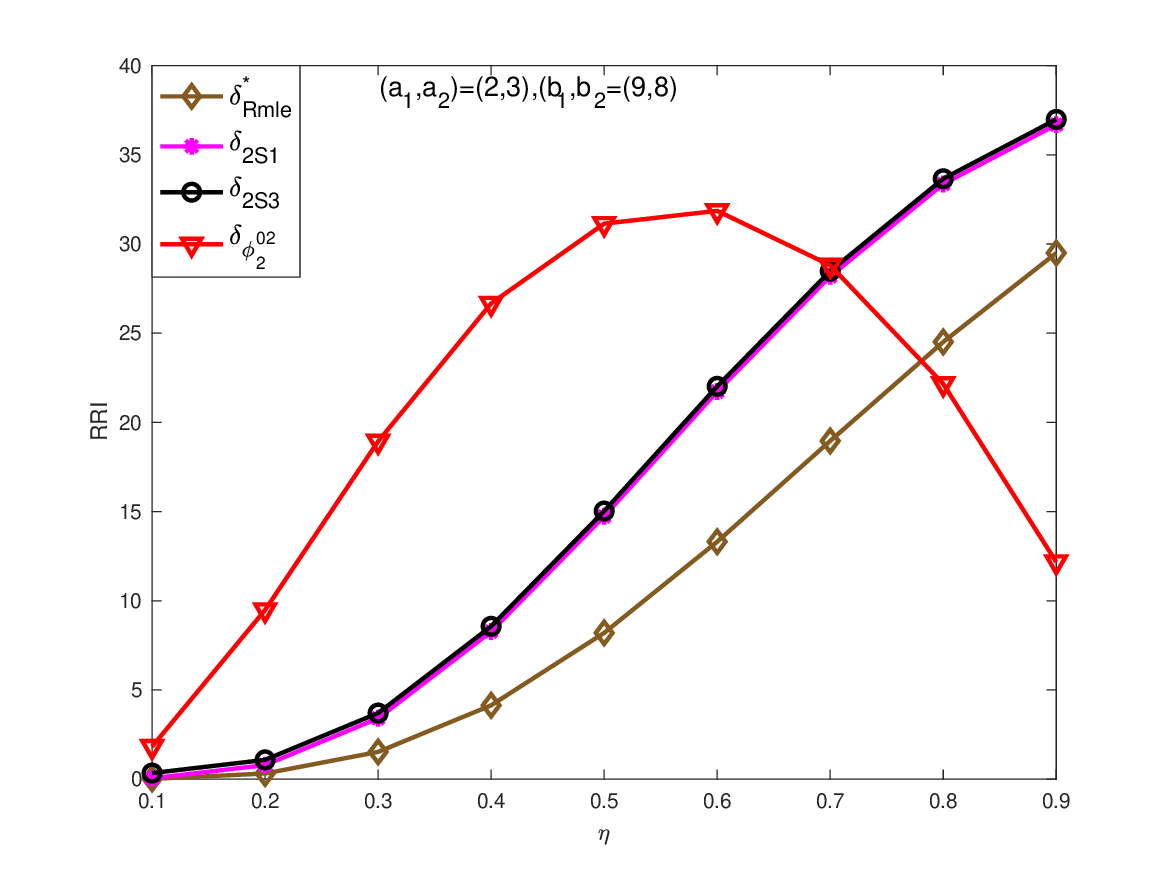}} 
								\hspace{1cm}
								\subfigure[\tiny{$(n_1,n_2)=(14,15) ,(\mu_1,\mu_2)=(0.4,0.7) $}]{\includegraphics[height=5.5cm,width=8cm]{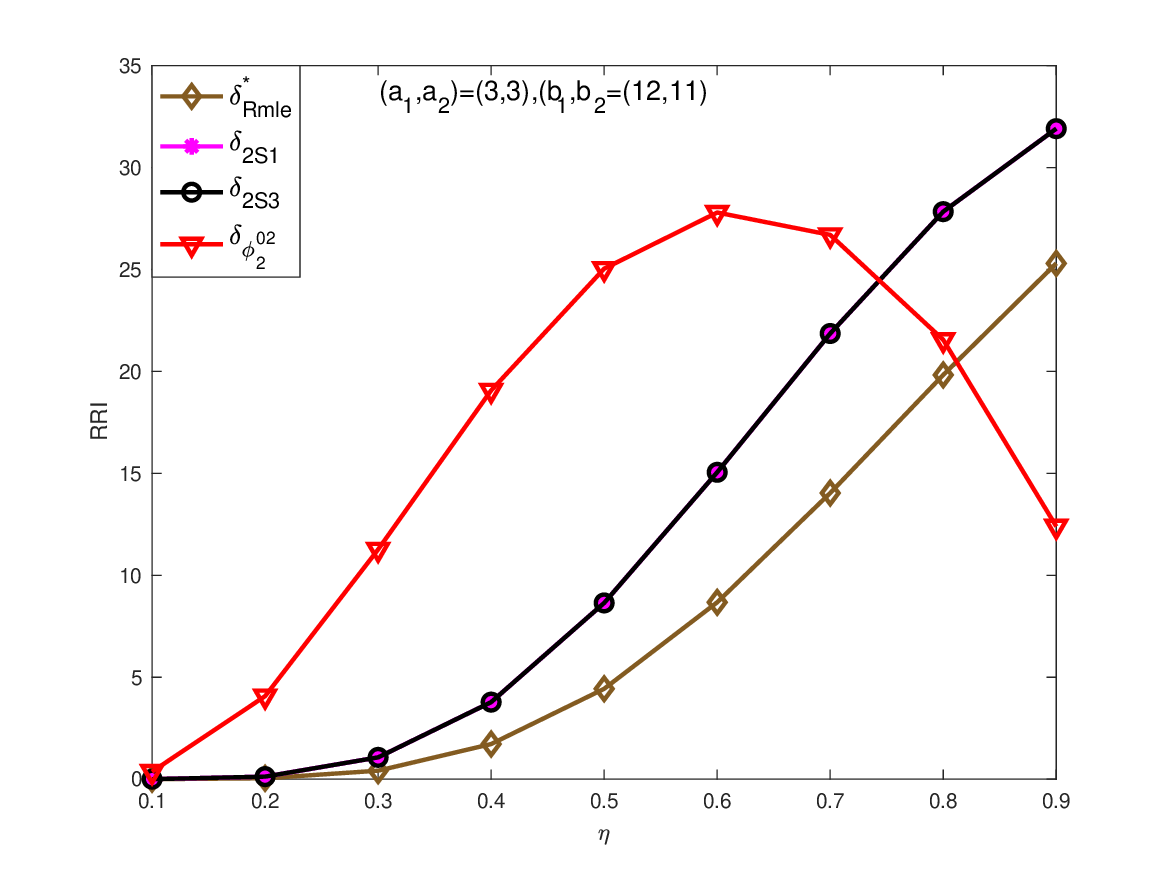}}
								\hspace{1cm}
								\subfigure[\tiny{$(n_1,n_2)=(14,15) ,(\mu_1,\mu_2)=(0.4,0.7)$}]{\includegraphics[height=5.5cm,width=8cm]{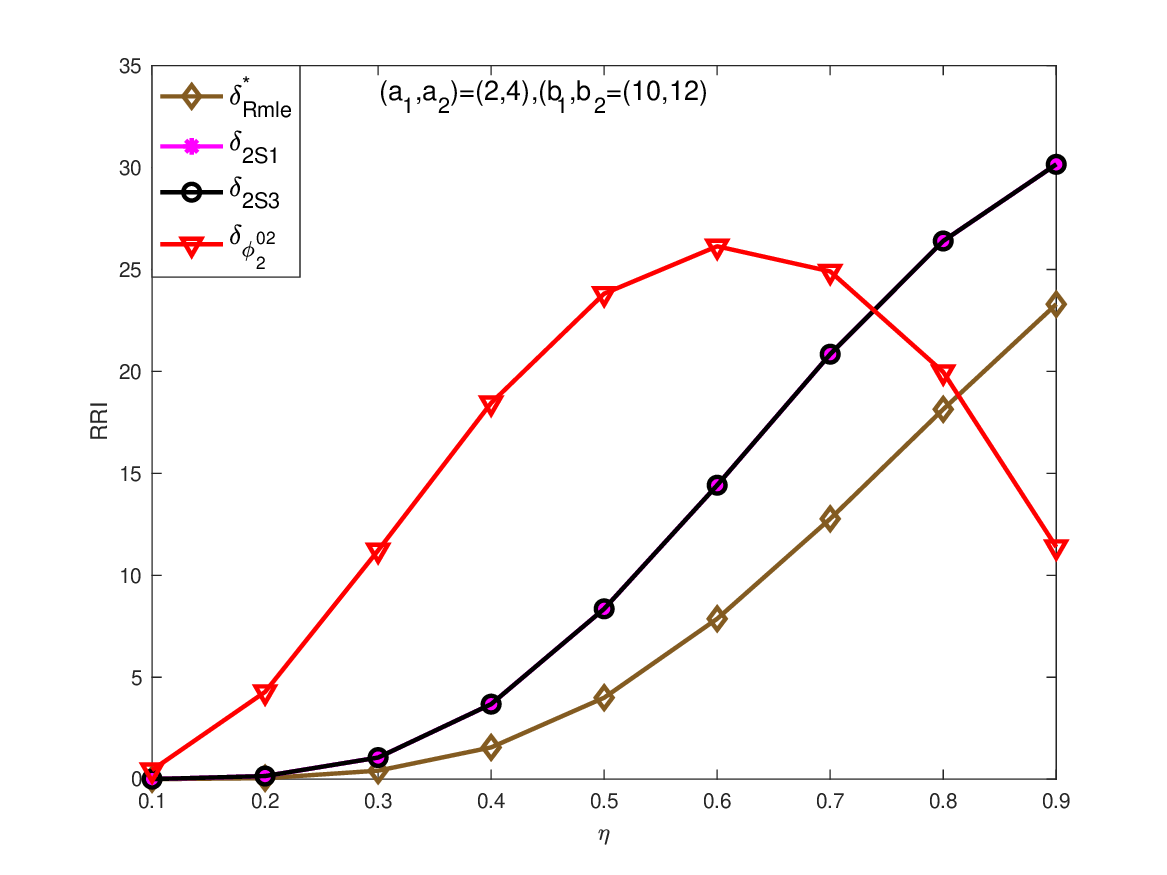}}
						\end{center}
						\caption{Relative risk improvement under entropy loss function $L_2(t)$ for $\sigma_{2}$}\label{figs22l2}
					\end{figure}

\begin{figure}[ht]
						\begin{center}

								\subfigure[\tiny{$(n_1,n_2)=(7,8) ,(\mu_1,\mu_2)=(-0.1,-0.2) $}]{\includegraphics[height=5.5cm,width=8cm]{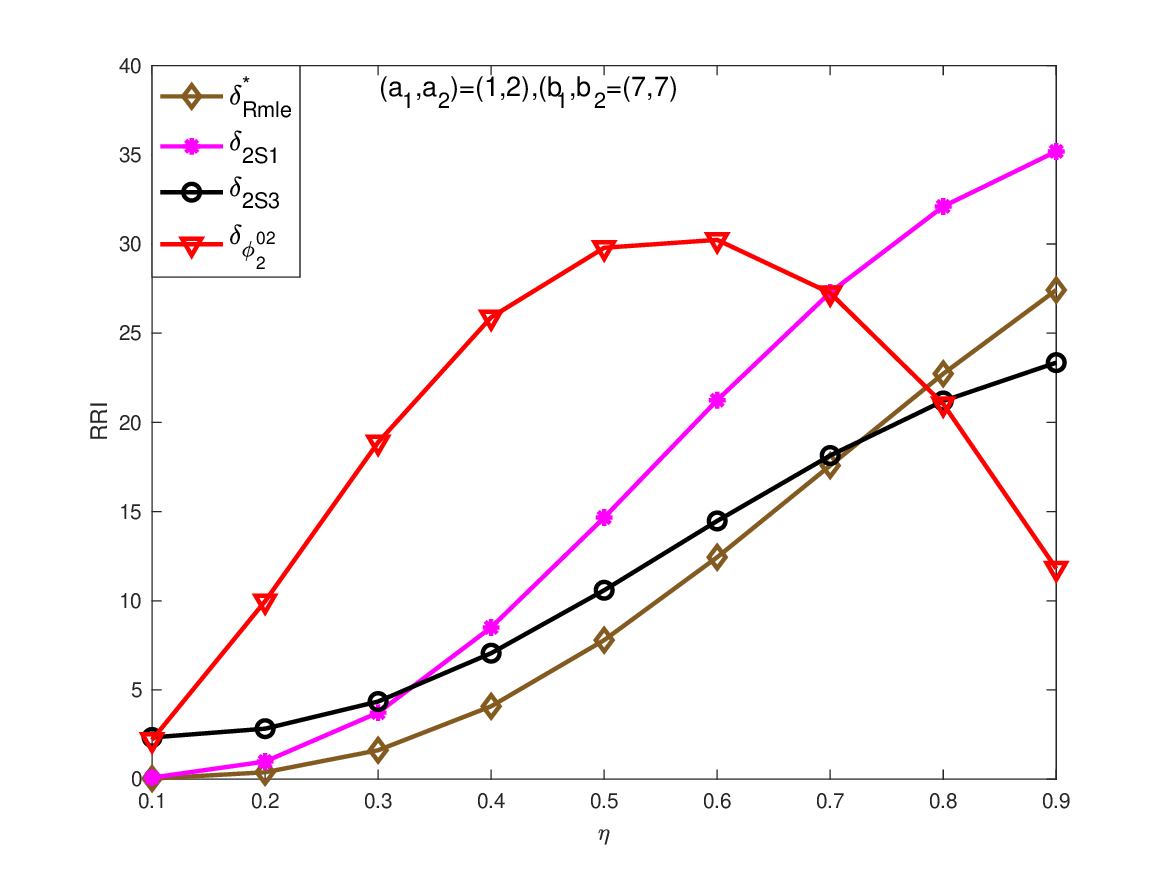}}
								\hspace{1cm}
								\subfigure[\tiny{$(n_1,n_2)=(7,8) ,(\mu_1,\mu_2)=(-0.1,-0.2) $}]{\includegraphics[height=5.5cm,width=8cm]{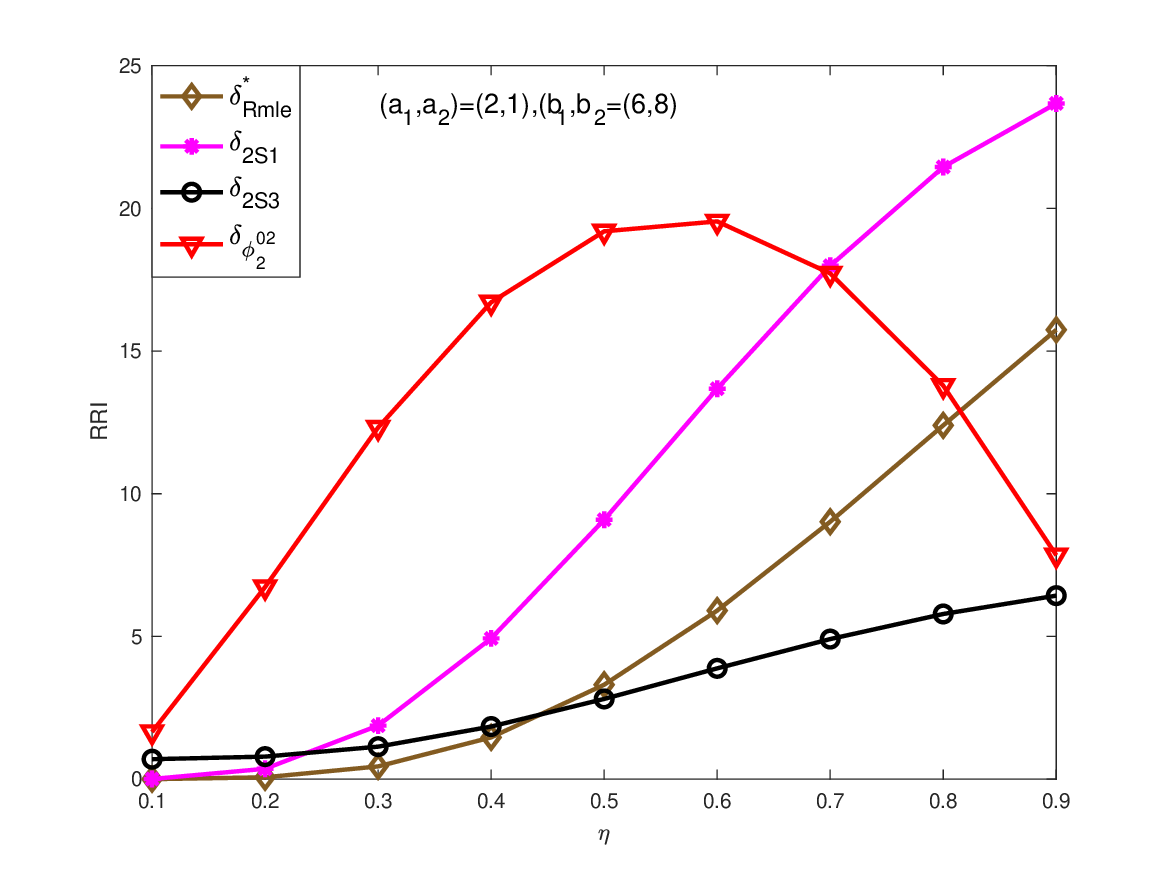}} 
							\end{center}
							\caption{Relative risk improvement under entropy loss function $L_2(t)$ for $\sigma_{2}$}\label{figs23L2}
						\end{figure}

\begin{figure}[ht]
							\begin{center}
								\subfigure[\tiny{$(n_1,n_2)=(8,10) ,(\mu_1,\mu_2)=(0,0) $}]{\includegraphics[height=5.5cm,width=8cm]{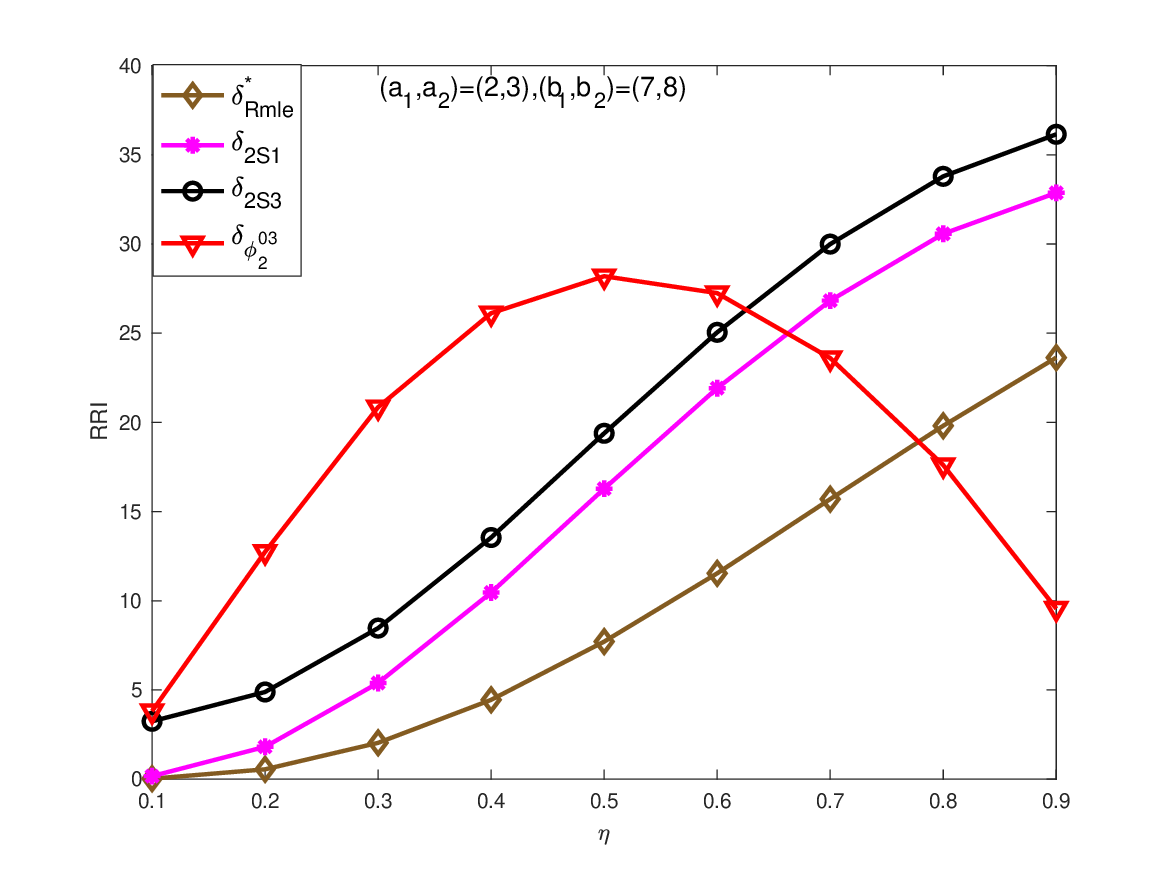}}
								\hspace{1cm} 
								\subfigure[\tiny{$(n_1,n_2)=(8,10) ,(\mu_1,\mu_2)=(0,0) $}]{\includegraphics[height=5.5cm,width=8cm]{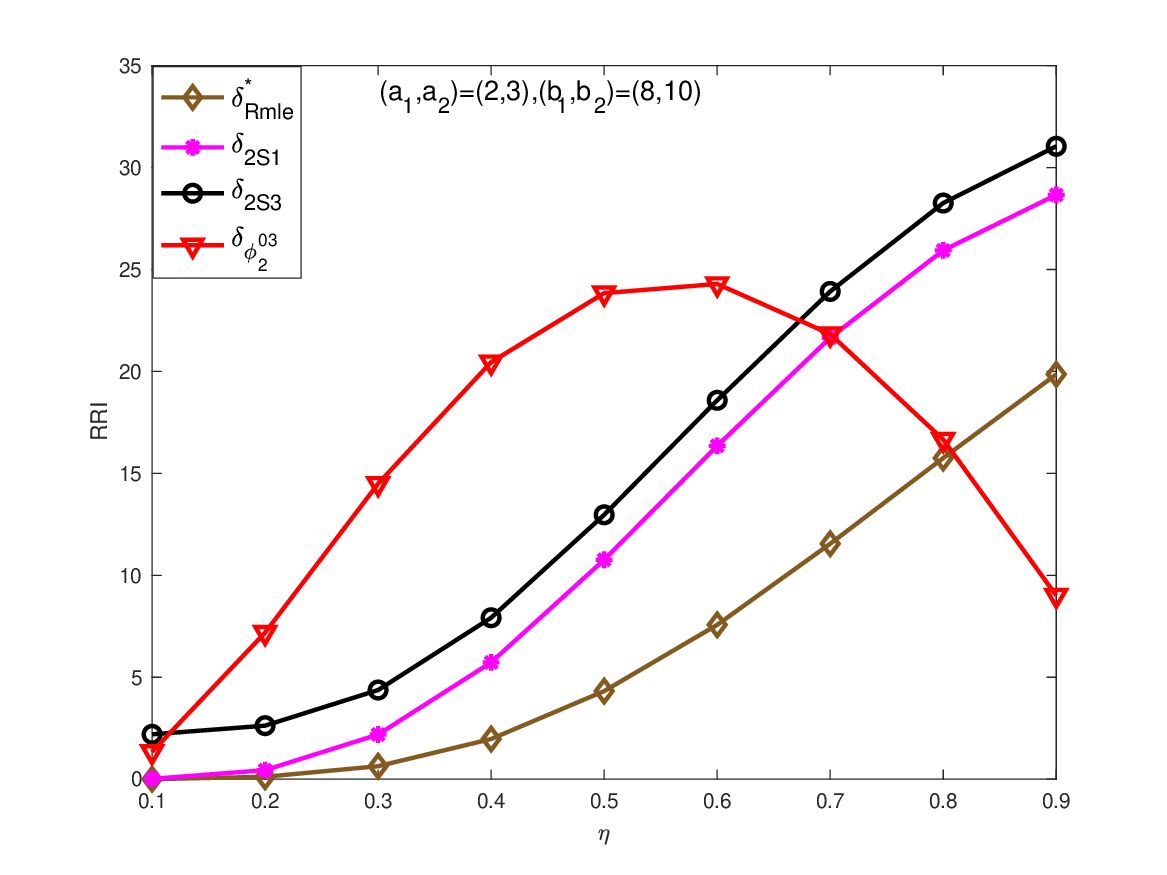}} 
								\hspace{1cm}
								\subfigure[\tiny{$(n_1,n_2)=(12,12) ,(\mu_1,\mu_2)=(0.05,0.03) $}]{\includegraphics[height=5.5cm,width=8cm]{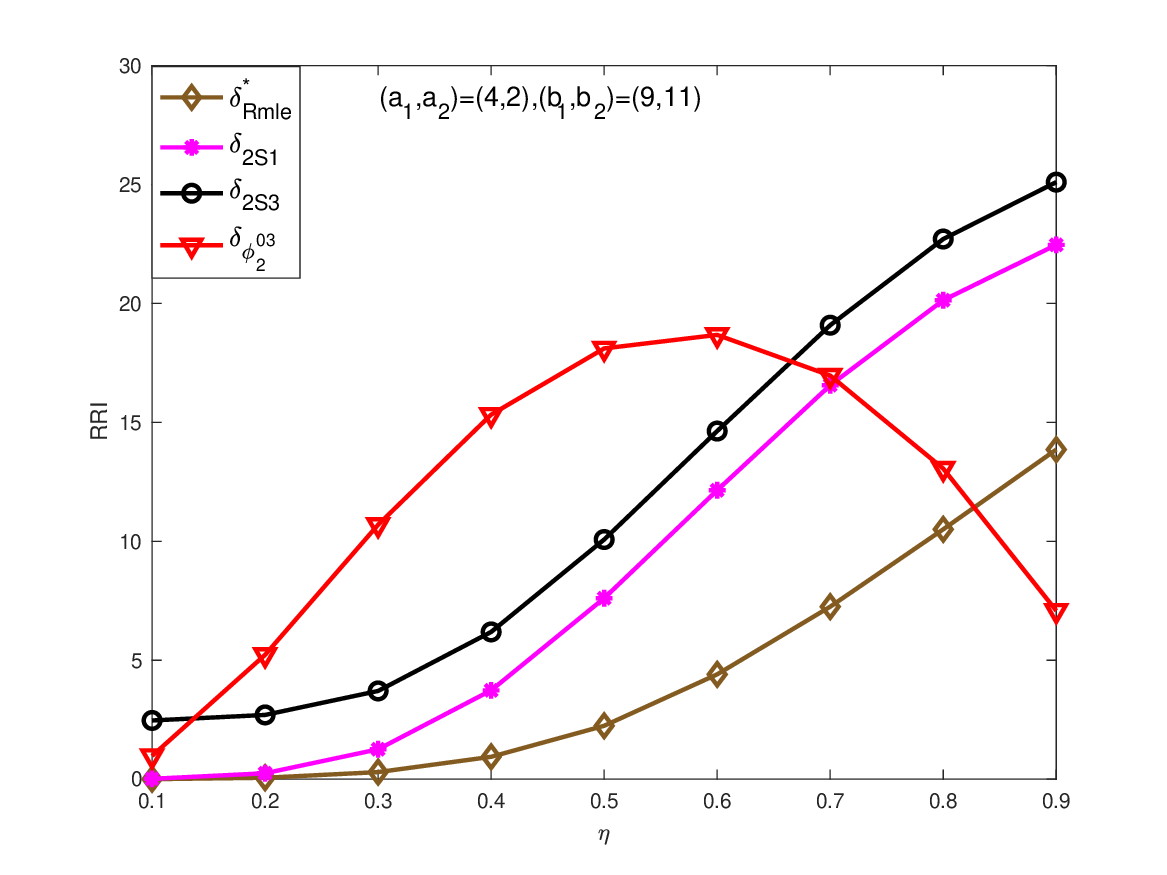}}
								\hspace{1cm}
								\subfigure[\tiny{$(n_1,n_2)=(12,12),(\mu_1,\mu_2)=(0.05,0.03)
									$}]{\includegraphics[height=5.5cm,width=8cm]{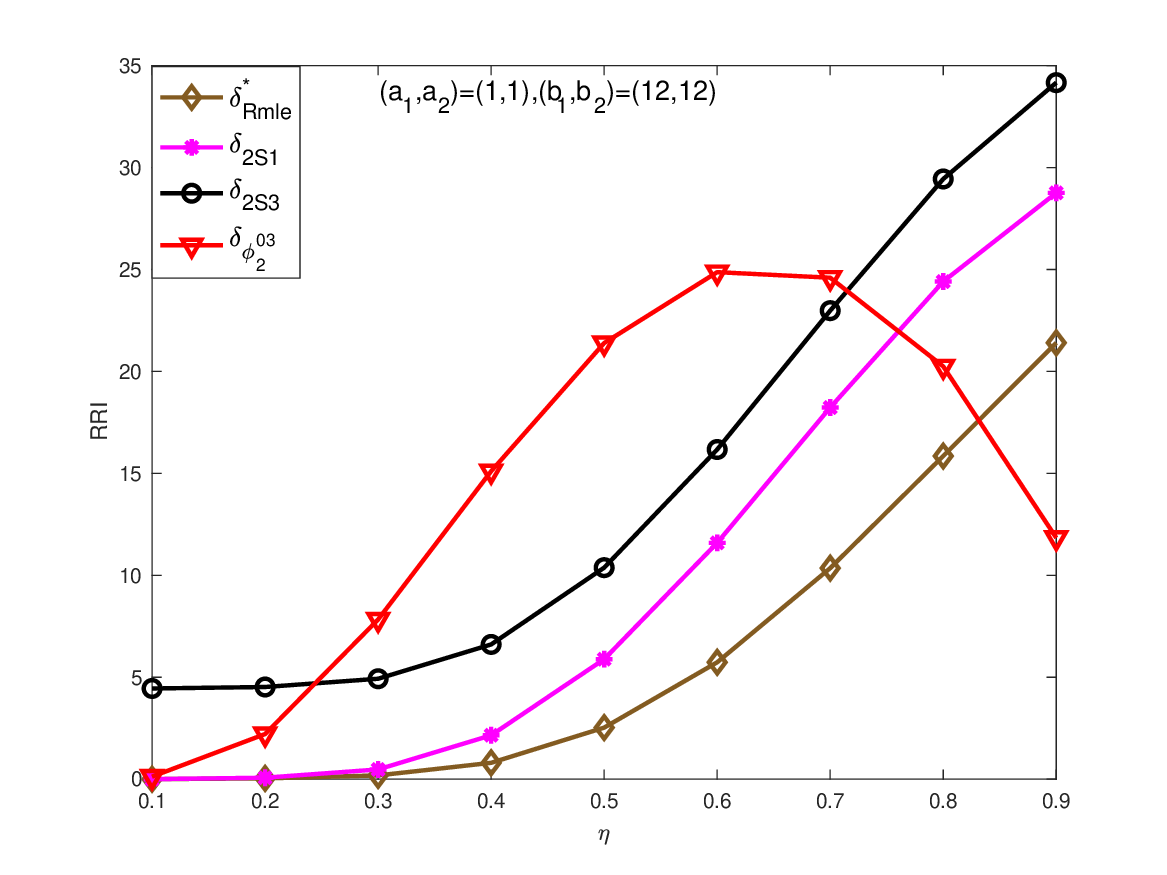}}
								\end{center}
								\caption{Relative risk improvement under symmetric loss function $L_3(t)$ for $\sigma_{2}$}\label{figs21L3}
							\end{figure}

\begin{figure}
							\begin{center}
								\subfigure[\tiny{$(n_1,n_2)=(10,9) ,(\mu_1,\mu_2)=(0.1,0.1) $}]{\includegraphics[height=5.5cm,width=8cm]{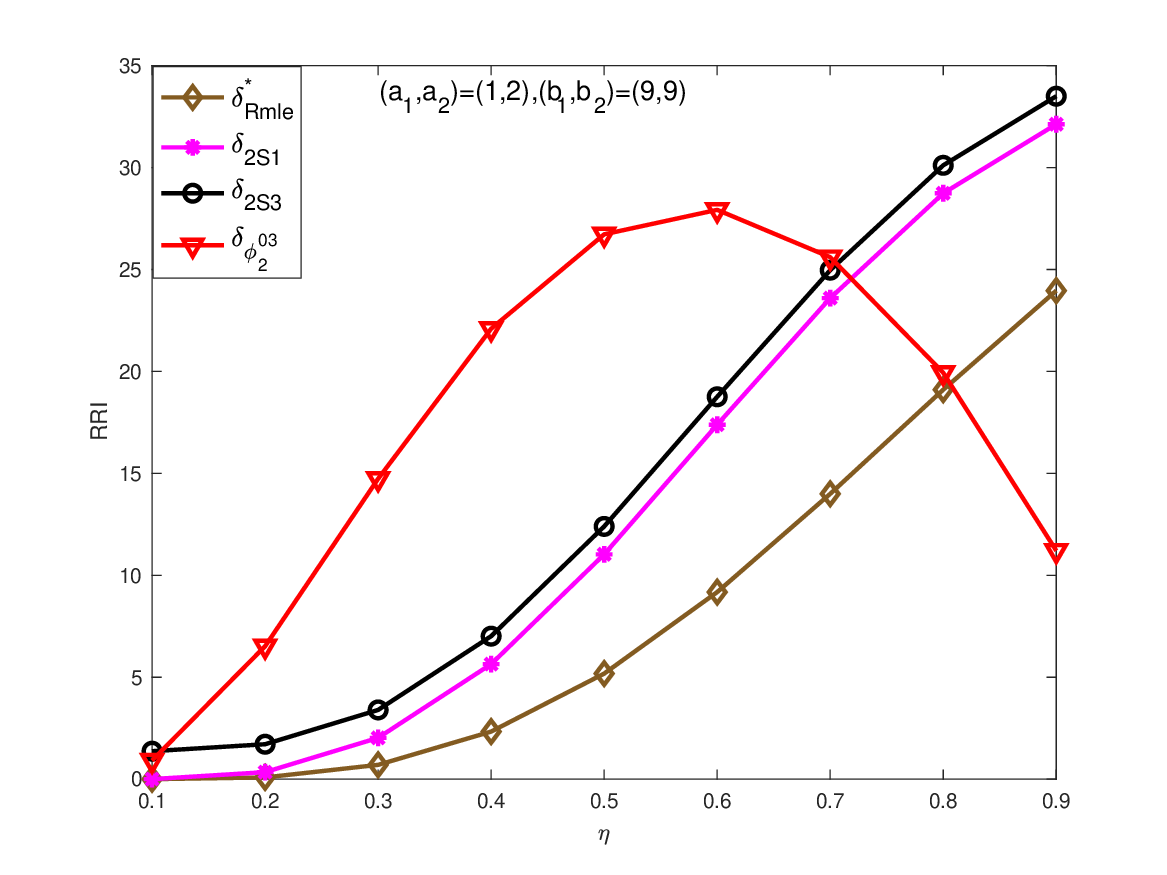}}
								\hspace{1cm}
								\subfigure[\tiny{$(n_1,n_2)=(10,9) ,(\mu_1,\mu_2)=(0.1,0.15) $}]{\includegraphics[height=5.5cm,width=8cm]{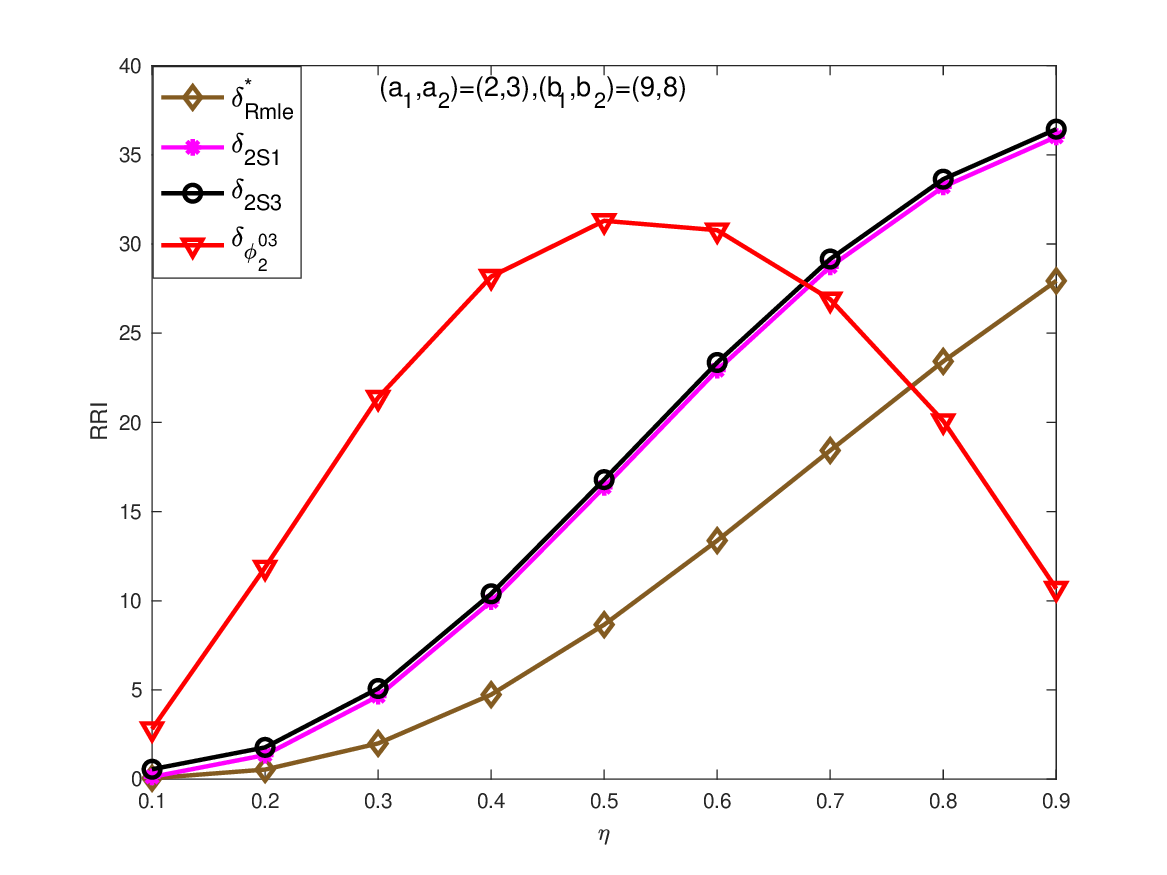}} 
								\hspace{1cm}
								\subfigure[\tiny{$(n_1,n_2)=(14,15) ,(\mu_1,\mu_2)=(0.4,0.7) $}]{\includegraphics[height=5.5cm,width=8cm]{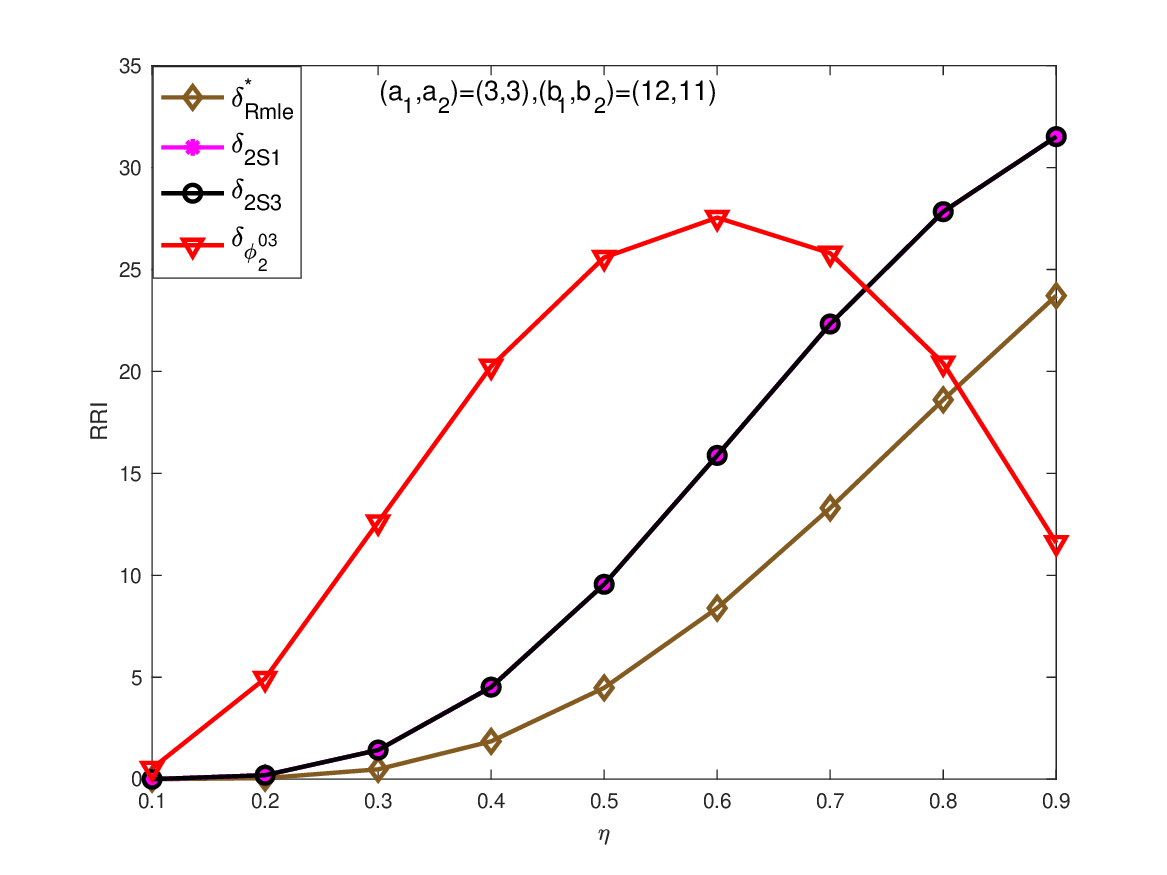}}
								\hspace{1cm}
								\subfigure[\tiny{$(n_1,n_2)=(14,15) ,(\mu_1,\mu_2)=(0.4,0.7)$}]{\includegraphics[height=5.5cm,width=8cm]{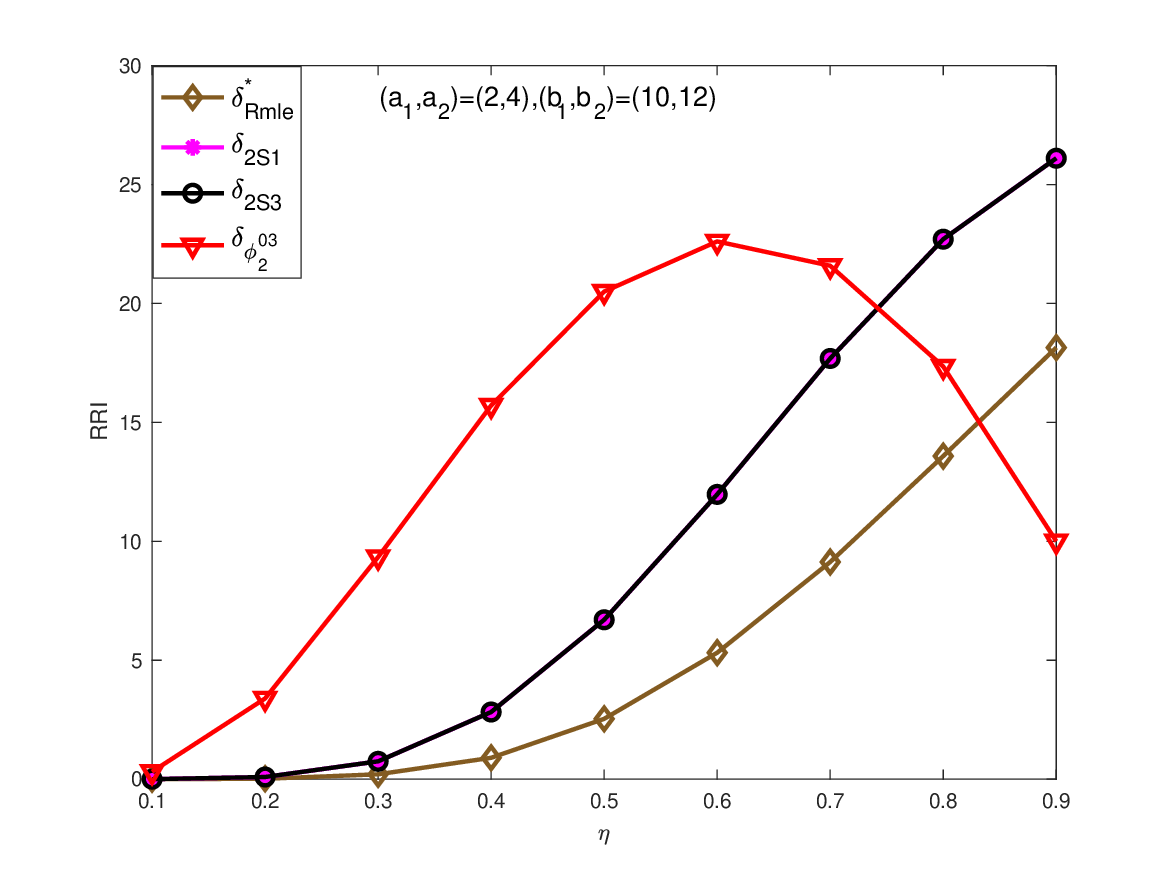}}
								\subfigure[\tiny{$(n_1,n_2)=(7,8) ,(\mu_1,\mu_2)=(-0.1,-0.2) $}]{\includegraphics[height=5.5cm,width=8cm]{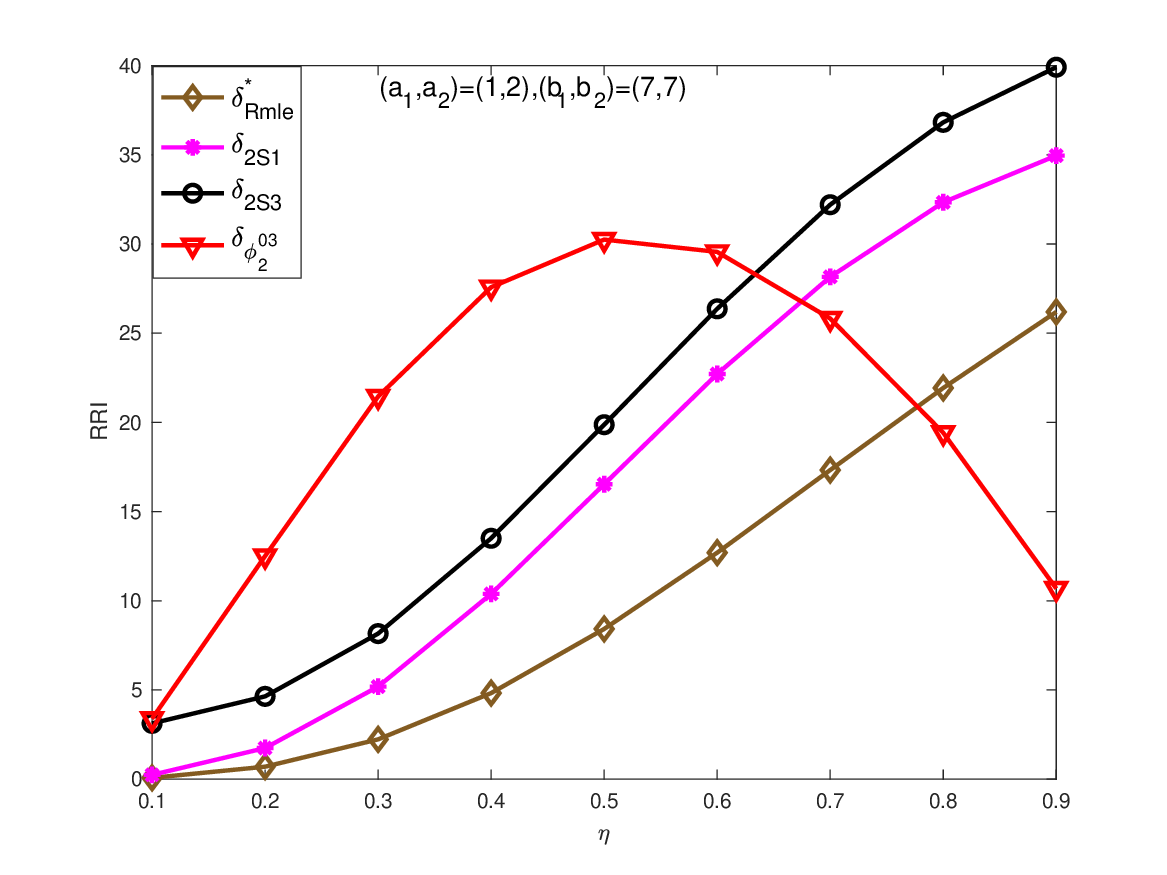}}
								\hspace{1cm}
								\subfigure[\tiny{$(n_1,n_2)=(7,8) ,(\mu_1,\mu_2)=(-0.1,-0.2) $}]{\includegraphics[height=5.5cm,width=8cm]{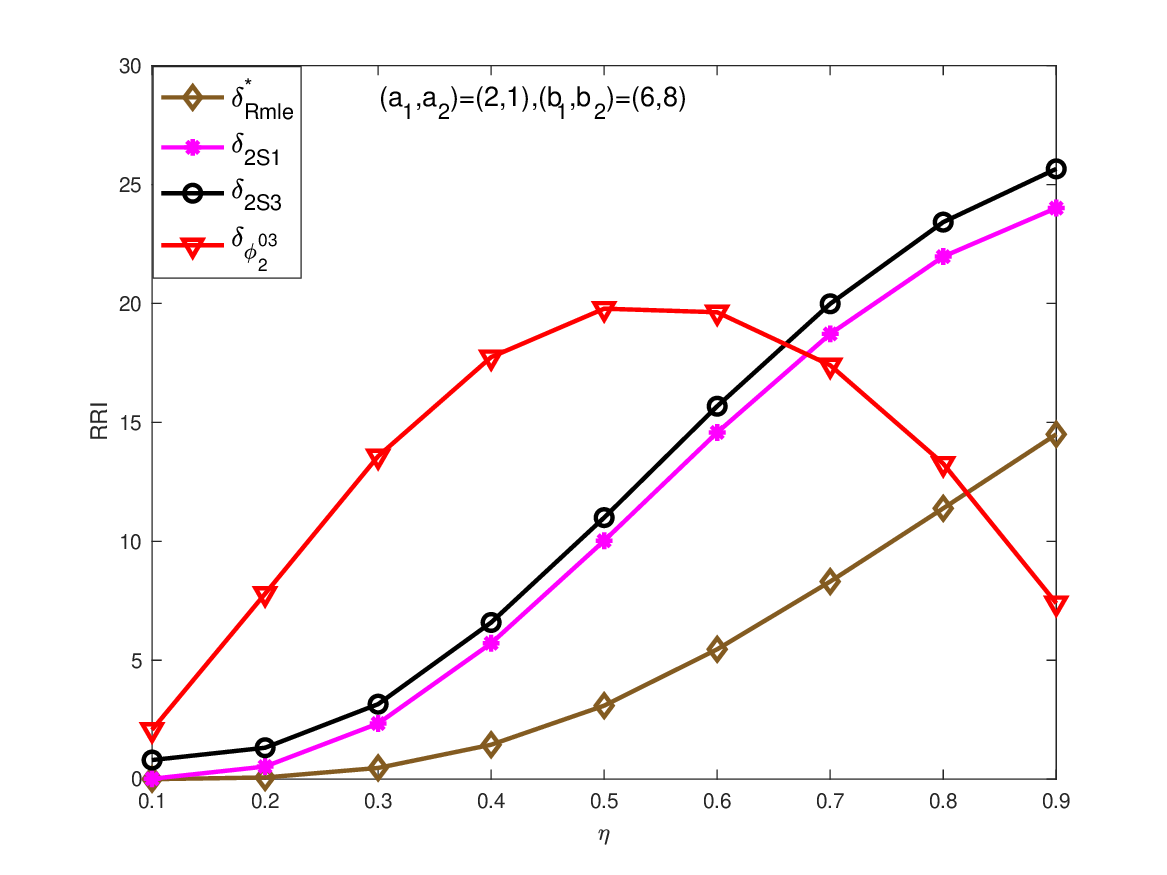}} 
							\end{center}
							\caption{Relative risk improvement under symmetric loss function $L_3(t)$ for $\sigma_{2}$}\label{figs22L3}
						\end{figure}

From the numerical comparison (Figure \ref{figs21L1} and \ref{figs22L1}) of proposed estimators of $\sigma_{2}$ we have derived following conclusion under quadratic loss.
 \begin{itemize}
 	\item [(i)] RRI of $\delta_{\phi_{2}^{01}}$ increases and then decreases as function of $\eta$.  The RRI of $\delta_{2S1},\delta_{2S3}$ and $\delta_{Rmle}^*$ is an increasing function of $\eta$.
 	\item[(ii)] The estimator $\delta_{\phi_{2}^{01}}$ performs better than $\delta_{2S1},\delta_{2S3}$  and $\delta_{Rmle}^*$ for small and moderate values of $\eta$ approximately 0.07. The region of improvement increases as sample size increases.
 	\item [(iii)] Risk performance of $\delta_{2S3}$ is better than $\delta_{2S1}$ and $\delta_{Rmle}^*$ for smaller values of  $\mu_1,\mu_{2}$.
 	 	\item [(iv)] The estimator $\delta_{2S1}$ and $\delta_{2S3}$ have equal risk performance for large and moderate values of $\mu_1$ and $\mu_{2}$, however for smaller values of $\mu_1,\mu_{2}$ The estimator $\delta_{2S3}$ performs better than $\delta_{2S1}$.
 	\item [(v)]The doubly shrinkage estimator denoted as $\delta_{2S3}$ always performs better than $\delta_{Rmle}^*$.
 	\item [(vi)] For small sample size and $\mu_1$, $\mu_2$ is close to $0$ all the estimators perform better than $\delta_{\phi_{2}^{01}}$, when $\eta$ close to $1$ that is approximately more than $0.8$.  
 \end{itemize} 
In the Figure \ref{figs21L2}, \ref{figs22l2} and \ref{figs23L2} we have plotted the RRI of the improved estimators with respect to entropy loss $L_2(t)$ and for symmetric loss $L_3(t)$ the RRI of improved estimators are plotted in the Figures \ref{figs21L3} and \ref{figs22L3}. We have similar type of observations for these type of loss functions too.

	
						\clearpage
	\subsection{Data analysis}
Here, we will present a real-life data analysis. For the data analysis, we consider the data of breaking strength of jute fiber of gauge length $20$mm and $5$ mm discussed in  \cite{xia2009study}. The data set contains 30 observations, by using one sample  Kolmogorov–Smirnov test we have verified that data following two parameter exponential distribution with location parameter $\mu_{1}=36.75$ and $\mu_2=129.08$ with p-values 0.6568 and 0.4908 respectively. Let $X_1\sim Exp(\mu_1,\sigma_1)$ and $X_2 \sim Exp(\mu_2,\sigma_2)$ denote the breaking strength of jute fiber of gauge length 20mm and 5mm, respectively, and we assume the that $\sigma_1 \le \sigma_2$. In the previous sections, we have proposed several improved estimators for $\sigma_1$ and $\sigma_2$. The values of the improved estimators are tabulated in the below tables for different values of $a_i$ and $b_i$, $i=1,2$. 
\begin{table}[h!]
	\centering
	\caption{gauge length $20$ mm }
	\vspace{0.2cm}
	\label{tabledata1}
	\begin{tabular}{llllllll}
		\hline
	71.46 & 419.02 & 284.64& 585.57 & 456.60 & 113.85&187.85& 688.16 \\ 
	 662.66 & 45.58 & 578.62 & 756.70& 594.29 & 166.49 &99.72 & 707.36 \\
	  765.14& 187.13 & 145.96 & 350.70 & 547.44 & 116.99 &375.81  &581.60 \\
	   119.86 & 48.01 & 200.16 & 36.75 & 244.53 &83.55      \\ \hline
	\end{tabular}
\end{table}
\begin{table}[h!]
	\centering
	\caption{gauge length $5$ mm }
	\vspace{0.2cm}
	\label{tabledata2}
	\begin{tabular}{llllllll}
		\hline
		566.31 & 270.79 &516.28 & 823.03 & 226.53 & 367.70   & 441.87 &618.57\\
		 546.11 & 268.20 & 315.33 & 809.23  & 218.86 & 583.97 & 304.84 &129.08 \\
		  537.45 &496.28 & 167.87 & 306.99 & 178.25 & 370.02 & 168.20 &554.61 \\
		   360.80  & 260.97  & 254.29 & 295.51  &187.68& &   \\     \hline
	\end{tabular}
\end{table}
	
	\begin{table}[h!]\label{da1}
		\centering
		\caption{Values of estimators of $\sigma_1$ under quadratic loss function $L_1(t)$}
		\begin{tabular}{ccccccccc}
			\hline\hline
			$(a_1,a_2)$ & $(b_1,b_2)$ & $\delta_{01}^1$ & $\delta_{Rmle}$ & $\delta_{1S1}$ & $\delta_{1S2}$ & $\delta_{1S3}$ & $\delta_{\phi_1^{01}}$ & $\delta_{\phi_{1.5,1}}$ \\ \hline\hline
			(1,1)       & (30,30)     & 303.99  & 279.64        & 284.38         & 298.01        & 303.99         & 264.68               & 352.90               \\ 
			(2,3)       & (27,28)     & 334.57         & 285.16         & 290.75         & 310.58         & 334.57        & 274.79               & 354.99             \\ 
			(1,1)       & (29,27)     & 314.18         & 290.56         & 295.84        & 310.24         & 314.18         & 278.76               & 360.81               \\ 
			(4,2)       & (30,30)     & 302.31          & 262.43         & 267.21         & 296.89        & 302.31         & 247.64                 & 333.06              \\ \hline\hline
		\end{tabular}
	\end{table}
	\begin{table}[h!]\label{da2}
		\centering
		\caption{Values of estimators of $\sigma_1$ under entropy loss function $L_2(t)$}
		\begin{tabular}{ccccccccc}
			\hline \hline
			$(a_1,a_2)$ & $(b_1,b_2)$ & $\delta_{01}^2$ & $\delta_{Rmle}$ & $\delta_{1S1}$ & $\delta_{1S2}$ & $\delta_{1S3}$ & $\delta_{\phi_1^{01}}$ & $\delta_{\phi_{1.5,2}}$ \\ \hline \hline
			(1,1)       & (30,30)     & 314.47 & 279.64        & 289.28         & 264.10         & 314.47         & 270.03                 & 360.61                 \\ 
			(2,3)       & (27,28)     & 347.96          & 285.16        & 296.57         & 245.18        & 347.96          & 280.85                & 363.16                 \\ 
			(1,1)       & (29,27)     & 325.40        &290.56          & 301.32        & 255.71         & 325.40          & 284.56                & 376.36                 \\ 
			(4,2)       & (30,30)     & 313.94         & 262.44        & 272.16         & 251.30         & 313.94          & 253.02                 & 333.57                 \\ \hline\hline
		\end{tabular}
	\end{table}
	\begin{table}[h!]\label{da3}
		\centering
		\caption{Values of estimators of $\sigma_1$ under symmetric loss function $L_3(t)$}
		\vspace{0.2cm}
		\begin{tabular}{ccccccccc}
			\hline\hline
			$(a_1,a_2)$ & $(b_1,b_2)$ & $\delta_{01}^3$ & $\delta_{Rmle}$ & $\delta_{1S1}$ & $\delta_{1S2}$ & $\delta_{1S3}$ & $\delta_{\phi_1^{01}}$ & $\delta_{\phi_{1.5,3}}$ \\ \hline\hline
			(1,1)       & (30,30)     &320.04 &279.64        &291.81        & 305.67        & 320.04         & 272.76              & 364.53               \\ 
			(2,3)       & (27,28)     &355.13     & 285.16         & 299.58         & 319.82         & 355.13        &  283.97         & 367.36                         \\ 
			(1,1)       & (29,27)     & 331.37      &290.56          & 304.15         & 318.80        & 331.37         & 287.54         &380.47                            \\ 
			(4.2)       & (30,30)     & 320.16      &262.44          & 274.71         & 305.07         & 320.16        & 255.78         & 337.49                           \\ \hline\hline
		\end{tabular}
	\end{table}
	
	\begin{table}[h!]
		\centering
		\caption{Values of estimators of $\sigma_2$ under quadratic $L_1(t)$ for $\sigma_2$}
		\vspace{0.1cm}
		\begin{tabular}{ccccccccc}
			\hline\hline
			$(a_1,a_2)$ & $(b_1,b_2)$ & $\delta_{02}^1$ & $\delta_{Rmle}^*$ & $\delta_{2S1}$ & $\delta_{2S2}$ & $\delta_{2S3}$ & $\delta_{\phi_{2}^{01}}$ & $\delta^*_{\phi_{1.5,1}}$  \\ \hline\hline
			$(1,1)$     & $(30,30)$    & 255.29         & 279.63           & 284.37         & 255.29         & 284.37         & 309.01                 & 460.22                        \\ 
			$(2,3)$     & $(27,28)$    & 235.75         & 285.15           & 290.74         & 235.75         & 290.74         & 314.15                  & 490.19                        \\ 
			$(1,1)$     & $(29,27)$     & 265.18         & 290.54           & 295.82        & 265.18          & 295.82         & 323.85                   & 491.51                         \\ 
			$(4,2)$     & $(30,30)$    & 225.31         & 262.42           & 267.19        & 225.31        & 267.19         & 287.84                 &  433.32                         \\ \hline\hline
		\end{tabular}
	\end{table}
	\begin{table}[h!]
		\centering
		\caption{Values of estimators of $\sigma_2$ under entropy loss $L_2(t)$ for $\sigma_2$}
		\vspace{0.2cm}
		\begin{tabular}{ccccccccc}
			\hline\hline
			$(a_1,a_2)$ & $(b_1,b_2)$ & $\delta_{02}^2$ & $\delta_{Rmle}^*$ & $\delta_{2S1}$ & $\delta_{2S2}$ & $\delta_{2S3}$ & $\delta_{\phi_{2}^{02}}$ & $\delta^*_{\phi_{1.5,2}}$ \\ \hline\hline
			$(1,1)$     & $(30,30)$    & 264.10       & 279.63         &289.27      & 264.10       & 289.27        & 316.12                & 471.33                       \\ 
			$(2,3)$     & $(27,28)$    & 245.18         & 285.15          & 296.55         & 245.18        & 296.55        &322.15                  & 503.01                         \\ \
			$(1,1)$     & $(29,27)$     & 275.38          & 290.54           &301.30        & 275.38          & 301.30         &332.07                   &504.64                        \\ 
			$(4,2)$     & $(30,30)$    & 233.35           &262.42         & 272.14        & 233.35         & 272.14         &294.59                  & 443.79                         \\ \hline\hline
		\end{tabular}
	\end{table}
	\begin{table}[h!]
		\centering
		\caption{Values of estimators of $\sigma_2$ under symmetric loss $L_3(t)$ for $\sigma_2$}
		\vspace{0.2cm}
		\begin{tabular}{ccccccccc}
			\hline\hline
			$(a_1,a_2)$ & $(b_1,b_2)$ & $\delta_{02}^3$ & $\delta_{Rmle}^*$ & $\delta_{2S1}$ & $\delta_{2S2}$ & $\delta_{2S3}$ & $\delta_{\phi_{2}^{03}}$ & $\delta^*_{\phi_{1.5,3}}$ \\ \hline\hline
			$(1,1)$     & $(30,30)$    & 268.77        & 279.63          & 291.80         & 268.77         & 298.80        & 319.88                  & 477.27                         \\ 
			$(2,3)$     & $(27,28)$    & 250.24         & 285.15         & 299.56        & 250.24      & 299.56         & 326.40                  & 509.85                        \\ 
			$(1,1)$     & $(29,27)$     & 280.84          & 290.54          & 304.13         &280.84          & 304.13        & 336.45                  & 511.73                         \\ 
			$(4,2)$     & $(30,30)$    & 237.64         & 262.42           & 272.70         & 237.64        & 274.70         & 298.16                  & 449.35                         \\ \hline\hline
		\end{tabular}
	\end{table}					
 \section{Conclusions}
We have population $\Pi_i$ having distribution $Exp(\mu_i,\sigma_i), i=1,2$ and $\sigma_1 \le \sigma_2$. In this paper, we have studied the estimation $\sigma_i$ based on a doubly type-II censored sample when location parameters $\mu_1$ and $\mu_2$ are unknown. This estimation problem has been studied with respect to a general scale invariant loss function. We have proposed several improved estimators using the techniques of \cite{stein1964}. We also obtained a class of improved estimators using the IERD approach of \cite{kubokawa1994unified}. We showed that the boundary estimator of the proposed class is a generalized Bayes estimator. We have obtained the improved estimators for three special loss functions, namely quadratic loss, entropy loss, and symmetric loss function. Maruyama type estimators are obtained for $\sigma_1$ and $\sigma_2$. For $\sigma_1$, we have derived a \cite{strawderman1974minimax} type improved estimator with respect to quadratic and entropy loss functions. A double shrinkage estimator is obtained for $\sigma_2$. We have seen that from these results, we get the results for i.i.d. sampling scheme, record values, type-II censoring scheme, and progressive type-II censoring scheme. So, our work provides a unified treatment for different important sampling schemes. We conduct a simulation study to compare the risk performance of the improved estimators. Finally, we gave a data analysis for implementation purposes. 
\section*{Acknowledgment}
\noindent Lakshmi Kanta Patra thanks the National Board of Higher Mathematics, India for providing financial support to carry out this research  with project number 02011/38/2023 NBHM (R.P)/R$\&$DII/13409. 
\bibliography{censor_exp}

\end{document}